\documentclass[12pt,a4paper]{article}
 
\usepackage{graphicx}
\usepackage[shortlabels]{enumitem} %%%changed on 22-09-2020 
\usepackage{ulem}
\usepackage{amssymb,amsmath}

\usepackage{hyperref}

\hypersetup{
    colorlinks=true,
    linkcolor=blue,
    filecolor=magenta,      
    urlcolor=cyan,
    }
\newcommand{\ignore}[1]{}  %% {to remove the content and {#1} to include them}
\def\one{\mbox{1\hspace{-4.25pt}\fontsize{12}{14.4}\selectfont\textrm{1}}} % 11pt %% to define indicator function 1
\usepackage[shortlabels]{enumitem} %%%changed on 22-09-2020 
\usepackage{ulem}
\usepackage{amssymb}

\newcommand{\RSCost}[3]{J^{#1}_{{#2},{#3}}(x)} %% new one for each cost function..
\newcommand{\RSCostVec}[3]{J^{#1}_{{#2},{#3}}} %% new one for each cost function..
\newcommand{\LDCost}[2]{v^{#1}_{#2}(x)} %% new one for each cost function..
\newcommand{\LDCostVec}[2]{v^{#1}_{#2}} %% new one for each cost function..
\newcommand{\RSCostT}[4]{J^{#1}_{{#2},{#3},{#4}}(x)} %% new one for each cost function..
\newcommand{\RSCostTVec}[4]{J^{#1}_{{#2},{#3},{#4}}} %% new one for each cost function..
\newcommand{\LDCostT}[3]{v^{#1}_{{#2},{#3}}(x)} %% new one for each cost function..
\newcommand{\LDCostTVec}[3]{v^{#1}_{{#2},{#3}}} %% new one for each cost function..
\newcommand{\FeasibleSetOrgnal}{\mathcal{F}}
\newcommand{\FeasibleSetLDTruncBlwAti}[2]{\mathcal{C}^{-}_{{#1},{#2}}}
\newcommand{\FeasibleSetLDTruncBlw}[1]{{\mathcal{C}^{-}_{{#1}}}}
\newcommand{\FeasibleSetRSTruncBlw}[1]{\hat{\mathcal{C}}^{-}_{{#1}}}
\newcommand{\FeasibleSetLDTruncAbv}[1]{\mathcal{C}^{+}_{{#1}}}
\newcommand{\FeasibleSetRSTruncAbv}[1]{\hat{\mathcal{C}}^{+}_{{#1}}}
\newcommand{\FeasibleSetLDBlwAti}[1]{\mathcal{C}_{#1}}
\newcommand{\FeasibleSetRSTruncBlwAtj}[1]{\hat{\mathcal{C}}^{-}_{j,{#1}}}
\newcommand{\FeasibleSetRSBlwAtj}{\hat{\mathcal{C}}_{j}}
\newcommand{\FeasibleSetLDTruncAbvAti}[1]{{\mathcal{C}^{+}_{i,{#1}}}}
\newcommand{\FeasibleSetRSTruncAbvAtj}[1]{\hat{\mathcal{C}}^{+}_{j,{#1}}}
\newcommand{\FeasibleSetTrunBlw}[1]{\mathcal{F}^{-}_{{#1}}}
\newcommand{\FeasibleSetTrunAbv}[1]{\mathcal{F}^{+}_{{#1}}}
\newcommand{\ForVT}[2]{\ifthenelse{\(2<1\)}{#1}{#2}}
\newcommand{\Qfunc}[3]{Q_{#1,#2}^{#3}}
\newcommand{\CRSMDPFinHorOne}{\hyperlink{theequation}{(P^-_T)}}
\newcommand{\CRSMDPFinHorTwo}{\hyperlink{theequation}{(P^+_T)}}
\newtheorem{thm}{Theorem}[section]

\newtheorem{rem}[thm]{Remark}
\newtheorem{ex}{Example}[section]

\newtheorem{lem}[thm]{Lemma}
\newtheorem{prop}[thm]{Proposition}
\newtheorem{defn}[thm]{Definition}
\newcommand{\E}{\mathbb{E}}
\newcommand{\R}{\mathbb{R}}
\newcommand{\Z}{\mathbb{Z}}

\DeclareMathOperator{\cl}{cl}
\DeclareMathOperator{\FeasibleSetLD}{\mathcal{C}}
\newcommand{\FeasibleSetRS}{\hat{\mathcal{C}}}
\newcommand{\FeasibleSetLDFinHor}{\bar{\mathcal{C}}}
\newcommand{\FeasibleSetRSFinHor}{\check{\mathcal{C}}}
\newcommand{\eop}{\hfill{$\Box$}}
\numberwithin{equation}{section}
\begin{document}

%\frontmatter
\thispagestyle{empty}
 \title{Approximate Solutions To Constrained Risk-Sensitive Markov Decision Processes}
 \author{$^1$Uday Kumar M, $^1$Sanjay P Bhat,\\
 $^2$Veeraruna Kavitha, $^2$Nandyala Hemachandra \\
 $^1$TCS Research, Hyderabad, India  and \\
 $^2$IEOR, IIT Bombay}

 \maketitle
\begin{abstract}
This paper considers the problem of finding near-optimal Markovian randomized (MR) policies for finite-state-action, infinite-horizon, constrained risk-sensitive Markov decision processes (CRSMDPs). Constraints are in the form of standard expected discounted cost functions as well as expected risk-sensitive discounted cost functions over finite and infinite horizons. The main contribution is to show that the problem possesses a solution if it is feasible, and to provide two methods for finding an approximate solution in the form of an ultimately stationary (US) MR policy. The latter is achieved through two approximating finite-horizon CRSMDPs which are constructed from the original CRSMDP by time-truncating the original objective and constraint cost functions, and suitably perturbing the constraint upper bounds. The first approximation gives a US policy which is $\epsilon$-optimal and feasible  for the original problem,  while the second approximation gives a near-optimal  US policy whose violation of the original constraints is bounded above by a specified $\epsilon$. A key step in the proofs is an appropriate choice of a metric that makes the set of infinite-horizon MR policies and the feasible regions of the three CRSMDPs compact, and the objective and constraint functions continuous.  A linear-programming-based formulation for solving the approximating finite-horizon CRSMDPs is also given.   

\end{abstract}

% section wise tex file

\section{Introduction}
\label{sec_intro}
A Markov decision process (MDP) evolving  over a finite set of states at discrete decision epochs under the influence of a finite number of actions is specified in terms of immediate rewards (or costs) and  controlled state transition probabilities satisfying the Markov property.

 A policy is a sequence of decisions rules, one for each decision epoch. A policy is \textbf{stationary} if it applies the same decisions rule at all the epochs,  and \textbf{ultimately stationary} (US) if it applies the same decision rule beyond a certain epoch.  
The decision problem involving an MDP is to choose the policy that optimizes the cost obtained by appropriately aggregating immediate costs incurred over several decision epochs. Examples of aggregate costs include total, average or discounted costs over a horizon that may be finite or infinite. In standard MDPs, one seeks to optimize the expected aggregate cost. However, this results in a risk-neutral approach which does not take into account the risk preferences of possibly  risk-sensitive decision makers.

The literature on MDPs considers various means for incorporating risk sensitivity into the decision problem. These include penalizing the long run variance of the costs discussed in \cite{Flair_Kallenberg_Lee}, \cite{Abhijit} and \cite{Li_Xia}
attempting   mean-variance tradeoff by simultaneously imposing thresholds on the  expectation and variance of the aggregate cost discussed in    \cite{Mannor_Tsitsiklis},  maximizing the probability of maintaining the aggregate cost within a budget is considered by authors in \cite{Moreira_Delgado_Barros}, using the conditional value-at-risk of the aggregate cost either as a constraint is discussed in \cite{Borkar_Rahul} or as the objective is discussed in \cite{Yang}, and enforcing a threshold on the probability that the system enters a high-risk state is discussed in \cite{Geibel_Wysotzki}. The most widely followed approach for taking risk sensitivity into account is to consider optimization of the expected exponential utility of the aggregate cost or, equivalently, its certainty equivalent dealt in \cite{Nicole_Ulrich,HwrdMathson}.  We follow the same approach in this paper, and use the term risk-sensitive MDP (RSMDP) to refer to an MDP with expected exponential utility of the total discounted cost as the risk-sensitive (RS) cost to be optimized.  For the sake of clarity, we will use the terminology ``standard discounted cost'' to refer to the expectation of the discounted sum of immediate costs over the horizon of interest.

Work in \cite{SJaq1973,SJaq1976}  contain two of the earliest treatments of infinite-horizon RSMDPs. Authors in \cite{SJaq1973} examined the relationship between optimality in terms of  expected exponential utility  and  optimality in terms of higher moments of the discounted cost for infinite-horizon RSMDPs. The reference \cite{SJaq1976} showed that while optimal policies in infinite-horizon RSMDPS are not always stationary, US optimal policies do exist under certain assumptions. Work in \cite{Kavitha_NH_Atul} gave a solution to finite-horizon unconstrained RSMDPs with finite state and action spaces using dynamic and linear programming (LP).

Emphasizing ease of computation and  practical implementation over exact optimality in \cite{UdaySanjayKavitaNH} focused on $\epsilon$-optimal US policies for infinite horizon RSMDPS. The results of \cite{UdaySanjayKavitaNH} exploited two key ideas for obtaining $\epsilon$-optimal policies. The first idea, which was also considered earlier in \cite{Coraluppi_1997_PhD} and \cite{Coraluppi_Marcus_1999}, that involves truncating the tail cost beyond a finite horizon and using an arbitrarily chosen stationary policy thereafter. Discounting ensures that the resulting truncation error is small. Consequently, optimizing the truncated RS cost results in an $\epsilon$-optimal policy referred in \cite{UdaySanjayKavitaNH} as the ultimately stationary tail off (USTO) policy. The second idea is based on the fact that the certainty equivalent of the RS cost approaches the standard discounted cost as the risk factor $\gamma$ approaches $0$ as discussed by authors in \cite{SJaq1973,SJaq1976} and \cite[Thm. 3]{UdaySanjayKavitaNH}.  Authors in \cite{UdaySanjayKavitaNH} makes subtle use of this fact by observing that the effective risk factor for the immediate costs sufficiently far out into the tail can be taken to be arbitrarily small due to the effect of discounting. The RS tail cost may therefore be approximated by the exponential of the standard discounted cost of the tail of the policy.  These ideas lead to  an $\epsilon$-optimal US policy called ultimately stationary linear discounted (USLD), which is obtained by appending the stationary optimal policy for the infinite-horizon discounted cost to the optimal policy for a finite-horizon terminal-cost  RSMDP whose  terminal cost is determined by the optimal infinite-horizon standard discounted cost.   

Real-life optimization problems invariably  involve constraints, and MDPs are no exceptions. Previous applications of constrained MDPs in risk-neutral as well as risk-sensitive settings include pavement management systems in \cite{Kamal_Ram_George}, multi-arm bandits in \cite{Eric_Park_Uriel} and delay torrent networks in \cite{Kavitha_NH_Atul}. 
The objective of this paper, therefore, is to extend the ideas of \cite{UdaySanjayKavitaNH} to constrained RSMDPS (CRSMDPs) that involve constraints based on standard discounted and RS costs over finite and infinite horizons. In particular, this paper seeks to extend the USTO idea in \cite{UdaySanjayKavitaNH} to infinite-horizon CRSMDPs.

While the literature on constrained MDPs (CMDPs) is not as extensive as that on unconstrained MDPs, constrained MDPs have been fairly well studied in the risk-neutral context. In one of the earliest developments, \cite{Derman_Klein} gave a LP-based method to compute the optimal policy for a finite-horizon, total cost CMDP. \cite{Kallenberg} gave a more detailed exposition of the same.  The influential paper \cite{Feinberg_Shwartz} showed that a finite-state infinite-horizon CMDP having both objective and constraints of the standard discounted type possesses stationary randomized optimal policies as well as ultimately deterministic US optimal policies, and gave a LP-based algorithm and an iterative procedure to compute the latter. \cite{Krishna_NH} analyzed the sensitivity of the optimal cost for the CMDP given in \cite{Feinberg_Shwartz}, and also gave a lower bound on the horizon beyond which the optimal policy is stationary and deterministic. The reader is referred to book \cite{Eitan_CMDP} for additional background and references on CMDPs in the risk-neutral case.

In contrast to CMDPs in the risk-neutral setting, CRSMDPs seem to have received very little attention in the literature. \cite{Kavitha_NH_Atul} provide a LP-based solution to finite-horizon CRSMDPs involving only finite-horizon standard discounted costs as constraints.  \cite{Haskell_Jain} consider infinite-horizon risk constrained MDPs involving the  minimization of a general law-invariant risk measure applied to infinite-horizon discounted cost. They also consider minimization of the expected utility of infinite-horizon discounted cost subject to a single constraint which is either a stochastic dominance constraint or a chance constraint. However,  
to the best of our knowledge, infinite-horizon CRSMDPs with multiple constraints of risk-neutral and risk-sensitive type over finite and infinite horizons have not been considered before. This motivates us to consider the problem of minimizing the expected exponential utility of total discounted cost over an infinite horizon scaled by a risk-factor  under a fairly general set of finite- and infinite-horizon RS and standard discounted cost constraints. Our work thus seeks to  extend that of \cite{SJaq1976} and  \cite{UdaySanjayKavitaNH} by including constraints. It also extends the work reported in \cite{Kavitha_NH_Atul} by allowing the horizon to be infinite and including RS cost constraints in addition to standard discounted cost constraints. 

The main contribution of the current paper is to show that the problem possesses a solution if it is feasible, and provide two approximation techniques for finding a near-optimal US policy for the problem. The  approximation is through two finite-horizon CRSMDPs which are constructed from the original CRSMDP by time-truncating the original objective and infinite-horizon constraint cost functions, and suitably perturbing the bounding constants defining the constraints. The construction of the approximating CRSMDPs is such that the feasible region of the first (second) approximating CRSMDPs provides an inner (outer) approximation to the feasible region of the original CRSMDP, while the optimal values of both converge to the optimal value of the original CRSMDP as the truncation horizon is extended. Consequently, a solution of the first approximating CRSMDP provides an $\epsilon$-optimal US solution to the original CRSMDP under an additional assumption. On the other hand, a solution to the second approximating CRSMDP provides a near-optimal US policy which is only guaranteed to be  $\epsilon$-feasible for the original problem in the sense that constraint violations, if any, do not exceed $\epsilon$. A key step in the development is the  introduction of a metric in which the set of MR policies is compact, and the objective and constraint functions are continuous. For the sake of completeness, we also provide a LP formulation for computing the solutions to the two approximating CRSMDPs.

The paper is organized as follows.  The mathematical framework and the formal problem statement are  given in Section \ref{sec_model_framework}. The  Section \ref{sec_main_results} gives the main results  on the existence of solutions to the infinite-horizon CRSMDP,  and the $\epsilon$-optimality and $\epsilon$-feasibility of the solutions obtained from the two finite-horizon approximating CRSMDPs described above.  The proofs of the main results along with related preliminaries are given in Section \ref{sec_proof_mainres}, while the LP formulation is described in Section  \ref{sec_algo}. The proofs of all subsidiary results are given in the appendices. 
\section{Model Framework}\label{sec_model_framework}
As in the case of an MDP with standard discounted cost, the description of a RSMDP involves  a state space, action space, immediate costs or rewards, and controlled transition probabilities. However, unlike the former, an RSMDP involves optimizing the expected exponential utility of the aggregated cost built up from costs collected over several decision epochs. In this paper, the aggregated cost is taken as the discounted sum of costs.
  
Let $\mathcal{S}=\{s_1,s_2,\ldots,s_m\}$ and $\mathcal{A}=\{a_1,a_2,\ldots,a_n\}$ denote the sets of all possible states and actions, respectively,  with cardinalities $|\mathcal{S}|=m$ and $|\mathcal{A}|=n$. The controlled transition probability, denoted by $p(s_j|s_i,a_k)$, represents the conditional probability that the system transitions to state $s_j$ given that action $a_{k}$ is taken in the current state $s_{i}$. Observe that the controlled transition probabilities are time homogeneous. Let $R: \mathcal{S}\times\mathcal{A} \rightarrow \R$ denote the immediate cost function. 
%We assume that $R$ is  absolutely bounded by $C>0$, that is, $|R(s,a)|\leq C$ for all $(s,a)\in \mathcal{S}\times\mathcal{A}$.  

Decisions are made at discrete epochs. In general, one can consider history-dependent decisions, in which the decision at any epoch requires the history of the process up to that epoch. However, in this paper, we only consider \textbf{Markovian} decisions, in which the decision at any epoch is based only on the current state at that epoch. More precisely, a \textbf{Markovian randomized (MR)} decision rule is a mapping $d:\mathcal{S}\rightarrow \mathcal{P}(\mathcal{A})$ from the state space $\mathcal{S}$ to the space $\mathcal{P}(\mathcal{A})$ of probability distributions on the action space. Here, $d(s)$ is the conditional probability mass  function of the  action chosen given that the state at the decision epoch is $s$. Note that decision rules which deterministically assign a unique action to each state are special cases of  MR decision rules. An MR policy is an infinite sequence of MR decision rules, one for each decision epoch. We will denote the set of MR policies by $\Pi_{\rm MR}$, and  use $\pi=\{d_{t}\}_{t=0}^{\infty}$ to denote an arbitrary policy $\pi\in\Pi_{\rm MR}$, with $d_{t}$ representing the decision rule applied at epoch $t$. An MR policy $\pi=\{d_{t}\}_{t=0}^{\infty}$ is \textbf{ultimately stationary} (US) if there exists $T\geq 0$ such that $d_{i}=d_{j}$ for all $i,j\geq T$, and \textbf{stationary} if $d_{i}=d_{j}$ for all $i,j\geq 0$.  

For every initial state $x\in\mathcal{S}$, each policy $\pi$ induces a probability measure on the set of sequences of state-action pairs $\Omega= \mathcal{(S\times A)\times(S\times A)\times\cdots }$  (see, for example, Chapter 2  of \cite{PUT}). We denote the state and action at a time $t\in\{0,1,\ldots\}$  by random variables $X_{t}$ and $A_{t}$, respectively, on $\Omega$. 
% and note that $X_{t}$ and $A_{t}$ are random variables . 
Let $E^\pi_x [\cdot]$ represent the expectation conditioned on the initial state $X_0=x$ under the probability measure induced by $\pi$. 
Throughout this paper, we fix a discount factor $\beta\in(0,1)$ and a risk factor  $\gamma \in \R\setminus\{0\}$. Unlike \cite{SJaq1973,SJaq1976}, we allow $\gamma$ to take both negative and positive values. 
%\subsection*{Notations and Definitions}
%\begin{definition}
\begin{defn}
\textbf{[Standard discounted and RS cost functions:]}  Given a policy $\pi\in\Pi_{\rm MR}$, an immediate cost function $R$ and an initial condition $x$, the RS discounted cost and the standard discounted cost over an infinite horizon as well as over a finite horizon $T\geq 1$ are given by 
\begin{eqnarray}
\RSCost{\pi}{\gamma}{R} & := \E^\pi_x\left[e^{ \gamma\sum\limits_{t=0}^{\infty}\beta^t R(X_t,A_t) }\right] \mbox{, }     \RSCostT{\pi}{\gamma}{R}{T} &:= \E^\pi_x\left[e^{ \gamma\sum\limits_{t=0}^{T-1}\beta^t R(X_t,A_t) }\right], \label{eq_RSCost_FinandinfHor}\\
\LDCost{\pi}{R} &:= \E^\pi_x\left[\sum\limits_{t=0}^{\infty}\beta^t R(X_t,A_t)\right]
\mbox{, }       
\LDCostT{\pi}{R}{T} &:= \E^\pi_x\left[\sum\limits_{t=0}^{T-1}\beta^t R(X_t,A_t)\right].\label{eq_LDCost_FinandinfHor}
\end{eqnarray}\label{costdef}
\end{defn}
%Define  $\RSCostVec{\pi}{\gamma}{R} := (\RSCost{\pi}{\gamma}{R})_{x \in {\cal X} }$,  $\RSCostTVec{\pi}{\gamma}{R}{T}:= (\RSCostT{\pi}{\gamma}{R}{T})_{x \in {\cal X}}$, $\LDCostVec{\pi}{R}:=(\LDCost{\pi}{R})_{x \in {\cal X} }$, and $\LDCostTVec{\pi}{R}{T}:=(\LDCostT{\pi}{R}{T})_{x \in {\cal X} }$ the vector of corresponding costs over all initial states $x\in\mathcal{S}$.

Let  $\RSCostVec{\pi}{\gamma}{R}$, $\RSCostTVec{\pi}{\gamma}{R}{T}$, $\LDCostVec{\pi}{R}$ and $\LDCostTVec{\pi}{R}{T}$, represent the respective $m$-dimensional cost vectors whose elements indexed by the state $x\in\mathcal{S}$ are given by  \eqref{eq_RSCost_FinandinfHor} and \eqref{eq_LDCost_FinandinfHor}. 
%Define below the vector of corresponding costs over all initial states $x\in\mathcal{S}$:
%\begin{eqnarray}
%\RSCostVec{\pi}{\gamma}{R} &:= (\RSCost{\pi}{\gamma}{R})_{x \in {\cal X} }, \ \ \  \RSCostTVec{\pi}{\gamma}{R}{T}&:= (\RSCostT{\pi}{\gamma}{R}{T})_{x \in {\cal X}}, \label{eq_RSCost_FinandinfHor_Vec}\\
%\LDCostVec{\pi}{R}&:=(\LDCost{\pi}{R})_{x \in {\cal X} }, \mbox{ and, } \ \ \ \LDCostTVec{\pi}{R}{T}&:=(\LDCostT{\pi}{R}{T})_{x \in {\cal X} }.\label{eq_LDCost_FinandinfHor_Vec}
%\end{eqnarray}
%This paper considers risk-sensitive cost only for $\gamma\neq 0$.   

%Denote the vectorized linear and risk-sensitive cost over all the initial conditions as $\RSCostVec{\pi}{\gamma}{R},\LDCostVec{\pi}{R},\RSCostTVec,\LDCostTVec,\RSCostUSLD,\RSCostUSLDVec$, where the $x$ in the bracket is suppressed.
\subsection{Problem Statement}\label{sec_prb_stmt}
In this paper, we consider CRSMDPs involving a fairly general set of constraints. Specifically, we consider constraints involving RS and standard discounted costs over both finite and infinite horizons. To this end, given nonnegative integers $M$, $\hat{M}$, $\bar{M}$ and $\check{M}$ at least one of which is nonzero, we introduce a set of immediate cost functions $C_i$, $\hat{C}_j$, $\bar{C}_k$, $\check{C}_l$, constants $b_i$, $\hat{b}_j$, $\bar{b}_k$,  $\check{b}_l$, and finite times $\bar{T}_k$, $\check{T}_l$ for $i\in\{1,2,\ldots, M\}$,  $j\in\{1,2,\ldots, \hat{M}\}$, $k\in\{1,2,\ldots, \bar{M}\}$, and $l\in\{1,2,\ldots, \check{M}\}$.  We define our constraints in terms of subsets of $\Pi_{\rm MR}$ given by 
\begin{eqnarray}
 \FeasibleSetLD &:=& \{\pi\in\Pi_{\rm MR}:\LDCost{\pi}{C_i}\leq b_i \mbox{ for all } i=1,2,\ldots, M\},\label{eq_Infin_Hor_Lin_Const}\\
\FeasibleSetRS &:=&\{\pi\in\Pi_{MR}:\RSCost{\pi}{\gamma}{\hat{C}_j}\leq \hat{b}_j \mbox{ for all }j=1,2,\ldots, \hat{M}\}\label{eq_Infin_Hor_RSConst},\\
 \FeasibleSetLDFinHor &:=& \{\pi\in\Pi_{MR}:\LDCostT{\pi}{\bar{C}_k}{\bar{T}_k}  \leq \bar{b}_k \mbox{ for all }k=1,2,\ldots, \bar{M}\}\label{eq_Fin_Hor_Lin_Const}, \\
 \FeasibleSetRSFinHor&:=& \{\pi\in\Pi_{MR}: \RSCostT{\pi}{\gamma}{\check{C}_l}{\check{T}_l} \leq \check{b}_l \mbox{ for all }l=1,2,\ldots, \check{M}\}, \label{eq_Fin_Hor_RSConst}
\end{eqnarray}
where the notation follows the definitions in  \eqref{eq_RSCost_FinandinfHor}-\eqref{eq_LDCost_FinandinfHor}.
% {\color{red}\sout{Throughout, we append the symbols \hspace{0.05 cm}$\hat{}$\hspace{0.05 cm}, $\bar{}$ \hspace{0.05 cm} and \hspace{0.05 cm} $\check{}$\hspace{0.1 cm} for representing RS infinite horizon, standard discounted finite horizon and RS finite horizon based objects respectively.}}
Our assumption that $\mathcal{S}$ and $\mathcal{A}$ are finite implies that that all the immediate cost functions  appearing in  \eqref{eq_Infin_Hor_Lin_Const}-\eqref{eq_Fin_Hor_RSConst} as well as the immediate cost function $R$ are uniformly bounded above in absolute value by a common bound, which we denote by $C>0$.

We are interested in the infinite-horizon CRSMDP problem given by %~\eqref{eq_CRSMDP_Prob}:
\begin{equation} \label{eq_CRSMDP_Prob}\tag{$P$}
\min\limits_{\pi\in\FeasibleSetOrgnal}\RSCost{\pi}{\gamma}{R},
% %\begin{split}
% &\min\RSCost{\pi}{\gamma}{R}   \\
% &\mbox{subject to }   \pi\in \FeasibleSetLD \cap \FeasibleSetRS \cap \FeasibleSetLDFinHor\cap \FeasibleSetRSFinHor.
% %\end{split}
\end{equation}
where 
%Denote the feasible region of \eqref{eq_CRSMDP_Prob} by $\FeasibleSetOrgnal$, that is, 
$\FeasibleSetOrgnal:=\FeasibleSetLD \cap \FeasibleSetRS \cap \FeasibleSetLDFinHor\cap \FeasibleSetRSFinHor$.  Note that the CRSMDP problem \eqref{eq_CRSMDP_Prob} involves minimizing an infinite-horizon RS cost function subject to $M+\bar{M}+\hat{M}+\check{M}$ number of constraints involving finite- and infinite-horizon RS and standard discounted constraints.
 
Our approach is to approximate the infinite-horizon problem \eqref{eq_CRSMDP_Prob} by a finite-horizon problem where the objective function as well as the constraints are appropriately time truncated. For this purpose, denote $K=C/(1-\beta)$. For each  $T\in\{1,2,\ldots\}$, let $K_{T} = e^{|\gamma| K\beta^T}$,  and define 
\begin{eqnarray}
\FeasibleSetLDTruncBlw{T} &:=& \{\pi\in\Pi_{\rm MR}:\LDCostT{\pi}{C_i}{T}\leq b_i - K\beta^T \mbox{ for all } i=1,2,\ldots, M\}, 
\label{eq_defn_LDTrunBlw}\\
\FeasibleSetRSTruncBlw{T} &:=&\left\{\pi\in\Pi_{MR}:\RSCostT{\pi}{\gamma}{\hat{C}_j}{T}\leq \frac{\hat{b}_j}{K_{T}} \mbox{ for all }j=1,2,\ldots, \hat{M}\right\},\label{eq_defn_RSTrunBlw}\\
\FeasibleSetLDTruncAbv{T} &:=& \{\pi\in\Pi_{\rm MR}:\LDCostT{\pi}{C_i}{T}\leq b_i + K\beta^T \mbox{, for all } i=1,2,\ldots, M\}, \label{eq_defn_LDTrunAbv}\\
\FeasibleSetRSTruncAbv{T}&:=&\{\pi\in\Pi_{\rm MR}: \RSCostT{\pi}{\gamma}{\hat{C}_j}{T}\leq \hat{b}_j K_{T} \mbox{, for all }j=1,2,\ldots, \hat{M}\}.\label{eq_defn_RSTrunAbv}%\label{eq_defn_ovrlnDT_RSconst_trunc}
\end{eqnarray}

The subsets of $\Pi_{\rm MR}$ defined in \eqref{eq_defn_LDTrunBlw}-\eqref{eq_defn_RSTrunAbv} represent constraints defined in terms of the time-truncated versions of the cost functions appearing in  \eqref{eq_Infin_Hor_Lin_Const}-\eqref{eq_Infin_Hor_RSConst}. Taken together with the time-truncated version of the infinite-horizon RS cost of problem  \eqref{eq_CRSMDP_Prob}, the constraints \eqref{eq_defn_LDTrunBlw}-\eqref{eq_defn_RSTrunAbv} allow us to define, for each $T\geq 1$,  the finite-horizon problems  $\CRSMDPFinHorOne$ and $\CRSMDPFinHorTwo$ given by
%We approximating the problem \eqref{eq_CRSMDP_Prob} with two finite horizon CRSMDPs $\CRSMDPFinHorOne$ and $\CRSMDPFinHorTwo$ (given below), by truncating the infinite horizon objective and constraints as given in \eqref{eq_defn_LDTrunBlw} to \eqref{eq_defn_RSTrunAbv}. The two problems $\CRSMDPFinHorOne$ and $\CRSMDPFinHorTwo$ are defined as follows:
\vspace{-4mm}
%\begin{equation} \label{eq_CRSMDP_Prob_FinHor1}\tag{\mbox{$P_T^-$}}
%\min\limits_{\pi\in \FeasibleSetTrunBlw{T}} \RSCostT{\pi}{\gamma}{R}{T},
%\end{equation}

%\begin{equation} \label{eq_CRSMDP_Prob_FinHor2}\tag{\mbox{$P_T^+$}}
%\min\limits_{\pi\in \FeasibleSetTrunAbv{T}} \RSCostT{\pi}{\gamma}{R}{T},
%\end{equation}
%
\hypertarget{theequation}{\begin{equation}
(P_T^-): \min\limits_{\pi\in \FeasibleSetTrunBlw{T}} \RSCostT{\pi}{\gamma}{R}{T}\ \mbox{ and } \ (P_T^+):\min\limits_{\pi\in \FeasibleSetTrunAbv{T}} \RSCostT{\pi}{\gamma}{R}{T},\label{eq_test2}\vspace{-2mm}%\tag{\mbox{$P_T^+$}}
\end{equation}}
% Denote the feasible regions of the problems $\CRSMDPFinHorOne$ and $\CRSMDPFinHorTwo$ by $\FeasibleSetTrunBlw{T}$ and $\FeasibleSetTrunAbv{T}$ respectively, that is, 
where $\FeasibleSetTrunBlw{T}:=\FeasibleSetLDTruncBlw{T} \cap \FeasibleSetRSTruncBlw{T} \cap \FeasibleSetLDFinHor \cap \FeasibleSetRSFinHor$ and $\FeasibleSetTrunAbv{T}:=\FeasibleSetLDTruncAbv{T} \cap \FeasibleSetRSTruncAbv{T} \cap \FeasibleSetLDFinHor \cap \FeasibleSetRSFinHor$.  We will show that the finite-horizon problems $\CRSMDPFinHorOne$ and $\CRSMDPFinHorTwo$ serve as approximations to the infinite-horizon problem \eqref{eq_CRSMDP_Prob}. More precisely, we will show that the optimal values of the problems $\CRSMDPFinHorOne$ and $\CRSMDPFinHorTwo$ converge to the value of \eqref{eq_CRSMDP_Prob} as $T\rightarrow \infty$ under appropriate conditions. Moreover, the feasible sets of the problems are related in that, the feasible sets of $\CRSMDPFinHorOne$ for different values of $T$ form an increasing nested sequence contained in the feasible set of \eqref{eq_CRSMDP_Prob}, while the feasible sets of $\CRSMDPFinHorTwo$ form an decreasing nested sequence containing the feasible set of \eqref{eq_CRSMDP_Prob}. As we will show in the next section, these properties allow us to conclude that solutions of $\CRSMDPFinHorOne$ and $\CRSMDPFinHorTwo$ for sufficiently large $T$ serve as approximate solutions to \eqref{eq_CRSMDP_Prob}.

% {\color{red} I will come back to this paragraph later. Reading further might help me understand why this point is needed.}
As mentioned in the introduction, reference \cite{UdaySanjayKavitaNH} followed an identical approach to approximate the  unconstrained version of the problem \eqref{eq_CRSMDP_Prob} by the unconstrained version of the problem  $\CRSMDPFinHorOne$ or, equivalently, the unconstrained version of $\CRSMDPFinHorTwo$. 
 This paper extends the approach of   \cite{UdaySanjayKavitaNH}  to constrained RSMDPs. 
 
% We showed that \eqref{eq_CRSMDP_Prob} can be approximated uniformly by the unconstrained version of the problems $\CRSMDPFinHorOne$ as $T\to\infty$. Policies considered can again be seen as US policies, where the stationary part is any arbitrary policy and the non-stationary part is the first $T$ components of the optimal policy of $\CRSMDPFinHorOne$. The approximation is such that the stationary part doesn't have any impact on the objective function and this is because all immediate stage-wise costs from time $T$  onwards are neglected. Moreover, the rate of convergence of $\CRSMDPFinHorOne$ to \eqref{eq_CRSMDP_Prob} as $T\to\infty$ is also derived. In the current paper we consider the same problem with constraints. 

We state the main results in the next section.
\section{Main results}\label{sec_main_results}

%In this section, we state the three main results of the paper. 
%Recall that we assume \eqref{eq_CRSMDP_Prob} is feasible. 
Our first result below asserts that the CRSMDPs \eqref{eq_CRSMDP_Prob} and $\CRSMDPFinHorTwo$, for every $T$, possess solutions if \eqref{eq_CRSMDP_Prob} is feasible, while the CRSMDP $\CRSMDPFinHorOne$ has a solution if it is feasible. The proof is given in section \ref{sec_proof_mainres}. 

\begin{thm}\label{thm_existence_result}\textbf{[Existence of optimal solutions:]}
If the feasible region $\FeasibleSetOrgnal\subseteq \Pi_{\rm MR}$ of the original problem \eqref{eq_CRSMDP_Prob} is non-empty, then \eqref{eq_CRSMDP_Prob} has a solution and, for each $T>0$, $\CRSMDPFinHorTwo$ has a  solution. On the other hand, %feasible region of the problem $\CRSMDPFinHorOne$, 
if $\FeasibleSetTrunBlw{T^*} \neq \varnothing$ for some $T^*>0$, then $\FeasibleSetOrgnal\neq \varnothing$,  and $\CRSMDPFinHorOne$ has a solution for each $T\geq T^*$.
\end{thm}

While Theorem \ref{thm_existence_result} guarantees the existence of a solution to \eqref{eq_CRSMDP_Prob}, computing a solution is a  challenge.  One reason is that the solution is generally non-stationary (see, for instance, \cite{SJaq1976}). To the best of our knowledge, there is no numerical method to compute such a solution. One expects that if the optimal policies are stationary or of the US type, then it may be relatively easier to compute the solution.
Hence, it is natural to ask: i) can US policies approximate the optimal value of the problem \eqref{eq_CRSMDP_Prob}; and ii) can a numerical method be devised to compute these approximate solutions? In this paper we provide affirmative answers to  both these questions. Specifically, we show below that any policy that solves either  $\CRSMDPFinHorOne$ or $\CRSMDPFinHorTwo$ for sufficiently large $T$ can be extended to a US policy that approximately solves the original problem \eqref{eq_CRSMDP_Prob}. Furthermore, we discuss the solution of $\CRSMDPFinHorOne$ and $\CRSMDPFinHorTwo$ using a LP-based approach in section \ref{sec_algo}. Our next two results show how  $\CRSMDPFinHorOne$ and $\CRSMDPFinHorTwo$ approximate \eqref{eq_CRSMDP_Prob}.

Theorem \ref{thm_finhorapprx_frmblw} given below is our second main result. It gives sufficient conditions under which a solution of $\CRSMDPFinHorOne$ yields an approximate solution to \eqref{eq_CRSMDP_Prob} for $T$ sufficiently large.  The sufficient condition involves checking if 0 is a local minimum value for the {\textit{ maximum constraint violation}} map $h:\Pi_{\rm MR} \mapsto \R$ which maps a policy $\pi\in\Pi_{\rm MR}$ to the largest constraint violation  achieved by $\pi$ across all the $M+\bar{M}+\hat{M}+\check{M}$ constraints that form part of the problem \eqref{eq_CRSMDP_Prob}. The qualifier ``local" above refers to the topology on $\Pi_{\rm MR}$ that we introduce in  section \ref{sec_Top_Markvpolc}. 
\begin{thm}\label{thm_finhorapprx_frmblw}
\textbf{[Finite-horizon approximation with feasibility:]} 
% {\color{blue} Define $h:\Pi_{\rm MR}\to \R$ by
% \begin{equation}
%     h(\pi):=\max\{\theta(\pi), \hat{\theta}(\pi),\bar{\theta}(\pi),\check{\theta}(\pi)\}, \label{eq_defn_of_h}
% \end{equation}
% where the real-valued maps $\theta,\hat{\theta},\bar{\theta},\check{\theta}$ on $\Pi_{\rm MR}\to \R$ are given by $\theta(\pi):=& \max_{i\in\{1,\ldots,M\}} \{\LDCost{\pi}{C_i}- b_i  \}$, $\hat{\theta}(\pi):=& \max_{j\in\{1,\ldots,\hat{M}\}} \{\RSCost{\pi}{\gamma}{\hat{C}_j} -\hat{b}_j\}$, $\bar{\theta}(\pi):=& \max_{k\in\{1,\ldots,\bar{M}\}}\{\LDCostT{\pi}{\bar{C}_k}{\bar{T}_k} - \bar{b}_k\}$ and $\check{\theta}(\pi):=& \max_{l\in\{1,\ldots,\check{M}\}}\{\RSCostT{\pi}{\gamma}{\check{C}_l}{\check{T}_l}- \check{b}_l\}$. }
Define the maps $\theta,\hat{\theta},\bar{\theta},\check{\theta},h:\Pi_{\rm MR}\to \R$ by
\begin{eqnarray}
\theta(\pi):=&\hspace{-0.4cm} \max\limits_{i\in\{1,\ldots,M\}} \{\LDCost{\pi}{C_i}- b_i  \}; \mbox{ }
\hat{\theta}(\pi):=&\hspace{-0.4cm} \max\limits_{j\in\{1,\ldots,\hat{M}\}} \{\RSCost{\pi}{\gamma}{\hat{C}_j} -\hat{b}_j\},\nonumber\\
\bar{\theta}(\pi):=& \hspace{-0.4cm}\max\limits_{k\in\{1,\ldots,\bar{M}\}}\{\LDCostT{\pi}{\bar{C}_k}{\bar{T}_k} - \bar{b}_k\}; \mbox{ }
\check{\theta}(\pi):=&\hspace{-0.4cm} \max\limits_{l\in\{1,\ldots,\check{M}\}}\{\RSCostT{\pi}{\gamma}{\check{C}_l}{\check{T}_l}- \check{b}_l\}, \nonumber\\
\mbox{and}
\hspace{1 cm}h(\pi):=&\hspace{-0.7 cm}\max\{\theta(\pi), \hat{\theta}(\pi),\bar{\theta}(\pi),\check{\theta}(\pi)\}.\label{eq_defn_of_h}%\hspace{4 cm}\nonumber
 \end{eqnarray}
Suppose  $\mathcal{F}\neq \varnothing$ and $0$ is not a local minimum value of the map $h$. Then, for every $\epsilon>0$, there exists $T^*>0$ such that, for every $T\geq T^*$ and every  policy $\eta\in\FeasibleSetTrunBlw{T}$ that solves the problem $\CRSMDPFinHorOne$, $\eta$ is $\epsilon$-optimal for the problem \eqref{eq_CRSMDP_Prob}, that is,  $\eta\in\mathcal{F}$ and $\eta$ satisfies
\begin{equation}
0\leq \RSCost{\eta}{\gamma}{R}-\inf\limits_{\pi\in\mathcal{F}}\RSCost{\pi}{\gamma}{R}<\epsilon.
    \label{eq_truncapprox}
\end{equation} 
\end{thm}

Theorem \ref{thm_finhorapprx_frmblw} guarantees a feasible and $\epsilon$-optimal solution to the original infinite-horizon problem \eqref{eq_CRSMDP_Prob}. However, the theorem requires checking a  local minimum condition on the function $h$ given by \eqref{eq_defn_of_h}, which could be difficult to check in practice. At the same time, we show through a counterexample in \ref{appendixD} that the conclusion of the theorem need not hold if the condition is violated. Our next result provides an alternative means of obtaining an approximate solution to problem \eqref{eq_CRSMDP_Prob} in situations where the local minimum condition on $h$ either does not hold or is  difficult to verify.  The theorem shows that, for sufficiently large $T$,  a solution of  $\CRSMDPFinHorTwo$ provides a policy which is approximately optimal and approximately feasible for the problem \eqref{eq_CRSMDP_Prob}, where approximate feasibility, made precise below, means that the maximum constraint violation can be bounded by an arbitrarily specified small quantity.
%a possibly infeasible  policy becomes feasible if  all the constraints are relaxed by a small amount. 
%The following definition makes this notion precise. 
% \footnote{Note that a policy in $\Pi_{\rm MR}$ is feasible, if and only if, it $\epsilon$-feasible for every $\epsilon>0$.}.
\begin{defn}\label{def_eps_feasb}
\textbf{[$\epsilon$-feasibility:] } Given $\epsilon>0$, a policy $\pi\in\Pi_{\rm MR}$ is \textbf{$\epsilon$-feasible} for the problem \eqref{eq_CRSMDP_Prob}  if it satisfies  
$$\LDCost{\pi}{C_i}\leq b_i+\epsilon, \hspace{0.2 cm} \RSCost{\pi}{\gamma}{\hat{C}_j}\leq \hat{b}_j+\epsilon,  \hspace{0.2 cm}\LDCostT{\pi}{\bar{C}_k}{\bar{T}_k}  \leq \bar{b}_k+\epsilon,\hspace{0.2 cm} \mbox{ and } \RSCostT{\pi}{\gamma}{\check{C}_l}{\check{T}_l} \leq \check{b}_l+\epsilon$$
for all $i=1,\ldots,M$, $j=1,\ldots,\hat{M}$, $k=1,\ldots,\bar{M}$ and $l=1,\ldots,\check{M}$.
\end{defn}

% \begin{thm}\label{thm_finhorapprx_frmabv}
% \textbf{[Finite Horizon Approximation From Above:]} If $\mathcal{F}\neq \varnothing$ then, for every $\epsilon>0$, there exists $T>0$ and an US policy $\eta\in\FeasibleSetTrunAbv{T}$ such that $\eta$ is a solution to the problem
% $\CRSMDPFinHorTwo$, is $\epsilon$-feasible for the problem \eqref{eq_CRSMDP_Prob}, and satisfies
% \begin{equation}
%     \left|\RSCost{\eta}{\gamma}{R}-\inf\limits_{\pi\in\mathcal{F}}\RSCost{\pi}{\gamma}{R}\right|<\epsilon.\label{eq_truncapprox2}
% \end{equation} 
% \end{thm}

\begin{thm}\label{thm_finhorapprx_frmabv}
\textbf{[Finite-horizon approximation with $\epsilon$-feasibility:]} Suppose $\mathcal{F}$ is nonempty.  Then, for every $\epsilon>0$, there exists $T^{*}>0$ such that, for every $T\geq T^{*}$ and every  policy $\eta\in\FeasibleSetTrunAbv{T}$ such that $\eta$ is a solution to the problem
$\CRSMDPFinHorTwo$, $\eta$ is $\epsilon$-feasible for the problem \eqref{eq_CRSMDP_Prob}, and satisfies
\begin{equation}
    \left|\RSCost{\eta}{\gamma}{R}-\inf\limits_{\pi\in\mathcal{F}}\RSCost{\pi}{\gamma}{R}\right|<\epsilon.\label{eq_truncapprox2}
\end{equation} 
\end{thm}

It is worthwhile to compare theorems \ref{thm_finhorapprx_frmblw} and \ref{thm_finhorapprx_frmabv}. Both theorems show that an approximate solution to the infinite-horizon CRSMDP  \eqref{eq_CRSMDP_Prob} can be obtained by solving either of the finite-horizon CRSMDPs $\CRSMDPFinHorOne$ or  $\CRSMDPFinHorTwo$ for a sufficiently large finite horizon $T$. While Theorem \ref{thm_finhorapprx_frmblw} requires an additional condition to hold, the approximate solution  it provides is feasible for the original problem  \eqref{eq_CRSMDP_Prob}. In contrast, Theorem \ref{thm_finhorapprx_frmabv} does not involve any additional conditions, but yields an approximate solution that is only approximately feasible in the sense of Definition \ref{def_eps_feasb}. 
% {\color{red}
% \begin{rem} \sout{Observe from \eqref{eq_defn_LDTrunBlw}-\eqref{eq_defn_RSTrunAbv} that, given any policy $\pi\in\Pi_{\rm MR}$, the feasibility as well as optimality of $]pi$ for the problems $\CRSMDPFinHorOne$ and  $\CRSMDPFinHorTwo$     %given any feasible policy $\pi = (d_0, d_1, \ldots, d_{T-1}, d_T, \ldots)$ of
%  are determined only by the first $T$ decision rules from $\pi$. 
% Therefore, without any loss of generality, we may restrict our attention to only those solutions of $\CRSMDPFinHorOne$ and  $\CRSMDPFinHorTwo$ that are US, with stationarity beginning from the decision epoch $T$. }
% \label{rem_us}
% \end{rem}

\begin{rem} Observe from \eqref{eq_defn_LDTrunBlw}-\eqref{eq_defn_RSTrunAbv} that, given any policy $\pi\in\Pi_{\rm MR}$, the feasibility as well as optimality of $\pi$ for the problems $\CRSMDPFinHorOne$ and  $\CRSMDPFinHorTwo$     %given any feasible policy $\pi = (d_0, d_1, \ldots, d_{T-1}, d_T, \ldots)$ of
 are determined only by the first $T$ decision rules of $\pi$. Consequently, appending arbitrary decision rules after the horizon $T$ to a solution of $\CRSMDPFinHorOne$ or  $\CRSMDPFinHorTwo$ does not affect its feasibility or optimality with respect to $\CRSMDPFinHorOne$ or  $\CRSMDPFinHorTwo$. In particular, one may extend a solution of either problem to a US policy by choosing all decision rules beyond the horizon $T$ to be the same. Hence, the approximate solutions to the problem \eqref{eq_CRSMDP_Prob} guaranteed by theorems \ref{thm_finhorapprx_frmblw}-\ref{thm_finhorapprx_frmabv} may be chosen to be US policies. \label{rem_us1}
% Therefore, the existence and approximated solutions that are given in Theorems   \ref{thm_existence_result}-\ref{thm_finhorapprx_frmabv} to the problems $\CRSMDPFinHorOne$ and  $\CRSMDPFinHorTwo$ are characterized upto only a finite time $T$. One can append any arbitrary decisions after time $T$. However, for implementation and simplicity purposes, one can choose the same decision after time $T$ making the approximated policy an ultimately stationary type.\label{rem_us1}
\end{rem}

\begin{rem}
Theorems \ref{thm_existence_result}, \ref{thm_finhorapprx_frmblw} and \ref{thm_finhorapprx_frmabv} hold as stated even in the case where the discount factors and risk factors used in \eqref{eq_RSCost_FinandinfHor} and \eqref{eq_Infin_Hor_Lin_Const}-\eqref{eq_Fin_Hor_RSConst} are different. While the modifications needed to accommodate the added generality is only minor, we have sacrificed the slight increase in generality in favour of simplicity  
%restricted the treatment to the case of a single discount and risk factor for simplicity. 
The case of multiple discount factors is considered in \cite{Piunovskiy} for CMDPs in the  standard-discounted-cost setting.
\end{rem}

% \begin{rem}
% Note that the results of this section are stated for the case where the risk-sensitive objective function is to be minimized, and the constraints \eqref{eq_Infin_Hor_Lin_Const}-\eqref{eq_Fin_Hor_RSConst} involve upper bounds. Problems involving maximization or lower bounds as constraints can be easily recast to fit our setting by reversing the signs of the immediate cost functions, and suitably redefining the constants defining the bounds in the constraints. 
% \end{rem}

The next section provides proofs of the three main results stated in this section. 
\section{Proofs of the Main Results}\label{sec_proof_mainres}

In this section, we provide proofs of theorems \ref{thm_existence_result}, \ref{thm_finhorapprx_frmblw} and \ref{thm_finhorapprx_frmabv}. In particular, we motivate the main steps in the proofs, and provide the auxiliary results which are needed to complete the proofs.

\subsection{Proof Outline of Theorem \ref{thm_existence_result}}\label{thm1ssec}% for Proof of Main Theorems}\label{thm1ssec}%{Proof of Theorem \ref{thm_existence_result}}

The proof of Theorem \ref{thm_existence_result} uses the well known fact that a continuous function achieves its infimum on a nonempty compact set. Thus the major steps in proving Theorem \ref{thm_existence_result} are: (i) constructing a topology on the set of  policies $\Pi_{\rm MR}$, (ii) showing that the objective functions in \eqref{eq_CRSMDP_Prob}, $\CRSMDPFinHorOne$ and $\CRSMDPFinHorTwo$ are continuous in this topology, and (iii) showing that the feasible regions are compact in the same topology.
Steps (i), (ii) and (iii) above are achieved in sub-sections \ref{sec_Top_Markvpolc} - \ref{subsec_cpctsssec} below.  Before proceeding, we introduce the required notations. 

\subsubsection{Notations}\label{sec_matrxNot}
The set $\R^{m\times n}$ of all real matrices of order $m\times n$ is a real vector space under element-wise addition and scalar multiplication. The matrix $\one_{m\times n}$  is the $m\times n$ matrix with all elements equal to $1$. Any vector in $\R^n$ is represented by a column vector.

For any $B=[b_{i,j}]\in \R^{m\times n}$ and  $y=[y_i]\in\R^{n}$, $y^{\prime}$ and $B^{\prime}$, respectively represent their transposes. The $l_\infty$ norm (or max norm) and the $l_{1}$ norm of $y\in\R^{n}$ are defined by  $\|y\|_\infty := \max\limits_{i=1,\ldots,n} |y_i|$ and $\|y\|_{1}=\sum_{i=1}^{n}|y_{i}|$, respectively. We denote by $\|B\|_{\infty}$ the matrix norm of $B\in\R^{m\times n}$ induced by $l_\infty$ which equals  $\max\limits_{i=1,\ldots, m} \sum_{j=1}^n|b_{i,j}|$ (Refer \cite[Ex. 5.6.5, p. 345]{HJ}).
We also define $\|B\|_{\rm max} = \max\limits_{i,j}|b_{i,j}|$. 
For any two matrices $B_1, B_2\in \R^{m\times n}$, the Schur product, denoted by $B_1\odot B_2$, is the matrix obtained by element-wise multiplication of matrices $B_1$ and $B_2$. Similarly, the Schur exponential, denoted  by $e^{\odot B}$, is the element-wise exponential of the matrix $B$ . 
%\label{eq_matrixnorminfty_eqls_maxrowsum}
%\footnotetext{This notation is used in \eqref{eq_matrixnorminfty_eqls_maxrowsum} many books, }
%Denote, 
%For any matrix $A$ and (column) vector $x$, we denote the following norms:
% \begin{eqnarray}
% \|A\|_{\rm max} = \max\limits_{i,j}|a_{ij}|, \ \ \ \ \|A\|_\infty = \max_{i}\sum_{j} |a_{ij}|\nonumber\\
% \|x\|_1 = \sum_i|x_i|, \ \ \ \ \|x\|_\infty = \max\limits_{i}|x_i|\nonumber
% \end{eqnarray}
\subsubsection{Construction of a Topology on Policy Space }\label{sec_Top_Markvpolc}%$\Pi_{\rm MR}$}\label{sec_Top_Markvpolc}
The topological structure that we define on $\Pi_{\rm MR}$ is based on the observation that an MR policy is a sequence of MR decision rules, one for each decision epoch, while each MR decision rule can be represented by a row-stochastic matrix of dimension $m \times n$, (recall that $m = |\mathcal{S}|$ and $n = |\mathcal{A}|$). This allows us to construct a metric, and hence a topology, on $\Pi_{\rm MR}$ by  using the matrix norm $\|\cdot\|_{\infty}$. To this end, let $\mathcal{R} \subset \R^{m\times n}$ denote the set of row-stochastic matrices, that is, \vspace{-2mm}
\begin{equation}
\mathcal{R} := \bigg\{B=[b_{i,j}]\in \R^{m\times n} : b_{i,j}\geq 0\hspace{0.1 cm} \forall i, j \mbox{ and } \sum\limits_{j=1}^{n} b_{i,j} = 1 \hspace{0.1 cm}\forall i\bigg \}.\vspace{-2mm}
\end{equation}

As a subset of $\mathbb{R}^{m\times n}$, $\mathcal{R}$ is  closed and bounded, and hence compact,  in the topology induced by the norm $\|\cdot\|_{\infty}$. 
% Closedness follows from the nature of the inequalities defining elements of $\mathcal{R}$, while boundedness follows from noting that the distance between any two row-stochastic matrices under the  maximum  absolute row-sum norm  $\|\cdot\|_{\infty}$ is bounded above by 2. Thus, $\mathcal{R}$ is a compact subset of $\mathbb{R}^{m\times n}$. {\color{red} Sanjay to modify this slightly to fit in with the proof of compactness of $\Pi_{\rm MR}$ better.}

% {\color{blue}
% Note that $\mathcal{R}$ is a closed and bounded and hence a compact subset of $\mathbb{R}^{m\times n}$ in the topology induced by norm $\|\cdot\|_{\infty}$.
% Closedness follows from the nature of the inequalities defining elements of $\mathcal{R}$, while boundedness follows from noting that the distance between any two row-stochastic matrices under the {\color{red}norm } $\|\cdot\|_{\infty}$ is bounded above by 2. 
 %Thus, $\mathcal{R}$ is a compact subset of $\mathbb{R}^{m\times n}$. %{\color{red} Sanjay to modify this slightly to fit in with the proof of compactness of $\Pi_{\rm MR}$ better.}
%}

Next, we define a sequence of equivalent metrics on $\mathcal{R}$. 
Fix $\delta\in(\beta,1)$. For each  $t\in\{0,1,\ldots\}$, define a metric $\mu_t$ on $\mathcal{R}$  by letting 
%\begin{equation}
  $\mu_t(B_1,B_2)=\delta^t\|B_1-B_2\|_\infty$ for $B_1,B_2\in\mathcal{R}$.
  %\label{eq_defn_metric_on_Markvdcns}
%\end{equation}
Note that, for each $t>0$,  
% the metric $\mu_{t}$ is equivalent to the metric induced on $\mathcal{R}$ by the norm $\|\cdot\|_{\infty}$, and hence
the metric space  $(\mathcal{R},\mu_t)$ is a compact metric space with diameter $2\delta^t$.

% {\color{blue} For any $B_1,B_2\in\mathcal{R}$ and $t=0,1,\ldots$, define an equivalent metric $\mu_t(B_1,B_2):=\delta^t\|B_1-B_2\|_\infty$ where $\delta\in(\beta,1)$ is a fixed constant. It is easy to see that $(\mathcal{R},\mu_t)$ is a compact metric space with diameter $2\delta^t$.

% }

As used by \cite{Derman_Klein}, we characterize an MR decision rule by a row-stochastic matrix. We identify an MR decision rule $d:\mathcal{S}\rightarrow \mathcal{P}(\mathcal{A})$ with the unique matrix in $\mathcal{R}$ whose $(i,j)$th element is the probability of taking action $a_{j}$ at state $s_{i}$ under the decision rule $d$. It is clear that the identification just described sets up a bijection between $\mathcal{R}$ and the set of MR decision rules. As a departure from the convention of using upper case letters for matrices,  we denote the matrix corresponding to an MR decision rule $d$ by $d$ again, and use the suggestive notation $d(a_{j}|s_{i})$ denote its $(i,j)$th element. %For all other matrices, other than MR decision rule, we follow, as usual, upper case letters.   %{\color{red} Again with slight abuse of notation, we represent the single subscript MR decision rule such as $d_t$ for decision rule at time epoch $t$ and at the same time double subscript $d_{i,j}$ as its $\{i,j\}^{th}$ element.}

% by An element $d\in\mathcal{R}$ can be identified with the unique MR  decision rule under which, for every $i=1,\ldots,m$ and $j=1,\ldots,n$,  the probability of taking action $a_j$ at  the current state $s_i$ equals $d_{i,j}$. Thus, any MR decision rule can be uniquely identified with an element of  $\mathcal{R}$ and vice versa. Hence, we can write an MR decision rule $d$ in matrix form as:
% \begin{equation}
% d=[d_{i,j}]_{m\times n}=\quad
% \begin{bmatrix}
% d(a_1|s_1) & d(a_2|s_1) & \cdots & d(a_n|s_1) \\
% d(a_1|s_2) & d(a_2|s_2) & \cdots & d(a_n|s_2)\\
% \vdots     & \vdots     &        & \vdots\\
% d(a_1|s_m) & d(a_2|s_m) & \cdots & d(a_n|s_m)
% \end{bmatrix}.
% \quad
% \end{equation}

Next, it is easy to see that a MR policy can be identified with a unique sequence having elements in $\mathcal{R}$. We therefore identify $\Pi_{\rm MR}$ with the set  $\mathcal{R}^\infty$  of all sequences having elements in $\mathcal{R}$. In the rest of the paper, we will use $\Pi_{\rm MR}$ interchangeably with  $\mathcal{R}^\infty$. Using this identification, we now define a metric on $\Pi_{\rm MR}$ as follows.  Given policies $\pi_1 =\{d_t\}_{t=0}^{\infty}$ and $\pi_2=\{f_t\}_{t=0}^{\infty}$ in $\Pi_{\rm MR}(=\mathcal{R}^\infty)$, define
%\begin{equation}
$\mu(\pi_1,\pi_2):=\sup\limits_{t\geq 0}\mu_t(d_t,f_t)$.
%\label{eq_defn_metric_mu}. 
%\end{equation}
Our next result shows that $\mu$ is a metric on $\Pi_{\rm MR}$, and provides the foundation for the framework that we use to prove our main results.  
%{\color{red} There is a mild resemblance of $\Pi_{\rm MR}$ with a subset of $l^\infty$ on $\R^{m\times n}$, the set of bounded sequences. However, the standard metric defined in $l^\infty$ is not equivalent to $\mu$.  } 

%\begin{thm}\label{thm_metric_on_Pi_metrizes_ProdTop}
%The map $\mu$ is a metric on $\Pi_{\rm MR}$, and metrizes the product topology of the Cartesian product $\prod\limits_{t= 0}^{\infty}(\mathcal{R},\mu_t)$. Moreover, $\Pi_{\rm MR}$ is a compact metric space.
%\end{thm}

\begin{thm}\label{thm_metric_on_Pi_metrizes_ProdTop}
The map $\mu$ is a metric on $\Pi_{\rm MR}$, and $(\Pi_{\rm MR},\mu)$ is a compact metric space.
\end{thm}

It is  well-known that any Cartesian product of compact spaces is compact in the product topology. Theorem \ref{thm_metric_on_Pi_metrizes_ProdTop}  simply follows by showing  that the metric $\mu$ constructed above metrizes the product topology on $\mathcal{R}^\infty$. The choice of the metric is thus crucial to Theorem \ref{thm_metric_on_Pi_metrizes_ProdTop}. For instance, compactness fails to hold if $\mu$ is chosen to be the more familiar $l^{\infty}$ metric on $\mathcal{R}^{\infty}$ obtained by setting $\delta=1$.

\subsubsection{Continuity of standard discounted cost functions} \label{cctsssec}
To show that the finite- and infinite-horizon standard discounted cost functions appearing in problems \eqref{eq_CRSMDP_Prob}, $\CRSMDPFinHorOne$ and $\CRSMDPFinHorTwo$ are continuous, we exploit well-known expressions for such a cost function  which use matrix-vector notation to  explicitly bring out the dependence on the decision rules that constitute the policy (see, for instance, \cite[eqn. (6.1.2)]{PUT}). To do so, we will need some additional notation. 

 First, given $(s,a)\in\mathcal{S}\times \mathcal{A}$, we denote by $p(\cdot|s,a)\in\R^{m}$ the vector whose $i$th element is the transition probability $p(s_{i}|s,a)$.  Given a immediate cost function $R:\mathcal{S}\times \mathcal{A}\rightarrow \R$, we abuse the notation slightly and denote by $R\in\R^{m\times n}$ the matrix whose $(i,j)$th element is given by $R(s_{i},a_{j})$.  
 Following  \cite[sec. 5.6]{PUT}, the  vector of expected costs $R_{d}\in\R^{m}$ and the transition probability matrix $P_{d}\in\R^{m\times m}$ corresponding to a MR decision rule $d\in\mathcal{R}$ are defined element-wise by 
\begin{eqnarray}
(R_d)_{i} := \sum\limits_{j=1}^{n}R(s_{i},a_{j})d(a_{j}|s_{i}),\   
(P_d)_{i,j} := \sum\limits_{k=1}^n p(s_j|s_i,a_k)d(a_k|s_i).
\label{eq_def_of_vec_rd}
\end{eqnarray}
Observe that $R_{d}=(R\odot d)\one_{n\times 1}$.

%For $d\in\mathcal{R}$, one can treat $R_d$ as vector valued mapping considered over state space and $P_d$ as matrix valued mapping.

Next, given a MR policy $\pi=\{d_t\}_{t=0}^{\infty}$ define the $t$-step transition probability matrix $P^\pi_t :=  \prod_{0 < l < t}  P_{d_l}$ for each $t=1,2,\ldots$, and set $P^{\pi}_{0}$ to be the $m\times m$ identity matrix. 
%Our next lemma shows that the vector $R_d$ and matrix $P_d$ appearing in \eqref{eq_expect_in_matrxiandvector_form} and \eqref{eq_expect_in_matrxiandvector_form} are Lipschitz continuous as functions of the MR decision rule $d\in\mathcal{R}$.
The vector expressions for the standard discounted costs $\LDCostVec{\pi}{R}$ and $\LDCostTVec{\pi}{R}{T}$ defined in \eqref{eq_LDCost_FinandinfHor} may now be written as \vspace{-2mm}
\begin{eqnarray}
\LDCostTVec{\pi}{R}{T} =  \sum\limits_{t=0}^{T-1}\beta^t P^{\pi}_{t}R_{_{d_t}}, \ %\label{eq_expect_in_matrxiandvector_form} % R_{d_0} + \beta P_{{d_0}} R_{{d_1}}  +\cdots + \beta^{T-1} (P_{{d_0}}P_{{d_1}}\cdots P_{{d_{T-2}}} ) R_{{d_{T-1}}} 
\LDCostVec{\pi}{R} = \sum\limits_{t=0}^{\infty}\beta^t P^{\pi}_{t}R_{{d_t}}. \label{eq_expect_in_matrxiandvector_form} % \label{eq_expect_Finhor_in_matrxiandvector_form} R_{d_0} + \beta P_{{d_0}} R_{{d_1}} + \beta^2 P_{{d_0}}P_{{d_1}} R_{{d_2}} +\cdots =
\end{eqnarray}%\vspace{-2mm}
%The expression \eqref{eq_expect_in_matrxiandvector_form} is well known, and appears in \cite[Ch. 6]{PUT}.
The expressions above are well known and given, for instance, in \cite[Ch. 6]{PUT}. Note that $P_t^{\pi}R_{d_t}$, the $t$th term  in the summation for  $\LDCostTVec{\pi}{R}{T}$ above, is the expectation of the immediate cost incurred at time $t$ under the policy $\pi$. The continuity of $\LDCostTVec{\pi}{R}{T}$  in \eqref{eq_expect_in_matrxiandvector_form} thus depends on how the  expected immediate cost at time $t$ changes when the decision rules up to time $t$ are changed. The bound in the next lemma answers this question.
 
\begin{lem}\label{lemma_norm_rdminrf_PdminPf}
Let $\pi_1=\{d_t\}_{t=0}^{\infty}$ and $\pi_2=\{f_t\}_{t=0}^{\infty}$ be two policies in $\Pi_{\rm MR}$, and let $t\geq 0$. Then, 
$\|P_t^{\pi_1}R_{d_t} -P_t^{\pi_2}R_{f_t}\|_\infty \leq C\sum\limits_{i=0}^t \|d_i-f_i\|_\infty$ holds.
\end{lem}

The inequality above
%asserted by  Lemma \ref{lemma_norm_rdminrf_PdminPf} 
along with the expression \eqref{eq_expect_in_matrxiandvector_form} leads to Lipschitz continuity of  finite-horizon standard discounted cost functions stated in Theorem \ref{thm_contuty_lin_disc_cost_prod_toplogy} below. The proof of continuity of infinite-horizon standard discounted cost functions depends on the following result which states that the finite-horizon standard discounted cost converges to the infinite-horizon standard discounted cost as the horizon $T$ tends to infinity, uniformly in the policy $\pi$. While the result appears to be well known, we include a proof in \ref{appendixA} for the sake of completeness.   

\begin{lem}\label{lem_linCost_finhor_cgsto_infinhor}
\textbf{[Convergence of finite-horizon standard discounted cost:] } Let $x\in\mathcal{S}$ and  $R:\mathcal{S}\times \mathcal{A}\rightarrow \R$ be an immediate cost function.
Then, for every  $\epsilon>0$, there exists $T^*>0$ such that $|\LDCostT{\pi}{R}{T} -\LDCost{\pi}{R}|< \epsilon$ for all $T\geq T^*$ and all $\pi\in\Pi_{\rm MR}$. 
\end{lem}

The continuity properties of finite- and infinite-horizon standard discounted costs are summarized in our next result. 

%\begin{thm}\label{thm_contuty_lin_disc_cost_prod_toplogy}
%Let $T>0$. Then the vector-valued finite-horizon standard discounted cost $\pi\to \LDCostTVec{\pi}{R}{T}$, mapping the metric space $(\Pi_{\rm MR},\mu)$ to the metric space $(\R^{m},\|\cdot\|_{\infty})$, is Lipschitz continuous with Lipschitz constant $\frac{C(\delta^{T}-\beta^{T})}{\delta^{T-1}(\delta-\beta)(1-\beta)}$. Furthermore, the vector-valued infinite-horizon standard discounted cost $\pi\to \LDCostVec{\pi}{R} $  mapping the metric space $(\Pi_{\rm MR},\mu)$ to the metric space $(\R^{m},\|\cdot\|_{\infty})$ is Lipschitz continuous with Lipschitz constant $C\delta (\delta-\beta)^{-1}(1-\beta)^{-1}$.
%\end{thm}
\begin{thm}\label{thm_contuty_lin_disc_cost_prod_toplogy}
\textbf{[Continuity of standard discounted costs:]}
Let $T>0$. The mappings $\pi\mapsto \LDCostTVec{\pi}{R}{T}$ and $\pi\mapsto \LDCostVec{\pi}{R} $ from the metric space $(\Pi_{\rm MR},\mu)$ to the metric space $(\R^{m},\|\cdot\|_{\infty})$ are Lipschitz continuous with Lipschitz constants $\frac{K(\delta^{T}-\beta^{T})}{\delta^{T-1}(\delta-\beta)}$ and $K\delta (\delta-\beta)^{-1}$, respectively.
\end{thm}

\subsubsection{Continuity of RS cost functions}\label{rsctsssec}
% \sout{ As in the case of standard discounted cost functions in sub-subsection \ref{cctsssec} above, the starting point for proving continuity of risk-sensitive discounted cost functions is to express such a cost function in a manner that makes the dependence on the decision rules of the policy explicit. While matrix-vector expressions in \eqref{eq_expect_in_matrxiandvector_form} are  well-known in the case of  standard discounted cost \cite{PUT}[Chapter 5 and 6], similar expressions do not appear to be available for the risk-sensitive case. Our first result below fills this gap by giving a backward-recursion-based expression for the finite-horizon risk-sensitive discounted cost defined in Definition \ref{costdef}. The result may be viewed as a policy evaluation procedure for finite-horizon risk-sensitive discounted cost. }

 As in the previous subsection, we begin by expressing RS cost functions in a way  that makes the dependence on decision rules explicit.
 To this end, we introduce the risk-sensitive version  of the familiar  state-action value function or $Q$-factor. Given a finite-horizon $T>0$ and a time instant $t\in\{0,\ldots,T-1\}$, the RS state-action value function $\Qfunc{t}{T}{\pi}:\mathcal{S}\times \mathcal{A}\rightarrow \R$ of a policy $\pi\in\Pi_{\rm MR}$ at the state-action pair $(s,a)\in\mathcal{S}\times \mathcal{A}$ is the expected RS cost incurred between $t$ to $T$ when the action $a$ is applied at time $t$ in the state $s$, and the policy $\pi$ is applied thereafter. More precisely, given $(s,a)\in\mathcal{S}\times \mathcal{A}$ and $t,T$ as above, we have 
\begin{equation}
    \Qfunc{t}{T}{\pi}(s,a):=\E^\pi\left[e^{\gamma\sum\limits_{\tau=t}^{T-1}\beta^\tau R(X_\tau,A_\tau)}\bigg|X_{t}=s,A_{t}=a\right].
\label{eq_qvaluedef}
\end{equation}
As in the case of the immediate cost function $R$ in sub-subsection \ref{cctsssec}, it will be convenient to treat $\Qfunc{t}{T}{\pi}$ as the $ m\times n$ matrix having $\Qfunc{t}{T}{\pi}(s_{i},a_{j})$ as its $(i,j)$th entry. 

% As in the case of standard discounted cost (e.g., \cite{PUT}), 
The next proposition provides recursive expressions for  the risk-sensitive state-action value function as well as the finite-horizon RS cost for a given policy analogous to well-known recursive expressions for the standard discounted cost as given in \cite[eqn. (4.2.6)]{PUT}. 
% (notations in sub-section \ref{sec_matrxNot}) 
% \sout{ For stating the result, recall that $e^{\odot\gamma\beta^t R}$ denotes the  Schur exponential of the matrix $\gamma \beta^{t}R$ for any given $t\geq 0$, where $R$ is the immediate cost function expressed in matrix notation.}

%%%%%%%%%
%In this section, we derive a recursive, matrix-vector expression for the finite horizon risk-sensitive cost $\RSCostTVec{\pi}{\gamma}{R}{T}$ for any given policy $\pi\in\Pi_{\rm MR}$. 
%We begin by rewriting some of the elements of the CRSMDP using matrix-vector notation as in \cite{PUT}. We immediately have the following result.

\begin{prop}\label{thm_matrix_formulaton_FinHor_RSCost} Suppose $T\geq 1$, and let $\pi=\{d_t\}_{t=0}^{\infty}\in\Pi_{\rm MR}$. Then the following equations hold. 
\begin{eqnarray}
\Qfunc{T-1}{T}{\pi} &=& e^{\odot\gamma\beta^{T-1} R}, \label{eq_zeta0_eql_ones_matrix}\\
\Qfunc{t-1}{T}{\pi}(s,a) &=& e^{\gamma\beta^{t-1} R(s,a)} p(\cdot|s,a)^\prime (d_{t}\odot \Qfunc{t}{T}{\pi})\one_{n\times 1},
\nonumber\\& & \mbox{ for }  
 t=1,2,\ldots, T-1,\ (s,a)\in\mathcal{S}\times \mathcal{A}, \label{zetat_interms_zetatmin1} \\
\RSCostTVec{\pi}{\gamma}{R}{T} &=& 
 (d_0\odot \Qfunc{0}{T}{\pi})\one_{n\times 1}.\label{eq_JpiT_In_Matrix_Form_eql_expgammar_d0_zetaTmin1}
\end{eqnarray}
\end{prop}

Equations \eqref{eq_zeta0_eql_ones_matrix}-\eqref{eq_JpiT_In_Matrix_Form_eql_expgammar_d0_zetaTmin1}  provide a backward-recursion-based procedure for policy evaluation of finite-horizon RS  cost similar to that for standard discounted cost (see equation (22) of \cite{PUT}), and could be of independent interest. 
%{\color{red} We need to find something to call the recursive equations and connect them to something that appears in the literature. Can these be called or related to Bellman equation for RS cost? Didn't a previous paper by Professor Kavitha's student have some recursive equations, but written using expectations? }

It is clear from \eqref{eq_zeta0_eql_ones_matrix}-\eqref{eq_JpiT_In_Matrix_Form_eql_expgammar_d0_zetaTmin1} that continuity properties of the finite-horizon RS cost is determined by those of the matrices $\Qfunc{0}{T}{\pi},\ldots, \Qfunc{T-1}{T}{\pi}\in\R^{m\times n}$.   Our next result is a lemma that gives explicit bounds on the elements of these  matrices and their differences. Both bounds will be used for  proving that the finite-horizon RS cost defined in \eqref{eq_RSCost_FinandinfHor} depends continuously on the policy. 
The proof is provided in \ref{appendixA}.

\begin{lem}\label{lemma_zeta_and_its_differ_bounds}
%Let $\pi_1=\{d_t\}_{t\geq 0},\pi_2=\{f_t\}_{t\geq 0}\in \Pi_{\rm MR}$. 
Let $\pi_{1}=\{d_{t}\}_{t=0}^{\infty}$ and  $\pi_{2}=\{f_{t}\}_{t=0}^{\infty}$  be two policies in $\Pi_{\rm MR}$. 
Then, for every $T\geq 1$ and every  $0\leq t \leq T-1$, we have 
\begin{eqnarray}
%\|\zeta_{T}^{\pi_1}\|_{\rm max} & \leq & \exp\left(\frac{\gamma C(1-\beta^T)}{1-\beta}\right) \label{eq_maxnorm_zetaT_bound}\\
%\|\zeta_{T-\tau}^{\pi_1}\|_{\rm max} & \leq & e^{\gamma K(\beta^\tau-\beta^T)}, THIS IS LAST ZETA VERSION
\|\Qfunc{t}{T}{\pi_1}\|_{\max} &\leq & e^{|\gamma|K(\beta^{t} -\beta^T)},
\label{eq_maxnorm_zetaT_bound}\\
%\|\zeta_{T}^{\pi_1}-\zeta_{T}^{\pi_2}\|_{max}& \leq &  \exp\left(\frac{\gamma C(1-\beta^T)}{1-\beta}\right)\sum\limits_{t=0}^{T-1}\|d_t-f_t\|_\infty.\label{eq_maxnorm_zetaTOfpi1_min_norm_zetaTOfpi2_bound}\\
%\|\zeta_{T-\tau}^{\pi_1}-\zeta_{T-\tau}^{\pi_2}\|_{\rm max}& \leq &  e^{\gamma K(\beta^\tau-\beta^T)}\sum\limits_{t=\tau}^{T-1}\|d_t-f_t\|_\infty. THIS IS LAST ZETA VERSION
%\|\Qfunc{t}{T}{\pi_1}-\Qfunc{t}{T}{\pi_2}\|_{\max} &\leq & e^{|\gamma|K(\beta^{t} -\beta^T)}\sum\limits_{\tau=t+1}^{T-1}\|d_\tau-f_\tau\|_\infty.\label{eq_maxnorm_zetaTOfpi1_min_norm_zetaTOfpi2_bound}\\
\|\Qfunc{t}{T}{\pi_1}-\Qfunc{t}{T}{\pi_2}\|_{\max} &\leq & \begin{cases}
        0, &\mbox{\rm  for }t=T-1,\\
        e^{|\gamma|K(\beta^{t} -\beta^T)}\sum\limits_{\tau=t+1}^{T-1}\|d_\tau-f_\tau\|_\infty, &\mbox{\rm for }0\leq t\leq T-2.
\end{cases}\label{eq_maxnorm_zetaTOfpi1_min_norm_zetaTOfpi2_bound}
\end{eqnarray}
\end{lem}
% {\color{blue}
% \begin{lem}\label{lemma_zeta_and_its_differ_bounds1}
% %Let $\pi_1=\{d_t\}_{t\geq 0},\pi_2=\{f_t\}_{t\geq 0}\in \Pi_{\rm MR}$. 
% Fix $T\geq 1$. Let $\pi_{1}=\{d_{t}\}_{t=0}^{\infty}$, $\pi_{2}=\{f_{t}\}_{t=0}^{\infty}\in\Pi_{\rm MR}$. For $0\leq \tau\leq T-1$, denote $E_\tau=\exp\left(\frac{\gamma C(\beta^\tau-\beta^T)}{1-\beta}\right)$. Then, for every  $0\leq \tau\leq T-1$,\vspace{-0.3 cm}
% \begin{eqnarray}
% %\|\zeta_{T}^{\pi_1}\|_{\rm max} & \leq & \exp\left(\frac{\gamma C(1-\beta^T)}{1-\beta}\right) \label{eq_maxnorm_zetaT_bound}\\
% \|\zeta_{T-\tau}^{\pi_1}\|_{\rm max}  \leq  E_\tau,\mbox{ and, } \ \
% %\|\zeta_{T}^{\pi_1}-\zeta_{T}^{\pi_2}\|_{max}& \leq &  \exp\left(\frac{\gamma C(1-\beta^T)}{1-\beta}\right)\sum\limits_{t=0}^{T-1}\|d_t-f_t\|_\infty.\label{eq_maxnorm_zetaTOfpi1_min_norm_zetaTOfpi2_bound}\\
% \|\zeta_{T-\tau}^{\pi_1}-\zeta_{T-\tau}^{\pi_2}\|_{\rm max} \leq   E_\tau\sum\limits_{t=\tau}^{T-1}\|d_t-f_t\|_\infty.\label{eq_maxnorm_zetaTOfpi1_min_norm_zetaTOfpi2_bound1}
% \end{eqnarray}
% \end{lem}

% }

 The inequalities \eqref{eq_maxnorm_zetaT_bound}-\eqref{eq_maxnorm_zetaTOfpi1_min_norm_zetaTOfpi2_bound} along with the backward recursive expression given in Proposition \ref{thm_matrix_formulaton_FinHor_RSCost} combine to yield continuity of the finite-horizon RS cost function. The extension to infinite-horizon RS cost comes from Theorem 1 of \cite{UdaySanjayKavitaNH}, which asserts that the  finite-horizon RS cost converges exponentially  to the infinite-horizon RS cost as the horizon increases. We state here a weaker, asymptotic version of the result by \cite[Thm 1]{UdaySanjayKavitaNH} for convenience. 
% Analogous to Lemma \ref{lem_linCost_finhor_cgsto_infinhor}, the result asserts the convergence of 

\begin{lem}\label{lem_RSCost_finhor_cgsto_infinhor}
\textbf{[Convergence of finite-horizon RS cost:]} Let $x\in\mathcal{S}$, and let  $R:\mathcal{S}\times \mathcal{A}\rightarrow \R$ be an immediate cost function. Then, for every  $\epsilon>0$ there exists $T^*>0$ such that $|\RSCostT{\pi}{\gamma}{R}{T}-\RSCost{\pi}{\gamma}{R}|< \epsilon$ for all $T\geq T^*$ and all $\pi\in\Pi_{\rm MR}$. 
\end{lem}

The following result summarizes the continuity properties of finite- and infinite-horizon RS costs.  
The proof is given in \ref{appendixA}.  

\begin{thm}\label{thm_unif_cont_of_USTOCost_prodtop}\textbf{[Continuity of RS costs:]}
 Let $T\geq 0$. The mapping $\pi \mapsto\RSCostTVec{\pi}{\gamma}{R}{T}$ from the metric space $(\Pi_{\rm MR},\mu)$ to the metric space $(\R^{m},\|\cdot\|_{\infty})$ is Lipschitz continuous with Lipschitz constant  $\delta^{-(T-1)}(1-\delta)^{-1}e^{|\gamma| K(1-\beta^{T})}$. Furthermore, the mapping $\pi\mapsto\RSCostVec{\pi}{\gamma}{R}$ from the metric space $(\Pi_{\rm MR},\mu)$ to the metric space $(\R^{m},\|\cdot\|_{\infty})$ is uniformly continuous. 
\end{thm}

Note that components of a vector-valued Lipschitz or uniformly continuous function are also Lipschitz or uniformly continuous, respectively. As a result, the assertions of theorems \ref{thm_contuty_lin_disc_cost_prod_toplogy} and \ref{thm_unif_cont_of_USTOCost_prodtop} also hold for the respective cost functions for a given initial condition. This is the form in which theorems \ref{thm_contuty_lin_disc_cost_prod_toplogy} and \ref{thm_unif_cont_of_USTOCost_prodtop} will be used in the sequel.

% The next corollary follows easily from theorems \ref{thm_contuty_lin_disc_cost_prod_toplogy} and \ref{thm_unif_cont_of_USTOCost_prodtop} on recalling that components of a vector-valued Lipschitz or uniform continuous function are also Lipschitz or uniform continuous.
% {\color{blue} We can delete the Corollary because it is obvious. We can directly call the theorems.}

% \begin{cor}\label{cor_unif_cont_USLD_prodtop_initial}
% For each initial state $x\in\mathcal{S}$ and $T\geq 0$, each of the maps $\pi\to \RSCostT{\pi}{\gamma}{R}{T}$, $\pi\to \LDCost{\pi}{R}$ and $\pi\to \LDCostT{\pi}{R}{T}$  is Lipschitz continuous and the map $\pi\to\RSCost{\pi}{\gamma}{R}$ is uniformly continuous in  $\pi$ on the metric space $(\Pi_{\rm MR},\mu)$.
% \end{cor}

\subsubsection{Compactness of feasible regions}\label{subsec_cpctsssec}
Theorem \ref{thm_metric_on_Pi_metrizes_ProdTop} already states that $\Pi_{\rm MR}$ is a compact space in the topology induced by the metric $\mu$. Hence compactness of the feasible regions of the problems \eqref{eq_CRSMDP_Prob}, $\CRSMDPFinHorOne$ and $\CRSMDPFinHorTwo$ will follow if it can be shown that these sets are closed subsets of $\Pi_{\rm MR}$. However, each of these sets is defined through non-strict inequalities involving finite- and infinite-horizon classical and RS cost functions. Since all of these functions have been shown to be continuous in the topology on $\Pi_{\rm MR}$ by theorems \ref{thm_contuty_lin_disc_cost_prod_toplogy} and  \ref{thm_unif_cont_of_USTOCost_prodtop}, it follows that the feasible regions are closed subsets of $\Pi_{\rm MR}$, and hence compact. This conclusion is stated as the first assertion of our next lemma. The remaining two  assertions of the lemma pertain to the behaviour of the feasible regions of the problems $\CRSMDPFinHorOne$ and $\CRSMDPFinHorTwo$ for different values of $T$, and will be useful for the proofs of theorems \ref{thm_finhorapprx_frmblw} and \ref{thm_finhorapprx_frmabv}.

% \begin{lem}
% \label{lem_fesblereg_prop}
% The following statements hold. 
% \begin{enumerate}[(i)]
% \item\label{itm_FT_nondec} $\{\FeasibleSetTrunBlw{T}\}_{T=1}^{\infty}$ is non-decreasing sequence of nested sets contained in $\FeasibleSetOrgnal$.
% \item\label{itm_barFT_noninc} $\{\FeasibleSetTrunAbv{T}\}_{T=1}^{\infty}$ is a non-increasing sequence of nested sets which converges to $\FeasibleSetOrgnal$ in the sense $\cap_{T=1}^{\infty}\FeasibleSetTrunAbv{T}=\FeasibleSetOrgnal$.
% \item \label{itm_compact_fsblsets}The set $\FeasibleSetOrgnal$ is a compact subset of $\Pi_{\rm MR}$. For each $T>0$, the sets $\FeasibleSetTrunBlw{T}$ and $\FeasibleSetTrunAbv{T}$ are compact subsets of $\Pi_{\rm MR}$.
% \end{enumerate}
% \end{lem}

\begin{lem}\label{lem_fesblereg_prop}
The following statements hold. 
\begin{enumerate}[(i)]
\item \label{itm_compact_fsblsets} The set $\FeasibleSetOrgnal$ is a compact subset of $\Pi_{\rm MR}$.
\item\label{itm_FT_nondec}  $\{\FeasibleSetTrunBlw{T}\}_{T=1}^{\infty}$ is a non-decreasing sequence of nested compact sets contained in $\FeasibleSetOrgnal$.
\item \label{itm_barFT_noninc} $\{\FeasibleSetTrunAbv{T}\}_{T=1}^{\infty}$ is a non-increasing sequence of nested compact sets satisfying \\ $\cap_{T=1}^{\infty}\FeasibleSetTrunAbv{T}=\FeasibleSetOrgnal$.% in the sense $\cap_{T=1}^{\infty}\FeasibleSetTrunAbv{T}=\FeasibleSetOrgnal$.
%\item \label{itm_compact_fsblsets1}The set $\FeasibleSetOrgnal$ is a compact subset of $\Pi_{\rm MR}$. For each $T>0$, the sets $\FeasibleSetTrunBlw{T}$ and $\FeasibleSetTrunAbv{T}$ are compact subsets of $\Pi_{\rm MR}$.
\end{enumerate}
\end{lem}
\subsubsection{Proof of Theorem \ref{thm_existence_result}}
The results of the preceding few sub-subsections now enable us to provide a formal proof of Theorem 
\ref{thm_existence_result}. 
%%%%%%%%%%%%%%%%%%% Proof of Theorem \ref{thm_existence_result} begins%%%%%%%%%%%%%%%%%%%%%%%%%%%%%%%%%%%%%%%%%%%%%%%%%%

%\textit{Proof of Theorem \ref{thm_existence_result}:}
% First, recall from Remark \ref{rem_us} that, if problem $\CRSMDPFinHorTwo$ possesses a solution for some $T>0$, then it also possesses a US solution.  Next, 
To prove the first assertion in Theorem \ref{thm_existence_result}, suppose $\mathcal{F}$ is nonempty. 
By Lemma \ref{lem_fesblereg_prop}.\ref{itm_barFT_noninc}, $\FeasibleSetTrunAbv{T}$ is nonempty for all $T$. Theorem \ref{thm_unif_cont_of_USTOCost_prodtop}
% Corollary \ref{cor_unif_cont_USLD_prodtop_initial} 
implies that the objective functions of the problems \eqref{eq_CRSMDP_Prob} and $\CRSMDPFinHorTwo$ are continuous functions on the metric space of policies $(\Pi_{\rm MR},\mu)$. Statements  \ref{itm_compact_fsblsets}   and \ref{itm_barFT_noninc} of Lemma \ref{lem_fesblereg_prop} guarantee  that the respective feasible regions $\FeasibleSetOrgnal$ and $\FeasibleSetTrunAbv{T}$ are compact subsets of $(\Pi_{\rm MR},\mu)$. The first assertion now follows from the well-known fact that a continuous function achieves its infimum over a non-empty compact set. 

%It follows that  the problems \eqref{eq_CRSMDP_Prob} and  $\CRSMDPFinHorTwo$ possess solutions. 

To prove the second assertion, let $T^*>0$ be such that  $\FeasibleSetTrunBlw{T^*}$ is nonempty. Statement \ref{itm_FT_nondec} of Lemma \ref{lem_fesblereg_prop}  implies that $\FeasibleSetOrgnal$ and $\FeasibleSetTrunBlw{T}$ are nonempty for all $T\geq T^{*}$. The rest of the assertion now follows by arguing as in the previous paragraph. \eop
%rguments under the assumption that the feasible region $\FeasibleSetTrunBlw{T}$ of $\CRSMDPFinHorOne$ is non-empty.  

\subsection{Proof Outline of Theorem \ref{thm_finhorapprx_frmblw}}\label{thm2ssec}
The starting point for the proof of Theorem \ref{thm_finhorapprx_frmblw} is Lemma \ref{lem_RSCost_finhor_cgsto_infinhor} 
%and Lemma \ref{lem_linCost_finhor_cgsto_infinhor}
which  states that, for a given initial state, the finite-horizon RS and standard discounted costs converge to corresponding infinite-horizon costs uniformly in the  policy  as the horizon $T\to\infty$. 
In view of Lemma \ref{lem_RSCost_finhor_cgsto_infinhor} and Theorem \ref{thm_existence_result}, to prove
%\ref{lem_linCost_finhor_cgsto_infinhor}, 
Theorem \ref{thm_finhorapprx_frmblw}, it suffices to show that the optimal value of  $\CRSMDPFinHorOne$ converges to that of \eqref{eq_CRSMDP_Prob} as $T\rightarrow \infty$, that is,
\begin{equation}
\lim\limits_{T\to\infty}\left(\inf\limits_{\pi\in\FeasibleSetTrunBlw{T}}\RSCostT{\pi}{\gamma}{R}{T}\right) = \inf\limits_{\pi\in\FeasibleSetOrgnal}\RSCost{\pi}{\gamma}{R}.\label{eq_PTminus_cgsto_P}
\end{equation}
Indeed, if \eqref{eq_PTminus_cgsto_P} holds, then the proof of Theorem \ref{thm_finhorapprx_frmblw} may be completed as follows: 
\renewcommand\theenumi{(\alph{enumi})}
\begin{enumerate}
\item  Suppose \eqref{eq_PTminus_cgsto_P} holds. Choose $\epsilon>0$. By Lemma \ref{lem_RSCost_finhor_cgsto_infinhor}, there exists $T^*>0$ such that $|\RSCostT{\pi}{\gamma}{R}{T}-\RSCost{\pi}{\gamma}{R}|< \epsilon/2$ for all $T\geq T^*$ and all $\pi\in\Pi_{\rm MR}$. 
\item By Theorem \ref{thm_existence_result}, the  right hand side of \eqref{eq_PTminus_cgsto_P} is well defined, which implies that $\FeasibleSetTrunBlw{T}\neq\varnothing$ eventually. 
\item We may choose $T^{*}$ in (a) above to be further larger, if required,  so that  every $T>T^{*}$ also satisfies $\FeasibleSetTrunBlw{T}\neq\varnothing$ and  
$$\left|\inf_{\pi\in\FeasibleSetTrunBlw{T}} \RSCostT{\pi}{\gamma}{R}{T}-\inf_{\pi\in\FeasibleSetOrgnal}\RSCost{\pi}{\gamma}{R}\right|<\epsilon/2.$$
% $$\left|\inf_{\pi\in\FeasibleSetTrunBlw{T}} \RSCostT{\pi}{\gamma}{R}{T}-\inf_{\pi\in\FeasibleSetOrgnal}\RSCost{\pi}{\gamma}{R}\right|<\epsilon/2.$$
\item Again by Theorem \ref{thm_existence_result}, for  every $T>T^*$, there exists  a policy $\eta\in\Pi_{\rm MR}$  that solves $\CRSMDPFinHorOne$, that is, a $\eta\in\FeasibleSetTrunBlw{T}$ such that 
%\begin{eqnarray}
$\RSCostT{\eta}{\gamma}{R}{T} = \inf_{\pi\in\FeasibleSetTrunBlw{T}}\RSCostT{\pi}{\gamma}{R}{T}$.
% and  $|\RSCostT{\eta}{\gamma}{R}{T}-\inf_{\pi\in\FeasibleSetOrgnal}\RSCost{\pi}{\gamma}{R}|<\frac{\epsilon}{2}$.
%\label{eq_JetaT_min_OptmJpi_leqeps}
%\end{eqnarray} 
\item  Finally, the inequalities in (a) and (c)  together show that every  policy $\eta\in\FeasibleSetTrunBlw{T}$ that solves problem $\CRSMDPFinHorOne$ for a given $T>T^{*}$ also satisfies
$\RSCost{\eta}{\gamma}{R} - \inf_{\pi\in\FeasibleSetOrgnal}\RSCost{\pi}{\gamma}{R} = (\RSCost{\eta}{\gamma}{R} -\RSCostT{\eta}{\gamma}{R}{T}) + ( \RSCostT{\eta}{\gamma}{R}{T} - \inf_{\pi\in\FeasibleSetOrgnal}\RSCost{\pi}{\gamma}{R})
<\epsilon,$ 
%\begin{eqnarray}
% $$\RSCost{\eta}{\gamma}{R} - \inf_{\pi\in\FeasibleSetOrgnal}\RSCost{\pi}{\gamma}{R} = (\RSCost{\eta}{\gamma}{R} -\RSCostT{\eta}{\gamma}{R}{T}) + ( \RSCostT{\eta}{\gamma}{R}{T} - \inf_{\pi\in\FeasibleSetOrgnal}\RSCost{\pi}{\gamma}{R})
% <\epsilon,$$ 
that is, $\eta$ satisfies \eqref{eq_truncapprox}.
%\end{eqnarray}
% Furthermore, arguing as in Remark \ref{rem_us} shows that we may choose the policy $\eta$ to be ultimately stationary.
\end{enumerate}
 Thus, the proof of Theorem \ref{thm_finhorapprx_frmblw} is achieved if \eqref{eq_PTminus_cgsto_P} is proved. Hence we next outline the steps involved in proving \eqref{eq_PTminus_cgsto_P}. 

 Lemma \ref{lem_RSCost_finhor_cgsto_infinhor} asserts that the objective function $\RSCostT{\pi}{\gamma}{R}{T}$ of the problem $\CRSMDPFinHorOne$ appearing on the left hand side of \eqref{eq_PTminus_cgsto_P} uniformly approximates  the objective function $\RSCost{\pi}{\gamma}{R}$ of the problem \eqref{eq_CRSMDP_Prob} appearing on the right hand side of \eqref{eq_PTminus_cgsto_P} as $T\to\infty$. Hence one intuitively expects \eqref{eq_PTminus_cgsto_P} to hold if the feasible region $\FeasibleSetTrunBlw{T}$ of problem  $\CRSMDPFinHorOne$ converges to the feasible region $\FeasibleSetOrgnal$ of the problem \eqref{eq_CRSMDP_Prob}  in some suitable sense as $T\to \infty$.

We already know from Lemma \ref{lem_fesblereg_prop} that  $\{\FeasibleSetTrunBlw{T}\}_{T=1}^{\infty}$ is an non-decreasing sequence of nested subsets of $\FeasibleSetOrgnal$. The following result shows that the additional condition under which \eqref{eq_PTminus_cgsto_P} holds is the equality $\cl\left(\cup_{T=1}^{\infty}\FeasibleSetTrunBlw{T}\right)= \FeasibleSetOrgnal$, where $\cl$ denotes topological closure. The result is stated below in the more general context of topological spaces, as it may be of independent interest.

\begin{prop}\label{prop_generalised_convgce_offunctions_over_seqofsets}
Let $\mathcal{U}$ be a topological space. Let $\{f_T\}_{T=1}^{\infty}$ be a sequence of real-valued functions on $\mathcal{U}$ converging uniformly to a continuous function $f:\mathcal{U}\to\R$. Let $\{\mathcal{V}_T\}_{T=1}^{\infty}$be a non-decreasing sequence of nested subsets of $\mathcal{U}$, and let  $\mathcal{V}:= \cl\left(\cup_{T=1}^\infty\mathcal{V}_T\right)$.  Suppose $f$ is bounded below on the set $\mathcal{V}$. Then,
\begin{equation}
\lim\limits_{T\to\infty}\left(\inf\limits_{y\in\mathcal{V}_T}f_T(y)\right) = \inf\limits_{y\in\mathcal{V}} f(y).\label{eq_limntoinf_infofgnoverFn_eqls_infofgoverF}
\end{equation} 
\end{prop}

In light of the above proposition, the limit in \eqref{eq_PTminus_cgsto_P} holds if 
\begin{equation}
\FeasibleSetOrgnal=\cl\left(\cup_{T=1}^{\infty}\FeasibleSetTrunBlw{T}\right)
\label{eq_cls_un_ovrT_intsc_ovi_FTi_eqls_intsc_ovi_Fi_1}
\end{equation}
holds.
Recall that, the sets $\{\FeasibleSetTrunBlw{T}\}_{T=1}^{\infty}$ are intersections of sublevel sets of finite-horizon standard discounted and RS cost functions, while the set $\FeasibleSetOrgnal$ is the intersection of sublevel sets of finite- and infinite-horizon standard discounted and RS cost functions. Moreover, by lemmas \ref{lem_linCost_finhor_cgsto_infinhor} and \ref{lem_RSCost_finhor_cgsto_infinhor}, the functions defining the subsets $\{\FeasibleSetTrunBlw{T}\}_{T=1}^{\infty}$ converge to the functions defining the set $\FeasibleSetOrgnal$. These properties lead one to expect \eqref{eq_cls_un_ovrT_intsc_ovi_FTi_eqls_intsc_ovi_Fi_1} to hold.  However, as a counterexample below shows, \eqref{eq_cls_un_ovrT_intsc_ovi_FTi_eqls_intsc_ovi_Fi_1} may not hold in general. 
% these features are themselves not sufficient for the equality \eqref{eq_cls_un_ovrT_intsc_ovi_FTi_eqls_intsc_ovi_Fi_1} to hold.
The following result provides the additional condition under which the desired equality above holds.  The additional condition requires the maximum constraint violation function to not have 0 as a local minimum value, 
% An additional condition requiring a certain function to not have $0$ as its local minimum value,
and is the source of the  similar-looking  condition appearing in Theorem \ref{thm_finhorapprx_frmblw}.  The result is again stated in the more general context of topological spaces, as it may be of independent interest when stated as such. 
\begin{prop}\label{prop_cl_uninovrT_intscovri_FTi_eql_intscovrFi}
Let $\mathcal{U}$ be a topological space, and let $N>0$. For each $i\in\{1,\ldots,N\}$, let $\{g_{T,i}\}_{T=1}^{\infty}$ be a sequence of real-valued functions on $\mathcal{U}$ converging pointwise to a continuous function  $g_i:\mathcal{U}\rightarrow \R$ as $T\rightarrow\infty$,
%For each \sout{$i=1,\ldots,N$} {\color{red} $i\in\{1,\ldots,N\}$}, 
and let $\{B_{T,i}\}_{T=1}^{\infty}$ be a sequence of real numbers converging to $B_i$ as $T\rightarrow \infty$. For each $T$ and $i$,  denote $\mathcal{G}_{T,i}=\{z\in\mathcal{U}:g_{T,i}(z)\leq B_{T,i}\}$,  $\mathcal{H}_{T}=\cap_{i=1}^{N}\mathcal{G}_{T,i}$,  $\mathcal{G}_i=\{z:g_{i}(z)\leq B_{i}\}$ and  $\mathcal{H}=\cap_{i=1}^{N}\mathcal{G}_{i}.$
% \begin{equation}
% \mathcal{G}_{T,i}=\{z\in\mathcal{U}:g_{T,i}(z)\leq B_{T,i}\},\ \mathcal{H}_{T}=\bigcap\limits_{i=1}^{N}\mathcal{G}_{T,i},\ 
% \mathcal{G}_i=\{z:g_{i}(z)\leq B_{i}\},\ \mathcal{H}=\bigcap\limits_{i=1}^{N}\mathcal{G}_{i}.\nonumber
% \end{equation}
Assume that $\mathcal{G}_{T,i}\subseteq \mathcal{G}_{i}$ for every $i$ and $T$. Define the map $h:\mathcal{U}\to\R$ by $h(z)=\max_{i\in\{1,\ldots,N\}}\{g_i(z)-B_i\}$. 
If $0$ is not a local minimum value of $h$, then 
$\cl\left(\cup_{T=1}^{\infty}\mathcal{H}_{T}\right)= \mathcal{H}$. 
%  \begin{equation} {\color{blue}
% \cl\left(\cup_{T=1}^{\infty}\mathcal{H}_{T}\right)= \mathcal{H}. }
% \label{eq_cls_un_ovrT_intsc_ovi_FTi_eqls_intsc_ovi_Fi} 
% \end{equation} 
\end{prop}

While the condition that $0$ should not be a local minimum value of the function $h$ in Proposition \ref{prop_cl_uninovrT_intscovri_FTi_eql_intscovrFi} appears cumbersome to check, it is easy to construct an example where the condition is not satisfied, and the conclusion of the proposition fails to hold. Indeed, let $\mathcal{U}=\R$,  $N=2$, $B_{1}=B_{2}=0$ and $g_{1}(x)=-g_{2}(x)=x$ for all $x\in\R$. For each $T$, let  $g_{T,1}(x)=g_{1}(x)$, $g_{T,2}(x)=g_{2}(x)$, and $B_{T,1}=B_{T,2}=-T^{-1}$. The function $h$ of Proposition \ref{prop_cl_uninovrT_intscovri_FTi_eql_intscovrFi} is then given by $h(x)=|x|$, and has $0$ as a local minimum value. It is easy to see that $\mathcal{G}_{T,1}=(-\infty,-T^{-1}]$, $\mathcal{G}_{1}=(-\infty,0]$, $\mathcal{G}_{T,2}=[T^{-1},\infty)$ and  $\mathcal{G}_{2}=[0,\infty)$. As a result,  $\mathcal{H}=\{0\}$, while  $\mathcal{H}_{T}=\varnothing $ for all $T$, showing that 
%\eqref{eq_cls_un_ovrT_intsc_ovi_FTi_eqls_intsc_ovi_Fi} 
 the assertion of Proposition \ref{prop_cl_uninovrT_intscovri_FTi_eql_intscovrFi} fails to hold. 

The obstruction posed by the failure of the local minimum condition goes beyond the purely topological setting of Proposition \ref{prop_cl_uninovrT_intscovri_FTi_eql_intscovrFi}.
%The failure of the local minimum condition poses an obstruction not only in the topological setting of Proposition \ref{prop_cl_uninovrT_intscovri_FTi_eql_intscovrFi}, but also in the CRSMDP setting of Theorem  \ref{thm_finhorapprx_frmblw}. 
Indeed, we present a more elaborate counter-example in \ref{appendixD} in a CRSMDP setting to show that the conclusion of Theorem 
\ref{thm_finhorapprx_frmblw} can fail to hold if the local minimum condition assumed therein is violated. 

\subsubsection{Proof of Theorem \ref{thm_finhorapprx_frmblw}}
 Problem \eqref{eq_CRSMDP_Prob} involves  a total of $\widetilde{M}:=M+\hat{M}+\bar{M}+\check{M}$ constraints.
%To apply Proposition \ref{prop_cl_uninovrT_intscovri_FTi_eql_intscovrFi}, 
For each  $T\in\{1,2,\ldots,\}$ and $i\in\{1,\ldots, \widetilde{M}\}$, define
% \begin{eqnarray}
% g_{T,i}(\pi)&:=&
% 			\LDCostT{\pi}{C_i}{T},  \  i=1,2,\ldots,M, \label{eq_gTi_finhorz} \\
% &:=&			\RSCostT{\pi}{\gamma}{\hat{C}_j}{T}, \  i=M+j,\  j=1,\ldots,\hat{M},
% \end{eqnarray}
% \begin{eqnarray*}
%     &:=&        \LDCostT{\pi}{\bar{C}_k}{\bar{T}_k}, \ i=M+\hat{M}+k,\ k=1,\ldots\bar{M},\\
%     &:=&         \RSCostT{\pi}{\gamma}{\check{C}_l}{\check{T}_l},  i=M+\hat{M}+\bar{M}+l,\ l=1,\ldots,\check{M},
% \end{eqnarray*}
% \begin{equation}\label{eq_gTi_finhorz}
% g_{T,i}(\pi):=
% \begin{cases}
% 			\LDCostT{\pi}{C_i}{T},  & i=1,2,\ldots,M,\\
% 			\RSCostT{\pi}{\gamma}{\hat{C}_j}{T}, & i=M+j,\  j=1,\ldots,\hat{M},\\
%             \LDCostT{\pi}{\bar{C}_k}{\bar{T}_k}, & i=M+\hat{M}+k,\ k=1,\ldots\bar{M},\\
%              \RSCostT{\pi}{\gamma}{\check{C}_l}{\check{T}_l}, & i=M+\hat{M}+\bar{M}+l,\ l=1,\ldots,\check{M},\\
% 		 \end{cases}
% \end{equation}
% \begin{equation}\label{eq_dfn_BTi}
% B_{T,i}:=
% \begin{cases}
% 		b_i-K\beta^T,  & i=1,2,\ldots,M,\\
% 		\hat{b}_j/K_{T}, & i=M+j,\  j=1,\ldots,\hat{M},\\
% 		\bar{b}_k, & i=M+\hat{M}+k,\ k=1,\ldots\bar{M},\\
% 		\check{b}_l,& i=M+\hat{M}+\bar{M}+l,\ l=1,\ldots,\check{M}.
% \end{cases}
% \end{equation}
\begin{eqnarray}
g_{T,i}(\pi)&:=&
\begin{cases}
			\LDCostT{\pi}{C_i}{T},  & i=1,2,\ldots,M,\\
			\RSCostT{\pi}{\gamma}{\hat{C}_j}{T}, & i=M+j,\  j=1,\ldots,\hat{M},\\
            \LDCostT{\pi}{\bar{C}_k}{\bar{T}_k}, & i=M+\hat{M}+k,\ k=1,\ldots\bar{M},\\
             \RSCostT{\pi}{\gamma}{\check{C}_l}{\check{T}_l}, & i=M+\hat{M}+\bar{M}+l,\ l=1,\ldots,\check{M},\\
		 \end{cases} \label{eq_gTi_finhorz} \\
B_{T,i}&:=&
\begin{cases}
		b_i-K\beta^T,  & i=1,2,\ldots,M,\\
		\hat{b}_j/K_{T}, & i=M+j,\  j=1,\ldots,\hat{M},\\
		\bar{b}_k, & i=M+\hat{M}+k,\ k=1,\ldots\bar{M},\\
		\check{b}_l,& i=M+\hat{M}+\bar{M}+l,\ l=1,\ldots,\check{M}.
\end{cases}\label{eq_dfn_BTi}
\end{eqnarray}
Lemmas \ref{lem_linCost_finhor_cgsto_infinhor} and \ref{lem_RSCost_finhor_cgsto_infinhor} show that, for each $i\in\{1,\ldots, \widetilde{M}\}$, the sequence of real-valued functions $\{g_{T,i}\}_{T=1}^{\infty}$ converges uniformly on $\Pi_{\rm MR}$ to the function $g_{i}:\Pi_{\rm MR}\rightarrow \R$ defined by  
\begin{equation}\label{eq_gi_infhorz}
g_{i}(\pi):=
\begin{cases}
			\LDCost{\pi}{C_i},  & i=1,2,\ldots,M,\\
			\RSCost{\pi}{\gamma}{\hat{C}_j}, & i=M+j,\  j=1,\ldots,\hat{M},\\
            \LDCostT{\pi}{\bar{C}_k}{\bar{T}_k}, & i=M+\hat{M}+k,\ k=1,\ldots\bar{M},\\
             \RSCostT{\pi}{\gamma}{\check{C}_l}{\check{T}_l}, & i=M+\hat{M}+\bar{M}+l,\ l=1,\ldots,\check{M}. \\
		 \end{cases}
\end{equation}
Also, for each $i\in\{1,\ldots, \widetilde{M}\}$ the sequence $\{B_{T,i}\}_{T=1}^{\infty}$ converges to $B_{i}$, where \begin{eqnarray}\label{eq_dfn_Binfi}
B_{i}:=
\begin{cases}
    	b_i,  & i=1,2,\ldots,M,\\
		\hat{b}_j, & i=M+j,\ j=1,\ldots,\hat{M},\\
		\bar{b}_k, & j=M+\hat{M}+k,\ k=1,\ldots,\bar{M},\\
		\check{b}_l,& i=M+\hat{M}+\bar{M}+l,\ l=1,\ldots,\check{M}.
\end{cases}
\end{eqnarray}

Next, for each $T\in\{1,2,\ldots,\}$ and $i\in\{1,\ldots, \widetilde{M}\}$, define the sets $\mathcal{G}_{T,i}:=\{\pi\in\Pi_{\rm MR}:g_{T,i}(\pi)\leq B_{T,i}\}$ and $\mathcal{G}_{i}:=\{\pi\in\Pi_{\rm MR}:g_{i}(\pi)\leq B_{i}\}$. Note that $\cap_{i=1}^{\widetilde{M}}\mathcal{G}_{T,i}= \FeasibleSetTrunBlw{T}$ for each  $T\in\{1,2,\ldots,\}$, while $\cap_{i=1}^{\widetilde{M}}\mathcal{G}_{i}=\FeasibleSetOrgnal$. Moreover, claims 1 and 2  established in the proof of Lemma \ref{lem_fesblereg_prop} imply that $\mathcal{G}_{T,i}\subseteq \mathcal{G}_{i}$ for all   $T\in\{1,2,\ldots,\}$ and $i\in\{1,\ldots, \widetilde{M}\}$. We also observe that the function $h:\Pi_{\rm MR}\to\R$ defined in \eqref{eq_defn_of_h} is equivalently given by   $h(\pi)= \max_{i\in\{1,2,\ldots,\widetilde{M}\}}\{g_i(\pi)-B_{i}\}$. Proposition \ref{prop_cl_uninovrT_intscovri_FTi_eql_intscovrFi} now applies with $\mathcal{U}:=\Pi_{\rm MR}$, $N:=\widetilde{M}$, $\mathcal{H}:=\FeasibleSetOrgnal$ and $\mathcal{H}_{T}:=\FeasibleSetTrunBlw{T}$ for each $T\in\{1,2,\ldots\}$, and allows us to conclude that \eqref{eq_cls_un_ovrT_intsc_ovi_FTi_eqls_intsc_ovi_Fi_1} holds.

Next, recall from Lemma \ref{lem_RSCost_finhor_cgsto_infinhor} that the sequence $\{\RSCostT{\pi}{\gamma}{R}{T}\}_{T=1}^{\infty}$  converges uniformly to  $\RSCost{\pi}{\gamma}{R}$ on $\Pi_{\rm MR}$. Moreover, the continuous function  $\pi\mapsto \RSCost{\pi}{\gamma}{R}$ is bounded below on the compact metric space $\Pi_{\rm MR}$.   Since  \eqref{eq_cls_un_ovrT_intsc_ovi_FTi_eqls_intsc_ovi_Fi_1} holds, Proposition \ref{prop_generalised_convgce_offunctions_over_seqofsets} applies by letting $\mathcal{U}=\Pi_{\rm MR}$, $f_{T}=\RSCostT{\pi}{\gamma}{R}{T}$ and $\mathcal{V}_{T}=\FeasibleSetTrunBlw{T}$ for each $T$, $f=\RSCost{\pi}{\gamma}{R}$ and $\mathcal{V}=\FeasibleSetOrgnal$. We conclude from Proposition \ref{prop_generalised_convgce_offunctions_over_seqofsets} that \eqref{eq_PTminus_cgsto_P} holds. The rest of the proof is completed by using the arguments given in the paragraph following \eqref{eq_PTminus_cgsto_P}. This completes the proof of Theorem \ref{thm_finhorapprx_frmblw}.
\eop

\subsection{Proof Outline of Theorem \ref{thm_finhorapprx_frmabv}} \label{sec_thm3ssec}
% The proof of Theorem \ref{thm_finhorapprx_frmabv} uses the following properties of the components of the problems  \eqref{eq_CRSMDP_Prob} and $\CRSMDPFinHorTwo$:
% \begin{itemize}
% \item continuity of the objective functions of both the problems \eqref{eq_CRSMDP_Prob} and $\CRSMDPFinHorTwo$,
% \item the nested property, its convergence and compactness of the feasible sets  $\{\FeasibleSetTrunAbv{T}\}_{T=1}^{\infty}$ of the problem $\CRSMDPFinHorTwo$ (see Lemma \ref{lem_fesblereg_prop}) and
% \item compactness of the original feasible region $\FeasibleSetOrgnal$.\end{itemize}
%  As discussed in the section \ref{sec_main_results}, the approximation of the problem \eqref{eq_CRSMDP_Prob} given in the theorem \ref{thm_finhorapprx_frmblw} requires a local minimum condition to be satisfied by the function $h$ defined by \eqref{eq_defn_of_h}, which is not always easy to check. In view of this, we provide an alternative approximation, when this condition is not satisfied or difficult to check, we give theorem \ref{thm_finhorapprx_frmabv} where both the feasibility and the objective are approximated to be as close as possible to that of original problem \eqref{eq_CRSMDP_Prob}. The closeness of the feasible region is measured by $\epsilon$-feasibility defined in section \ref{sec_main_results}.

The proof of Theorem \ref{thm_finhorapprx_frmabv} follows a pattern similar to that of Theorem \ref{thm_finhorapprx_frmblw}. The main step in the proof is to show that the value of the problem $\CRSMDPFinHorTwo$ approximates the value of \eqref{eq_CRSMDP_Prob} for large $T$, that is, 
\begin{equation}
\lim\limits_{T\to\infty}\left(\inf\limits_{\pi\in\FeasibleSetTrunAbv{T}}\RSCostT{\pi}{\gamma}{R}{T}\right) = \inf\limits_{\pi\in\FeasibleSetOrgnal}\RSCost{\pi}{\gamma}{R}.\label{eq_PTplus_cgsto_P}
\end{equation}
Theorem \ref{thm_finhorapprx_frmblw} then follows by arguing as in the paragraph following \eqref{eq_PTminus_cgsto_P}. 

Intuitively, one expects \eqref{eq_PTplus_cgsto_P} to hold as Lemma \ref{lem_RSCost_finhor_cgsto_infinhor} guarantees that $ \RSCostT{\pi}{\gamma}{R}{T}$ converges uniformly in $\pi$  to $\RSCost{\pi}{\gamma}{R} $ with increasing $T$, while statement \ref{itm_barFT_noninc} of Lemma \ref{lem_fesblereg_prop} asserts that the sequence of sets $\{\FeasibleSetTrunAbv{T}\}_{T=1}^{\infty}$ decrease to $ \FeasibleSetOrgnal$. The next result formalizes this intuition. The result is stated in the more general topological setting, as it might be of independent interest.

% Similar to Theorem \ref{thm_finhorapprx_frmblw}, to prove Theorem \ref{thm_finhorapprx_frmabv}, it suffices to show that the problem $\CRSMDPFinHorTwo$ approximates the problem \eqref{eq_CRSMDP_Prob} in value for a sufficiently large $T$, that is,

% The result \eqref{eq_PTplus_cgsto_P} implies there exists $\eta\in \FeasibleSetTrunAbv{T}$ such that \eqref{eq_truncapprox} is satisfied and convergence of $\FeasibleSetTrunAbv{T}$ to $\FeasibleSetOrgnal$ helps proving the $\epsilon$-feasibility. The rest of the proof is similar to that of Theorem \ref{thm_finhorapprx_frmblw}. 

% The following results helps in ascertaining the result \eqref{eq_PTplus_cgsto_P}. 

% \begin{prop}\label{prop_lim_of_inf_overGT_eql_inf_overG}
% Let $\mathcal{U}$ be a topological space, and let $g:\mathcal{U}\to\R$ be a continuous function. Let $\{\mathcal{G}_T\}_{T=0}^{\infty}$ be a non-increasing sequence of nested compact subsets of $\mathcal{U}$ converging to a subset $\mathcal{G}$ of $\mathcal{U}$ in the sense $\bigcap_{T=0}^{\infty}\mathcal{G}_{T}=\mathcal{G}$. Then,
% \begin{equation}
% \lim\limits_{T\to\infty}\inf\limits_{y\in\mathcal{G}_T}g(y) = \inf\limits_{y\in\mathcal{G}}g(y).\label{eq_limit_GT_infm_inf_G}
% \end{equation}
% \end{prop}

\begin{prop}\label{prop_lim_of_inf_overGT_eql_inf_overG}
Let $\mathcal{U}$ be a topological space, and let $g:\mathcal{U}\to\R$ be a continuous function. Let $\{\mathcal{G}_T\}_{T=1}^{\infty}$ be a non-increasing sequence of nested compact subsets of $\mathcal{U}$ converging to a subset $\mathcal{G}$ of $\mathcal{U}$ in the sense $\bigcap_{T=1}^{\infty}\mathcal{G}_{T}=\mathcal{G}$. Let $\{g_{T}\}_{T=1}^{\infty}$ be a sequence of real-valued functions on $\mathcal{U}$ converging to $g$ uniformly on  $\mathcal{G}_{1}$. Then,
\begin{equation}
\lim\limits_{T\to\infty}\left(\inf\limits_{y\in\mathcal{G}_T}g_{T}(y)\right) = \inf\limits_{y\in\mathcal{G}}g(y).\label{eq_limit_GT_infm_inf_G}
\end{equation}
\end{prop}

\subsubsection{Proof of Theorem \ref{thm_finhorapprx_frmabv}}

%We give here full proof of the Theorem \ref{thm_finhorapprx_frmabv}.

% the result \eqref{eq_limit_GT_infm_inf_G} holds, that is,
% \begin{equation}
% \lim\limits_{T\to\infty}\left(\inf\limits_{\pi\in\FeasibleSetTrunAbv{T}}\RSCostT{\pi}{\gamma}{R}{T}\right) = \inf\limits_{\pi\in\FeasibleSetOrgnal}\RSCost{\pi}{\gamma}{R}\label{eq_limit_FTpls_infm_inf_F}
% \end{equation}
% holds. 

To begin,  choose $\epsilon>0$. It is easy to see that, for each $i$ and $j$,  the constants $b_i + K\beta^T $ and $\hat{b}_j K_{T}$ defining the bounds appearing in the definitions \eqref{eq_defn_LDTrunAbv} and \eqref{eq_defn_RSTrunAbv} of $\FeasibleSetLDTruncAbv{T}$ and $\FeasibleSetRSTruncAbv{T}$ converge to $b_i$ and $\hat{b}_j$, respectively, from above as $T\to \infty$. It therefore follows from lemmas \ref{lem_linCost_finhor_cgsto_infinhor} and \ref{lem_RSCost_finhor_cgsto_infinhor} that there exists $T_{1}$ such that, for every $T>T_{1}$, every  policy in the feasible set  $\FeasibleSetTrunAbv{T}$ of problem $\CRSMDPFinHorTwo$ is $\epsilon$-feasible for the problem \eqref{eq_CRSMDP_Prob}.

Next, we apply Proposition \ref{prop_lim_of_inf_overGT_eql_inf_overG} by taking  $\mathcal{U}$, $g$, $\{\mathcal{G}_{T}\}_{T=1}^{\infty}$, $\mathcal{G}$ and $\{g_{T}\}_{T=1}^{\infty}$  in that proposition  to be $\Pi_{\rm MR}$, $\RSCost{\pi}{\gamma}{R}$,  $\{\FeasibleSetTrunAbv{T}\}_{T=1}^{\infty}$, $\FeasibleSetOrgnal$, and $\{\RSCostT{\pi}{\gamma}{R}{T}\}_{T=1}^{\infty}$, respectively. Lemma \ref{lem_fesblereg_prop}, Lemma \ref{lem_RSCost_finhor_cgsto_infinhor} and Theorem \ref{thm_unif_cont_of_USTOCost_prodtop} imply that all the hypotheses of Proposition  \ref{prop_lim_of_inf_overGT_eql_inf_overG} are satisfied,  and hence \eqref{eq_PTplus_cgsto_P} follows from \eqref{eq_limit_GT_infm_inf_G}. 

Using arguments similar to those laid out in the paragraph following  \eqref{eq_PTminus_cgsto_P}, we can conclude from \eqref{eq_PTplus_cgsto_P} that there exists $T^{*}\geq T_{1}$ such that, for every $T>T^{*}$, there exists a solution  $\eta\in\FeasibleSetTrunAbv{T}$ of $\CRSMDPFinHorTwo$ such that \eqref{eq_truncapprox2} holds. \eop 
% Furthermore, arguing as in Remark \ref{rem_us} shows that we may choose the policy $\eta$ to be ultimately stationary. This completes the proof. \eop

% Moreover, it is easy to observe that the  constants defining the bounds appearing in the definition of $\FeasibleSetLDTruncAbv{T}$ and $\FeasibleSetRSTruncAbv{T}$, respectively,  $b_i + K\beta^T $ and $\hat{b}_j K_{T}$ converge to respectively $b_i$ and $\hat{b}_j$ from above as $T\to \infty$, for each $i$ and $j$. Thus, the feasible set of problem $\CRSMDPFinHorTwo$, $\FeasibleSetTrunAbv{T}$ is such that for $T$ sufficiently large, the policies in  $\FeasibleSetTrunAbv{T}$ are $\epsilon$-feasible (Definition \ref{def_eps_feasb}) to the problem \eqref{eq_CRSMDP_Prob}.

% With \eqref{eq_limit_FTpls_infm_inf_F} and following the similar arguments given in paragraph immediately after \eqref{eq_PTminus_cgsto_P} we can show that there exists $T^*>0$ and a policy $\eta\in\FeasibleSetTrunAbv{T^*}$ which is $\epsilon$-feasible such that satisfies. 
% \eop

% } 
% \end{proof}
\section{Solution of Finite-Horizon CRSMDPs}\label{sec_algo}
Recall that, in light of Remark \ref{rem_us1},  theorems \ref{thm_finhorapprx_frmblw} and \ref{thm_finhorapprx_frmabv} provide a means of constructing US policies as approximate solutions to the infinite-horizon CRSMDP \eqref{eq_CRSMDP_Prob} from solutions of the finite-horizon CRSMDPs $\CRSMDPFinHorOne$ or $\CRSMDPFinHorTwo$, respectively. 
In this section, we provide a LP-based approach for computing the solution to a finite-horizon CRSMDP such as $\CRSMDPFinHorOne$ or $\CRSMDPFinHorTwo$. The approach involves first expressing finite-horizon RS cost functions as standard discounted cost functions of a terminal-cost MDP defined on an augmented state space. 

To illustrate the approach, choose $T>0$, and consider the RS objective cost function $\RSCostT{\pi}{\gamma}{R}{T}$ defined in  \eqref{eq_RSCost_FinandinfHor}. For each  $t\in\{0,1,\ldots,T\}$, define the random variable 
$\Psi^{\rm o}_{t}:=\exp(\gamma\sum_{j=0}^{t-1}\beta^j R(X_j,A_j))$. Note that, for each $t$, the random variable $\Psi^{\rm o}_{t}$ takes at most $(mn)^{t}$ values.  Consequently, the augmented process $\{Z_{t}\}_{t=0}^{T}$ defined by $Z_t := (X_t,\Psi^o_t)$, $t\in\{0,\ldots,T-1\}$, takes at most $m(mn)^{t}$ values at each time $t$, and is the state process of an augmented time-dependent MDP. The probability that the augmented MDP transitions from state  $z=(s,\psi)$ at time $t$ to $z^{\prime}=(s^{\prime},\psi^{\prime})$ at time $t+1$  is given by 
$p^{\rm A}_{t+1}(z^\prime|z,a)=1_{\{\psi^\prime=\psi\exp(\gamma \beta^t R(s,a))\}} p(s^\prime|s,a), 
$
where the $1_{S}$ denotes the indicator function of the set $S$, and the transition probabilities $p(\cdot|\cdot,\cdot)$  of the original MDP are as defined in section \ref{sec_model_framework}.
Next, note that $\RSCostT{\pi}{\gamma}{R}{T} =\E^\pi_x[\Psi^o_{T}]$. 
The RS cost function $\RSCostT{\pi}{\gamma}{R}{T}$ can therefore be viewed as the finite-horizon standard discounted cost of the augmented MDP with immediate and terminal cost functions given by 
$c_t((s,\psi),a):=0 \mbox{ for all } t=0,1,\ldots, T-1,$ and  $c_T((s,\psi)):= \psi$.
The same idea may be extended to the RS cost functions appearing in the constraints in problems $\CRSMDPFinHorOne$ and $\CRSMDPFinHorTwo$ to obtain a constrained terminal-cost MDP problem with standard discounted cost functions for the objective and constraint functions. Standard LP-based techniques for constrained standard discounted cost MDPs described, for instance in \cite{Derman_Klein}, may then be applied.   For the sake of completeness, we explicitly provide the resulting LP formulation below.

Recall that each of the problems $\CRSMDPFinHorOne$ and $\CRSMDPFinHorTwo$ involve $\hat{M}+\check{M}$ finite-horizon RS constraints and $M+\bar{M}$  finite-horizon standard discounted cost constraints. Of these, $\check{M}$ RS constraints are finite-horizon RS constraints with horizons $\check{T}_{1},\ldots,\check{T}_{\check{M}}$, $\bar{M}$ constraints are finite-horizon standard discounted cost constraints with horizons $\bar{T}_{1},\ldots,\bar{T}_{\bar{M}}$, and the rest arise by truncating infinite-horizon constraints of the RS or standard discounted type. For the sake of simplicity, we consider only the problem $\CRSMDPFinHorOne$, and restrict ourselves to the case where the truncation horizon $T$ is greater than the horizons of all the original finite-horizon constraints, that is, $T>\max_{1\leq l \leq \check{M}}\check{T}_{l}$ and $T>\max_{1\leq k \leq \bar{M}}\bar{T}_{k}$. Define random processes $\{\Psi^{1}_{t}\}_{t=0}^{T},\ldots,\{\Psi^{\hat{M}}\}_{t=0}^{T}$ corresponding to the truncated  RS constraints in a manner analogous to the construction of the process $\{\Psi^{\rm o}_{t}\}_{t=0}^{T}$ described above. Additionally, define random processes $\{\Psi^{\hat{M}+1}_{t}\}_{t=0}^{T},\ldots,$  $\{\Psi^{\hat{M}+\check{M}}\}_{t=0}^{T}$ corresponding to the finite-horizon  RS constraints by  \[\Psi^{ \hat{M}+l}_{t}:=\exp\left[\gamma\sum_{j=0}^{(t\wedge \check{T}_{l})-1}\beta^j \check{C}_{l}(X_j,A_j)\right],\ l\in\{1,\ldots,\check{M}\},\]
where $a\wedge b:=\min\{a,b\}$.
Finally, define the augmented state process $\{Z_{t}\}_{t=0}^{T}$ by letting $Z_{t}:= (X_t, \Psi^{\rm o}_t, \Psi^{1}_t, \ldots, \Psi^{\hat{M}+\check{M}}_t)$ for each $t\in\{0,\ldots,T\}$. For each $t$, denote by  ${\cal S}_t^{\rm A}$ the set of all possible values that the random tuple $Z_{t}$ can take, and note that ${\cal S}_t^{\rm A}$ is finite. Also, given $t\in\{0,\ldots,T\}$, $i\in\{1,\ldots, \hat{M}+\check{M}\}$ and $z\in {\cal S}_t^{\rm A}$, we abuse notation slightly to denote the second and $(i+2)$th components of the tuple $z$ by $\Psi^{\rm o}_t(z)$ and $\Psi^{i}_t(z)$, respectively.

Next, in order to state the LP formulation, we introduce a variable $y(t,z,a)$ for each $t\in\{0,\ldots,T-1\}$, $z\in {\cal S}_t^{\rm A}$ and $a\in \mathcal{A}$, and denote the tuple of all such  variables by ${\bf y}=\{   y (t, z_{t}, a): \  t\le T-1,  \ z_t \in {\cal S}^{\rm A}_t, \  a\in {\cal A}  \}$.  It can be observed that the variable $y(t,z,a)$ is the occupation measure which is interpreted as the probability of being in state $z$ and taking action $a$ at time $t$ (\cite{Derman_Klein,Eitan_CMDP}). The required LP formulation is then given by 
\begin{gather*}
    \min_{{\bf y}} 
%\begin{array}{ccc} {\bf y} = \big \{ y (t, z_{t}, a); \\  t\le T-1,  \\ z_t \in {\cal S}^A_t, \\ a\in {\cal A}  \big \} \end{array}} 
  \sum_{ \begin{array}{ccc} a \in {\cal A},\   z \in {\cal S}^{\rm A}_{T-1} \\  z^{\prime} \in {\cal S}^{\rm A}_{T} \end{array}}  
 \Psi^{\rm o}_T (z^{\prime}) y(T-1,z,a)  p^{\rm A}_{T}(z^{\prime} |z,a) \mbox{ subject to} \\
 \sum_{a\in\mathcal{A}} y(0,z,a)  =  1_{\{z=(x,1,\ldots,1)\in\mathcal{S}_{0}^{\rm A}\}}, \\
% \sum_{a\in\mathcal{A}} y(0,z,a)  &=&  \alpha(z_0) \mbox{ for  all $z \in {\cal S}^{\rm A}_0,$} \mbox{ (recall $z_0 = (x_0, 1, 1, \cdots, 1)$ for some ) } x_0 \in {\cal S},  \\
\sum_{a\in\mathcal{A}} y(t,z^{\prime},a)  =  \sum_{a\in\mathcal{A}} \sum_{z\in  {\cal S}^{\rm A}_{t-1} }   p^{\rm A}_{t} (z^{\prime} |z,a) y(t-1, z,a)  \nonumber  \mbox{ for all }  1\le t  \leq T-1,\  z^{\prime} \in {\cal S}^{\rm A}_t,     \\
%  y(T,z^{\prime})  &=&  \sum_{a\in\mathcal{A}} \sum_{z\in  {\cal S}^{\rm A}_{T-1} }   p^{\rm A} (z^{\prime} |z,a) y(T-1, z,a)  \nonumber  \mbox{ for all }    z^{\prime} \in {\cal S}^{\rm A}_T,     \\
 \sum_{a \in {\cal A}, \  z \in {\cal S}^{\rm A}_{T-1}, \  z^{\prime} \in {\cal S}^{\rm A}_{T} }  
 \Psi^i_T (z^{\prime}) y(T-1,z,a)  p^{\rm A}_{T}(z^{\prime} |z,a)  \leq b^{\rm R}_i  \mbox{ for all } i\in\{1,\ldots,\check{M}+\hat{M}\}, \\
%&& 
%\sum_{t \le T-1, z\in {\cal S}^{\rm A}_t, a\in\mathcal{A}} \beta^t c^i (z, a)   y(t,z,a) 
% \le b^{\rm L}_i  \mbox{ for all } i\in\{1,\ldots,M+\bar{M}\}, \label{lp1}
% \\
\sum_{t \le T-1, z\in {\cal S}^{\rm A}_t, a\in\mathcal{A}} \beta^t c^i (t,z, a)   y(t,z,a) 
 \le b^{\rm L}_i  \mbox{ for all } i\in\{1,\ldots,M+\bar{M}\},
\end{gather*}
  where $x$ is the initial condition of the finite-horizon CRSMDP, 
 \begin{eqnarray}
b^{\rm L}_i =
\begin{cases}
        &b_i-K\beta^T \mbox{ for }i=1,\ldots, M, \\
        &\bar{b}_k \mbox{ for }i=M+k \mbox{ where } k=1,\ldots, \bar{M}, 
\end{cases}
\end{eqnarray}
\begin{eqnarray}
b^{\rm R}_i =
\begin{cases}
        &\frac{\hat{b}_i}{K_T} \mbox{ for }i=1,\ldots, \hat{M}, \\
        &\check{b}_k \mbox{ for }i=\hat{M}+k \mbox{ where } k=1,\ldots, \check{M}, \mbox{ and}
\end{cases}
\end{eqnarray}
\begin{equation}
c^i(t,z, a)  = 
\begin{cases}
        & \hspace{-0.2 cm}C_i(z[1],a)  \mbox{ for } 1\leq i\leq M \mbox{ and for all } 0\leq t\leq T-1,\\
        & \hspace{-0.2 cm} \bar{C}_k(z[1],a) \mbox{ for }i=M+k \mbox{ where } 1\leq k\leq \bar{M},\mbox{ and } t\leq T_k-1,\\
        &  \hspace{-0.2 cm}0 \mbox{ for }i=M+k \mbox{ where } 1\leq k\leq \bar{M},\mbox{ and }  T_k \leq t \leq T-1,
\end{cases} \label{lp2}
\end{equation}
with $z[1]\in\mathcal{S}$ denoting the first component of vector $z$.
%\begin{eqnarray}
%c^i(z, a)  = 
%\begin{cases}
%        & C_i(x,a)  \mbox{ for }i=1,\ldots, M, \\
%        & \bar{C}_k(x,a) \mbox{ for }i=M+k \mbox{ where } k=1,\ldots, \bar{M},
%\end{cases} \label{lp2}
%\end{eqnarray}
%and 
% $b^{\rm L}_i$ and $b^{\rm R}_i$, $i\in\{1,\ldots,M+\bar{M}\},$ are the upper bounds appearing in the standard discounted and RS constraints, respectively, in the problem  $\CRSMDPFinHorOne$ or $\CRSMDPFinHorTwo$. For example, in the  case of $\CRSMDPFinHorOne$, we let 
% The functions $C_i$ and $\bar{C}_k$ appearing in \eqref{lp2} are introduced in the \cref{sec_prb_stmt}, while  $x\in\mathcal{S}$ in \eqref{lp2}  denotes the first component of $z$. 
% \begin{eqnarray}
% c^i(z, a)  = 
% \begin{cases}
%         & C_i(x,a)  \mbox{ for }i=1,\ldots, M \\
%         & \bar{C}_k(x,a) \mbox{ for }i=M+k \mbox{ where } k=1,\ldots, \bar{M},
% \end{cases}
% \end{eqnarray}
% where $x\in\mathcal{S}$ is the first component of $z$ and the functions $C_i$ and $\bar{C}_k$ are introduced in the \cref{sec_prb_stmt}.

The optimal policy $\pi^{*}$ for the CRSMDP may be constructed from the optimal solution $\bf y^{*}$  of the LP problem described above. 
%For this purpose, define $d_{t}:\mathcal{S}_{t}^{\rm A}\times \mathcal{A}\rightarrow [0,1]$ for each $t=0,\ldots,T$ by $ d_t (z, a) = [\sum_{a'\in\mathcal{A}} y^{*}(t, z, a')]^{-1}y^{*}(t, z, a)$, $(z,a)\in \mathcal{S}_{t}^{\rm A}$.
For this purpose, let $\mathcal{S}_{t}^{\rm A}(x)$ denote the set of tuples in $\mathcal{S}_{t}^{\rm A}$ having $x$ as their first component, for each $x\in\mathcal{X}$ and  $t\in\{0,\ldots,T-1\}$. Then the probability of taking action $a\in\mathcal{A}$ in state $s\in\mathcal{S}$ at time $t$ under the time-$t$ decision rule $d_t^*$ of the  optimal policy $\pi^{*}$ is given by \cite{Kavitha_NH_Atul}
\begin{equation}
 d_t^*(a|s) = \frac{y^{*}(t, s, a)}{\sum_{a'\in\mathcal{A}} y^{*}(t, s, a')}, \ \mbox{ where }y^*(t,s,a) = \sum_{z\in \mathcal{S}_{t}^{\rm A}(s)}  y^*(t,z,a). 
 %\mbox{ for } (s,a)\in \mathcal{S}\times \mathcal{A}.
%d^{*}_t (a|x) = \sum_{z\in \mathcal{S}_{t}^{\rm A}(x)}  d_t (z, a),\  x\in\mathcal{X}, a\in\mathcal{A}.
\end{equation}
%\begin{equation}
% d_t^*(a|s) = \frac{y^{*}(t, s, a)}{\sum_{a'\in\mathcal{A}} y^{*}(t, s, a')}, (s,a)\in \mathcal{S}\times \mathcal{A}.
%\end{equation}}
It is possible to  obtain  a LP formulation for a finite-horizon CRSMDP   without first converting it to the standard discounted setting \cite{Kavitha_NH_Atul}. This alternative LP formulation from \cite{Kavitha_NH_Atul} may also be applied to the finite-horizon CRSMDPs $\CRSMDPFinHorOne$ or $\CRSMDPFinHorTwo$. 

\section{Conclusions}\label{sec_conclusions}
We have shown that a CRSMDP involving a fairly general set of constraints in the form of upper bounds on finite- and infinite-horizon RS and standard discounted costs possesses a solution as long as it is feasible. Moreover, near-optimal and near-feasible US policies for the CRSMDP may be found by solving two approximating finite-horizon CRSMDPs obtained by time truncating the  original objective cost and constraint functions and suitably perturbing the bounding constants in the constraints. The approximating finite-horizon CRSMDPs may be solved by using a LP-based approach. 

\appendix

\section{Proofs for Results from Subsection \ref{thm1ssec}}\label{appendixA}
%%%%%%%%%%%%%%%%%%%%%%%%% Proof of Theorem \ref{thm_metric_on_Pi_metrizes_ProdTop} begins   %%%%%%%%%%%%%%

\textit{Proof of Theorem \ref{thm_metric_on_Pi_metrizes_ProdTop}.}
For each $t$, let $\mathcal{Y}_{t}$ denote the set $\mathcal{R}$ equipped with the topology induced by the metric $\mu_{t}$ defined as in  sub-section \ref{sec_Top_Markvpolc}. The topology on $\mathcal{Y}_{t}$ agrees with the subspace topology on $\mathcal{R}\subset \R^{m\times n}$. Since $\mathcal{R}$ is a compact subset of $\R^{m\times n}$, it follows that $\mathcal{Y}_{t}$ is a compact topological space for each $t$. Tychonoff's theorem (refer
\cite[p. 245]{RoyFtz}, \cite[p. 224]{JDug} 
%\cite[Sec. 12.2, p. 245]{RoyFtz}, \cite[Ch XI, Thm I.4, p. 224]{JDug} 
implies that the Cartesian product $\prod_{t=0}^{\infty}\mathcal{Y}_{t}$ is compact under the product topology. A direct application of Theorem 7.2 from \cite[p. 190]{JDug} 
%\cite[Ch. IX, p. 190]{JDug} 
now shows that $\mu$ is a metric which metrizes the product topology on $\prod_{t=0}^{\infty}\mathcal{Y}_{t}$. The last two assertions show that $\left(\prod_{t=0}^{\infty}\mathcal{Y}_{t},\mu\right)$ is a compact metric space. The result now follows by noting that, as a set, the Cartesian product above equals $\mathcal{R}^{\infty}$, which is identifiable with $\Pi_{\rm MR}$. 
\eop

% However, the infinite Cartesian product above equals $\mathcal{R}^{\infty}$ as a set.  

% $\mc{R}^{\infty}$ as a set, and thus 
% , proves $\mu$ is a metric on  (\mathcal{R},\mu_t)(=\Pi_{\rm MR})$ and metrizes the product topology.

% Since, $(\mathcal{R},\mu_t)$ is a compact space for each $t$, by Tychonoff's Theorem, refer ,  its infinite product, $\Pi_{\rm MR}$ is a compact set in the product topology.

%%%%%%%%%%%%%%%%%%%%%%%%% Proof of Theorem \ref{thm_metric_on_Pi_metrizes_ProdTop} ends   %%%%%%%%%%%%%%
%%%%%%%%%%%%%%%%%%%%%%%
%\hspace{-0.2 cm}\line(1,0){400}

The proofs of lemmas \ref{lemma_norm_rdminrf_PdminPf} and  \ref{lemma_zeta_and_its_differ_bounds}, and Theorem  \ref{thm_unif_cont_of_USTOCost_prodtop} below make use of the following inequalities for $A,B\in\R^{m\times n}$ and $y,z\in\R^{n}$. 
\begin{eqnarray}
\|A\odot B\|_\infty &\leq & \|A\|_{\rm max} \hspace{0.1 cm}\|B\|_\infty \leq \|A\|_{\infty} \hspace{0.1 cm}\|B\|_\infty, \label{eq_normAschurprodB_leq_normmaxA_times_norminfB}\\
%\|x^\prime A\|_\infty &\leq & \|x\|_{1} \hspace{0.1 cm}\|A\|_{\rm max} \leq \|x\|_{1} \hspace{0.1 cm}\|A\|_\infty,\\
\|AB\|_{\infty}&\leq &\|A\|_{\infty}\|B\|_{\infty},\ \|Ay\|_\infty \leq \|A\|_\infty \hspace{0.05 cm}\|y\|_\infty, \label{eq_norminfAx_leq_norminfA_mult_norminfx} \\
|z^\prime y| &\leq & \|z\|_1 \|y\|_\infty. \label{eq_Abs_xtrnspsy_leq_norm1x_times_norminfy}
\end{eqnarray}
We also note that, if $A,B\in \R^{m\times m}$ are nonnegative and row-stochastic, then $\|A\|_{\rm max}\leq 1=\|A\|_{\infty}$, and $AB$ is nonnegative and row-stochastic. 

%%%%%%%%%%%%%%%%%%%%% Proof of Lemma  \ref{lemma_norm_rdminrf_PdminPf} begins %%%%%%%%%%%%%%%%%%%%
\textit{Proof of Lemma \ref{lemma_norm_rdminrf_PdminPf}.}
First, we claim that the inequalities
% \begin{eqnarray}
% \|R_d-R_f\|_\infty &\leq& C\|d-f\|_\infty ,\label{eq_norm_rd_min_rf_leq_norm_dminf}\\
% \|P_d-P_f\|_\infty &\leq &  \|d-f\|_\infty, \label{eq_norm_Pd_min_Pf_leq_norm_dminf}
% \end{eqnarray} 
\begin{eqnarray}
\|R_d-R_f\|_\infty \leq C\|d-f\|_\infty \ \mbox{ and } \|P_d-P_f\|_\infty \leq  \|d-f\|_\infty \label{eq_normdiff_rAndP}
\end{eqnarray} 
hold for every  $d,f\in\mathcal{R}$.   Using $R_{d}=(R\odot d)\one_{n\times 1}$, \eqref{eq_norminfAx_leq_norminfA_mult_norminfx}, \eqref{eq_normAschurprodB_leq_normmaxA_times_norminfB}  and $\|\one_{n\times 1}\|_{\infty}=1$ in that order gives $\|R_d-R_f\|_\infty = \|(R\odot d - R\odot f)\one_{n\times 1}\|_\infty \leq \|R\odot(d-f)\|_\infty \|\one_{n\times 1}\|_{\infty}\leq \|R\|_{\rm max}\|d-f\|_\infty$. 
%$\|R_d-R_f\|_\infty \leq \|R\|_{\rm max}\|d-f\|_\infty$.} 
Next, using $\|R\|_{\rm max}\leq C$ yields the first equation of \eqref{eq_normdiff_rAndP}.

The sum of the absolute values of the entries of the $i$th row of the matrix $P_{d}-P_{f}$  
%\sout{ in \eqref{eq_normdiff_rAndP}} 
satisfies
\begin{eqnarray}
 \sum_{j=1}^{m}|(P_d)_{i,j} -(P_f)_{i,j}| &=& \sum_{j=1}^{m}\left|\sum_{k=1}^n p(s_j|s_i,a_k)\{d(a_{k}|s_{i})-f(a_{k}|s_{i})\}\right|\nonumber\\
 &\leq & \sum_{j=1}^{m}\sum_{k=1}^n |p(s_j|s_i,a_k)|\cdot|d(a_{k}|s_{i})-f(a_{k}|s_{i})|\nonumber\\
  &=& \sum_{k=1}^{n}\left(\sum_{j=1}^m p(s_j|s_i,a_k)\right)|d(a_{k}|s_{i})-f(a_{k}|s_{i})|  \nonumber\\
  &=&\sum_{k=1}^n |d(a_{k}|s_{i})-f(a_{k}|s_{i})| \leq  \|d-f\|_\infty. 
\end{eqnarray}
% \begin{eqnarray}
% \sum_{j=1}^{m}|(P_d)_{i,j} -(P_f)_{i,j}| &=& \sum_{j=1}^{m}\left|\sum\limits_{k=1}^n p(s_j|s_i,a_k)\{d(a_{k}|s_{i})-f(a_{k}|s_{i})\}\right|\nonumber \\ &\leq&  \sum_{j=1}^{m}\sum\limits_{k=1}^n |p(s_j|s_i,a_k)|\cdot|d(a_{k}|s_{i})-f(a_{k}|s_{i})|  \nonumber \\ 
% &=& {\color{red}\sum\limits_{k=1}^n |d(a_{k}|s_{i})-f(a_{k}|s_{i})| \leq  \|d-f\|_\infty. \label{eq_Pd_minus_Pf}}
% \end{eqnarray}
Thus every absolute row sum of the $m\times m$ matrix $P_d-P_f$ is bounded above by $\|d-f\|_\infty$. It now follows that the maximum absolute row sum norm $\|P_d-P_f\|_\infty$ satisfies 
%\sout{\eqref{eq_norm_Pd_min_Pf_leq_norm_dminf}} 
 the second inequality in \eqref{eq_normdiff_rAndP}. 

Let $\pi_1,\pi_2\in\Pi_{\rm MR}$ be as in the statement of the lemma. We claim that
\begin{eqnarray}
\|P^{\pi_1}_t-P^{\pi_2}_t\|_\infty &\leq& \sum\limits_{i=0}^{t-1}\|d_i-f_i\|_\infty, \label{eq_Ppi1t_min_Ppi2t_normbd}
\end{eqnarray}
for every $t>0$. To prove our claim, we first note that, for every $t>0$, the time-$t$ transition probability matrix  $P_{d_{t}}$ and the $t$-step transition probability matrix $P^{\pi_1}_{t-1}$ are both  nonnegative and row-stochastic. Next, for $t>0$, we have 
% $
% \|P^{\pi_1}_t-P^{\pi_2}_t\|_\infty = \|P^{\pi_1}_{t-1}P_{d_{t-1}}-P^{\pi_2}_{t-1}P_{f_{t-1}}\|_\infty 
% = {\color{cyan}\|P^{\pi_1}_{t-1}(P_{d_{t-1}}-P_{f_{t-1}})+(P^{\pi_1}_{t-1}-P^{\pi_2}_{t-1})P_{f_{t-1}}\|_\infty }
% \leq \|P^{\pi_1}_{t-1}(P_{d_{t-1}}-P_{f_{t-1}})\|_\infty + \|(P^{\pi_1}_{t-1}-P^{\pi_2}_{t-1})P_{f_{t-1}}\|_\infty 
% \leq  \|P^{\pi_1}_{t-1}\|_\infty\|P_{d_{t-1}}-P_{f_{t-1}}\|_\infty + \|P^{\pi_1}_{t-1}-P^{\pi_2}_{t-1}\|_\infty \|P_{f_{t-1}}\|_\infty
% \leq \|d_{t-1}-f_{t-1}\|_\infty +\|P^{\pi_1}_{t-1}-P^{\pi_2}_{t-1}\|_\infty, 
% $
 \begin{eqnarray}
\|P^{\pi_1}_t-P^{\pi_2}_t\|_\infty &=& \|P^{\pi_1}_{t-1}P_{d_{t-1}}-P^{\pi_2}_{t-1}P_{f_{t-1}}\|_\infty \nonumber \\ 
%  &=& {\color{cyan}\|P^{\pi_1}_{t-1}(P_{d_{t-1}}-P_{f_{t-1}})+(P^{\pi_1}_{t-1}-P^{\pi_2}_{t-1})P_{f_{t-1}}\|_\infty }\nonumber\\
 &\leq& \|P^{\pi_1}_{t-1}(P_{d_{t-1}}-P_{f_{t-1}})\|_\infty + \|(P^{\pi_1}_{t-1}-P^{\pi_2}_{t-1})P_{f_{t-1}}\|_\infty \nonumber\\
 &\leq & \|P^{\pi_1}_{t-1}\|_\infty\|P_{d_{t-1}}-P_{f_{t-1}}\|_\infty + \|P^{\pi_1}_{t-1}-P^{\pi_2}_{t-1}\|_\infty \|P_{f_{t-1}}\|_\infty \nonumber\\
&\leq& \|d_{t-1}-f_{t-1}\|_\infty +\|P^{\pi_1}_{t-1}-P^{\pi_2}_{t-1}\|_\infty, \label{eq_Ppi1t_min_Ppi2t_normbd_recsv}
\end{eqnarray}
where we have used the first inequality in \eqref{eq_norminfAx_leq_norminfA_mult_norminfx},  
%the inequality \sout{ \eqref{eq_norm_Pd_min_Pf_leq_norm_dminf}} 
 the second inequality in \eqref{eq_normdiff_rAndP}, and the fact that $\|P^{\pi_1}_{t-1}\|_\infty=\|P_{f_{t-1}}\|_\infty=1$. Solving the recursion \eqref{eq_Ppi1t_min_Ppi2t_normbd_recsv} and noting that $\|P^{\pi_1}_{0}-P^{\pi_2}_{0}\|_\infty=0$ shows that  \eqref{eq_Ppi1t_min_Ppi2t_normbd} holds.

Next, for $t\geq 0$, we have 
%$\|P_t^{\pi_1}R_{d_t} -P_t^{\pi_2}R_{f_t}\|_\infty \leq \|P_t^{\pi_1}(R_{d_t}-R_{f_t})\|_\infty  + \|(P_t^{\pi_1} -P_t^{\pi_2})R_{f_t}\|_\infty \leq \|P_t^{\pi_1}\|_\infty  \|R_{d_t}-R_{f_t}\|_\infty + \|P_t^{\pi_1} -P_t^{\pi_2}\|_\infty\|R_{f_t}\|_\infty \leq C \sum_{i=0}^{t}\|d_i-f_i\|_\infty$,
\begin{eqnarray}
 \|P_t^{\pi_1}R_{d_t} -P_t^{\pi_2}R_{f_t}\|_\infty &\leq& \|P_t^{\pi_1}(R_{d_t}-R_{f_t})\|_\infty  + \|(P_t^{\pi_1} -P_t^{\pi_2})R_{f_t}\|_\infty \nonumber\\
 &\leq&\|P_t^{\pi_1}\|_\infty  \|R_{d_t}-R_{f_t}\|_\infty + \|P_t^{\pi_1} -P_t^{\pi_2}\|_\infty\|R_{f_t}\|_\infty \nonumber\\
 &\leq&C \sum_{i=0}^{t}\|d_i-f_i\|_\infty,
\end{eqnarray}
%\begin{eqnarray}
% \|P_t^{\pi_1}R_{d_t} -P_t^{\pi_2}R_{f_t}\|_\infty &=& \|P_t^{\pi_1}(R_{d_t}-R_{f_t})+ (P_t^{\pi_1} -P_t^{\pi_2})R_{f_t}\|_\infty \nonumber\\
% &\leq&{\color{red} \|P_t^{\pi_1}(R_{d_t}-R_{f_t})\|_\infty  + \|(P_t^{\pi_1} -P_t^{\pi_2})R_{f_t}\|_\infty} \nonumber\\
% &\leq& \|P_t^{\pi_1}\|_\infty  \|R_{d_t}-R_{f_t}\|_\infty + \|P_t^{\pi_1} -P_t^{\pi_2}\|_\infty\|R_{f_t}\|_\infty {\color{red}\leq C \sum\limits_{i=0}^{t}\|d_i-f_i\|_\infty},\nonumber
% &\leq& C\|d_t-f_t\|_\infty + C \sum\limits_{i=0}^{t-1}\|d_i-f_i\|_\infty = C \sum\limits_{i=0}^{t}\|d_i-f_i\|_\infty,
% \end{eqnarray}
where we  have used the second inequality in \eqref{eq_norminfAx_leq_norminfA_mult_norminfx} along with   
%\sout{\eqref{eq_norm_rd_min_rf_leq_norm_dminf}} 
 the first inequality in \eqref{eq_normdiff_rAndP},  $\|P^{\pi_1}_{t}\|_\infty=1$ and $\|R_{f_t}\|_\infty\leq C$. This proves the lemma.
\eop
%%%%%%%%%%%%%%%%%%%%% Proof of Lemma  \ref{lemma_norm_rdminrf_PdminPf} ends %%%%%%%%%%%%%%%%%%%%

%\hspace{-0.4 cm}\line(1,0){420}

%%%%%%%%%%%%%%%%%%%%%%%%% Proof of Lemma \ref{lem_RSCost_finhor_cgsto_infinhor} begins %%%%%%%%%%%%%%%%%%%
\textit{Proof of Lemma \ref{lem_linCost_finhor_cgsto_infinhor}.}
Let $\pi\in\Pi_{\rm MR}$, $T>0$ and $x\in\mathcal{S}$ be arbitrary. We have  $ |\LDCost{\pi}{R} - \LDCostT{\pi}{R}{T}| = \left|\E^{\pi}_{x}\left[\sum_{t=T}^{\infty}\beta^t R(X_{t},A_{t})\right]\right|\leq \E^{\pi}_{x}\left[\sum_{t=T}^{\infty}\beta^t |R(X_t,A_t)|\right]
\leq \beta^T K. $
% {\color{red}\begin{eqnarray}
% |\LDCost{\pi}{R} - \LDCostT{\pi}{R}{T}| = \left|\E^{\pi}_{x}\left[\sum_{t=T}^{\infty}\beta^t R(X_{t},A_{t})\right]\right|\leq \E^{\pi}_{x}\left[\sum_{t=T}^{\infty}\beta^t |R(X_t,A_t)|\right]
% \leq \beta^T K. $}
%&=& \left|\E^{\pi}_{x}\left[\sum\limits_{t=T}^{\infty}\beta^t R(X_{t},A_{t})\right]\right|\nonumber\\
%&\leq& %\E^{\pi}_{x}\left[\left|\sum\limits_{t=T}^{\infty}\beta^t R(X_t,A_t)\right|\right]\leq
% \E^{\pi}_{x}\left[\sum\limits_{t=T}^{\infty}\beta^t |R(X_t,A_t)|\right]
% \leq \beta^T K.\label{eq_expec_inf_sum}
% \end{eqnarray}}
% The expectation is taken inside the infinite sum in \eqref{eq_expec_inf_sum} because the random variable $\sum\limits_{t=T}^{\infty}\beta^t R_t$ converges and using the \cite[Theorem 9.2]{JacodProtter}.
Since the bound in the previous inequality is independent of the policy $\pi$, the result follows by letting  $T\to\infty$. 
\eop

%%%%%%%%%%%%%%%%%%%%%%%%% Proof of Lemma \ref{lem_RSCost_finhor_cgsto_infinhor} ends %%%%%%%%%%%%%%%%%%%

%%%%%%%%%%%%%%%%%%%%%%55 Proof of Theorem \ref{thm_contuty_lin_disc_cost_prod_toplogy} begins%%%%%%%%%%%%
\textit{Proof of Theorem \ref{thm_contuty_lin_disc_cost_prod_toplogy}. }%For a given $\epsilon>0$, choose $\delta^\prime$ such that $0<\delta^\prime < \frac{\epsilon}{C\alpha S}$. Then, from the  inequality \eqref{eq_LDexpec_norm_diff_leq_dist_pi1pi2},  for $\mu(\pi_1,\pi_2) <\delta^\prime$ implies  $\|v^{\pi_1} - v^{\pi_2}\|_\infty < \epsilon$, proving the the map $\pi\to v^\pi$ is uniformly continuous under the product topology induced by $\mu$.
Let $\pi_{1}=\{d_{t}\}_{t=0}^{\infty}$ and  $\pi_{2}=\{f_{t}\}_{t=0}^{\infty}$ be two policies in $\Pi_{\rm MR}$. 
Using the vector expression in \eqref{eq_expect_in_matrxiandvector_form} for the finite-horizon standard discounted cost and  applying Lemma \ref{lemma_norm_rdminrf_PdminPf} gives 
%$
%\|\LDCostTVec{\pi_1}{R}{T} - \LDCostTVec{\pi_2}{R}{T}\|_\infty 
% \left\|\sum\limits_{t=0}^{T-1}\beta^t P^{\pi_1}_{t}R_{_{d_t}}-\sum\limits_{t=0}^{T-1}\beta^t P^{\pi_2}_{t}R_{_{f_t}}\right\|_\infty 
%=\left\|\sum\limits_{t=0}^{T-1}\beta^t (P^{\pi_1}_{t}R_{_{d_t}}-P^{\pi_2}_{t}R_{_{f_t}})\right\|_\infty
%\leq \sum\limits_{t=0}^{T-1}\beta^t \|P^{\pi_1}_{t}R_{_{d_t}} - P^{\pi_2}_{t}R_{_{f_t}}\|_\infty \leq \sum\limits_{t=0}^{T-1}\beta^t C \sum\limits_{i=0}^{t}\|d_i-f_i\|_\infty 
%=C\sum_{i=0}^{T-1}\sum_{t=i}^{T-1}\beta^t\|d_i-f_i\|_\infty =C\sum_{i=0}^{T-1}\frac{\beta^i(1-\beta^{T-i})}{1-\beta}\|d_i-f_i\|_\infty
%\leq K\sum_{i=0}^{T-1}\beta^i\|d_i-f_i\|_\infty =K\sum\limits_{i=0}^{T-1}\left(\frac{\beta}{\delta}\right)^i\delta^i\|d_i-f_i\|_\infty 
%\leq K\sum_{i=0}^{T-1}\left(\frac{\beta}{\delta}\right)^i \mu(\pi_1,\pi_2)
%= \frac{K(\delta^T-\beta^T)}{\delta^{T-1}(\delta-\beta)}\mu(\pi_1,\pi_2).$
 \begin{eqnarray}
 \|\LDCostTVec{\pi_1}{R}{T} - \LDCostTVec{\pi_2}{R}{T}\|_\infty &=& \left\|\sum\limits_{t=0}^{T-1}\beta^t P^{\pi_1}_{t}R_{_{d_t}}-\sum\limits_{t=0}^{T-1}\beta^t P^{\pi_2}_{t}R_{_{f_t}}\right\|_\infty\nonumber\\
 &=&\left\|\sum\limits_{t=0}^{T-1}\beta^t (P^{\pi_1}_{t}R_{_{d_t}}-P^{\pi_2}_{t}R_{_{f_t}})\right\|_\infty \nonumber\\
 &\leq& \sum\limits_{t=0}^{T-1}\beta^t \|P^{\pi_1}_{t}R_{_{d_t}} - P^{\pi_2}_{t}R_{_{f_t}}\|_\infty \leq \sum\limits_{t=0}^{T-1}\beta^t C \sum\limits_{i=0}^{t}\|d_i-f_i\|_\infty \nonumber\\
 &=&C\sum\limits_{i=0}^{T-1}\sum\limits_{t=i}^{T-1}\beta^t\|d_i-f_i\|_\infty =C\sum\limits_{i=0}^{T-1}\frac{\beta^i(1-\beta^{T-i})}{1-\beta}\|d_i-f_i\|_\infty\nonumber\\
 &\leq&K\sum\limits_{i=0}^{T-1}\beta^i\|d_i-f_i\|_\infty =K\sum\limits_{i=0}^{T-1}\left(\frac{\beta}{\delta}\right)^i\delta^i\|d_i-f_i\|_\infty \nonumber\\
 &\leq&K\sum\limits_{i=0}^{T-1}\left(\frac{\beta}{\delta}\right)^i \mu(\pi_1,\pi_2)
 = \frac{K(\delta^T-\beta^T)}{\delta^{T-1}(\delta-\beta)}\mu(\pi_1,\pi_2).
 \end{eqnarray}
Thus, the vector valued map $\pi\to \LDCostTVec{\pi}{R}{T}$ on $\Pi_{\rm MR}$ is Lipschitz continuous with Lipschitz constant $\frac{K(\delta^T-\beta^T)}{\delta^{T-1}(\delta-\beta)}.$

Next, choose $\epsilon>0$. We know from Lemma \ref{lem_linCost_finhor_cgsto_infinhor} that there exists $T^*>0$ such that $|\LDCostT{\pi}{R}{T} -\LDCost{\pi}{R}|< \epsilon/2$ for all $T\geq T^*$, $\pi\in\Pi_{\rm MR}$ and $x\in\mathcal{S}$. 
Choose $T\geq T^*$. We then have
 %$\|\LDCostVec{\pi_1}{R} - \LDCostVec{\pi_2}{R}\|_\infty \leq \|\LDCostVec{\pi_1}{R}- \LDCostTVec{\pi_1}{R}{T}\|_\infty +\|\LDCostTVec{\pi_1}{R}{T}-\LDCostTVec{\pi_2}{R}{T}\|_\infty + \|\LDCostTVec{\pi_2}{R}{T}- \LDCostVec{\pi_2}{R}\|_\infty
%\leq \frac{\epsilon}{2} + \frac{K(\delta^T-\beta^T)}{\delta^{T-1}(\delta-\beta)}\mu(\pi_1,\pi_2) + \frac{\epsilon}{2}\leq \epsilon+ \frac{K\delta}{(\delta-\beta)}\mu(\pi_1,\pi_2). $
 \begin{eqnarray}
 \|\LDCostVec{\pi_1}{R} - \LDCostVec{\pi_2}{R}\|_\infty &\leq& \|\LDCostVec{\pi_1}{R}- \LDCostTVec{\pi_1}{R}{T}\|_\infty +\|\LDCostTVec{\pi_1}{R}{T}-\LDCostTVec{\pi_2}{R}{T}\|_\infty + \|\LDCostTVec{\pi_2}{R}{T}- \LDCostVec{\pi_2}{R}\|_\infty\nonumber\\
 &\leq& \frac{\epsilon}{2} + \frac{K(\delta^T-\beta^T)}{\delta^{T-1}(\delta-\beta)}\mu(\pi_1,\pi_2) + \frac{\epsilon}{2}\leq \epsilon+ \frac{K\delta}{(\delta-\beta)}\mu(\pi_1,\pi_2).\nonumber
 %&=&\epsilon + \frac{K(\delta^T-\beta^T)}{\delta^{T-1}(\delta-\beta)}\mu(\pi_1,\pi_2) \leq \epsilon+ \frac{K\delta}{(\delta-\beta)}\mu(\pi_1,\pi_2).
\end{eqnarray}
Letting $\epsilon\to 0$ now shows that the map $\pi\to\LDCostVec{\pi}{R}$ is Lipschitz continuous with Lipschitz constant $\frac{K\delta}{(\delta-\beta)}$. \eop 

\textit{Proof of Proposition \ref{thm_matrix_formulaton_FinHor_RSCost}:} 
Note that \eqref{eq_zeta0_eql_ones_matrix} directly follows by letting $t=T-1$ in \eqref{eq_qvaluedef}. Next, choose $t\in\{1,\ldots,T-1\}$ and $(s,a)\in\mathcal{S}\times \mathcal{A}$. Applying \eqref{eq_qvaluedef}, we have
\begin{eqnarray}
  \Qfunc{t-1}{T}{\pi}(s,a)&=&\E^\pi\left[e^{\gamma\sum\limits_{\tau=t-1}^{T-1}\beta^\tau R(X_\tau,A_\tau)}\bigg|X_{t-1}=s,A_{t-1}=a\right] \nonumber \\
 % &=&e^{\gamma\beta^{t-1}R(s,a)}\E^\pi\left[e^{\gamma\sum\limits_{\tau=t}^{T-1}\beta^\tau R(X_\tau,A_\tau)}\bigg|X_{t-1}=s,A_{t-1}=a\right]\nonumber \\
 &=&e^{\gamma\beta^{t-1}R(s,a)}\E^\pi\left[\E^\pi\left.\left[e^{\gamma\sum\limits_{\tau=t}^{T-1}\beta^\tau R(X_\tau,A_\tau)}\bigg|X_t,A_t\right]\right|X_{t-1}=s,A_{t-1}=a\right] \nonumber \\
  &=&e^{\gamma\beta^{t-1}R(s,a)}\E^{\pi}[\Qfunc{t}{T}{\pi}(X_{t},A_{t})|X_{t-1}=s,A_{t-1}=a] \nonumber \\
 &=& e^{\gamma\beta^{t-1}R(s,a)}\sum_{i=1}^{m}\sum_{j=1}^{n}p(s_{i}|s,a)d_{t}(a_{j}|s_{i})\Qfunc{t}{T}{\pi}(s_{i},a_{j}).
 % \nonumber 
% \\ &=&e^{\gamma\beta^{t-1}R(s,a)}p(\cdot|s,a)^\prime (d_{t}\odot \Qfunc{t}{T}{\pi})\one_{n\times 1}. 
 \end{eqnarray}
The double summation above equals $p(\cdot|s,a)^\prime (d_{t}\odot \Qfunc{t}{T}{\pi})\one_{n\times 1}$, and \eqref{zetat_interms_zetatmin1} follows.  

To prove \eqref{eq_JpiT_In_Matrix_Form_eql_expgammar_d0_zetaTmin1}, choose $i\in\{1,\ldots,m\}$, and note that 
\begin{eqnarray*}
 \RSCostTVec{\pi}{\gamma}{R}{T}(s_i) 
 % CC&=& \E^\pi\left[e^{\gamma\sum\limits_{\tau=0}^{T-1}\beta^\tau R(X_\tau,A_\tau)}\bigg|X_{0}=s_i\right]\nonumber\\
  &=& \E^\pi\left[\E^\pi\left.\left[e^{\gamma\sum\limits_{\tau=0}^{T-1}\beta^\tau R(X_\tau,A_\tau)}\bigg|X_0,A_0\right]\right|X_{0}=s_i\right] \\
 % % CC {\color{red} = \E^\pi\left[\E^\pi\left[e^{\gamma\sum\limits_{\tau=0}^{T-1}\beta^\tau R(X_\tau,A_\tau)}\bigg|X_0=s_i,A_0\right]\right]} 
  %\nonumber\\
 &=&\E^\pi\left[\Qfunc{0}{T}{\pi}(X_{0},A_{0})\bigg|X_{0}=s_i\right] 
 % CCC= \E^\pi\left[\Qfunc{0}{T}{\pi}(s_i,A_{0})\right] 
 = \sum\limits_{j=1}^{n} d_{0}(a_j|s_i)\Qfunc{0}{T}{\pi}(s_{i},a_{j}). %CClabel{eq_ith_elmt_RSCOst_interms_zeta}
 \end{eqnarray*}
The summation above  
%in \eqref{eq_ith_elmt_RSCOst_interms_zeta} 
is exactly the $i$-th element of the vector $( d_0\odot \Qfunc{0}{T}{\pi})\one_{n\times 1}$. This completes the proof.
\eop

%\line(1,0){360}

%%%%%%%%%%%%%%%%%%%%% Proof of Lemma  \ref{lemma_zeta_and_its_differ_bounds} begins %%%%%%%%%%%%%%%%%%%%

\textit{Proof of Lemma \ref{lemma_zeta_and_its_differ_bounds}.}
To begin, note that \eqref{eq_maxnorm_zetaT_bound} follows directly from \eqref{eq_zeta0_eql_ones_matrix} in the case $t=T-1$. Next, 
 choose 
$(s,a)\in\mathcal{S}\times \mathcal{A}$ and $1\leq t\leq T-1$. 
% For notational convenience, temporarily denote $R(s,a)$ by $R$, and recall that $R\leq C$. 
From \eqref{zetat_interms_zetatmin1}, we have 
\begin{eqnarray}
|\Qfunc{t-1}{T}{\pi_1}(s,a)| &=&e^{\gamma \beta^{t-1}R(s,a)} |p(\cdot|s,a)^\prime(  d_{t}\odot \Qfunc{t}{T}{\pi_1})\one_{n\times 1}|\nonumber \\
&\leq & e^{|\gamma|\beta^{t-1}C}\|p(\cdot|s,a)\|_1\cdot \|(d_{t}\odot \Qfunc{t-1}{T}{\pi_1})\one_{n\times 1}\|_\infty \label{eq_splitmod_norm1_times_norminf} \\
&\leq &e^{|\gamma|\beta^{t-1}C}\|d_{t}\odot \Qfunc{t}{T}{\pi_1} \|_\infty \cdot \|\one_{n\times 1}\|_\infty \label{eq_splitnorminf}\\
&\leq & e^{|\gamma|\beta^{t-1}C}\|d_{t}\|_\infty\cdot \|\Qfunc{t}{T}{\pi_1}\|_{\rm max}
~~ \leq  e^{|\gamma|\beta^{t-1}C}\|\Qfunc{t}{T}{\pi_1}\|_{\rm max}. \label{eq_zetatofsa_bdd}
\end{eqnarray}
Inequalities \eqref{eq_splitmod_norm1_times_norminf}, \eqref{eq_splitnorminf} and \eqref{eq_zetatofsa_bdd} follow by using  \eqref{eq_Abs_xtrnspsy_leq_norm1x_times_norminfy}, \eqref{eq_norminfAx_leq_norminfA_mult_norminfx} and \eqref{eq_normAschurprodB_leq_normmaxA_times_norminfB}, respectively,  along with $|R(s,a)|\leq C$ and  $\|p(\cdot|s,a)\|_{1}=\|\one_{n\times 1}\|_{\infty}=\|d_{t}\|_{\infty}=1$. 
Since the bound in \eqref{eq_zetatofsa_bdd} is independent of $s$ and $a$, it follows that   $\|\Qfunc{t-1}{T}{\pi_1}\|_{\rm max} \leq e^{|\gamma|\beta^{t-1}C}\|\Qfunc{t}{T}{\pi_1}\|_{\rm max}$ for all  $1\leq t\leq T-1$. Solving this recursion and noting from \eqref{eq_zeta0_eql_ones_matrix} that $\|\Qfunc{T-1}{T}{\pi_{1}}\|_{\rm max}\leq e^{|\gamma|\beta^{T-1}C}$ gives $\|\Qfunc{t}{T}{\pi_{1}}\|_{\rm max}\leq \exp\left(|\gamma| C\beta^{t}\sum_{i=0}^{T-t-1}\beta^{i}\right)=\exp[|\gamma|K(\beta^{t}-\beta^{T})]$ for every $t$ satisfying $0\leq t < T-1$. This  proves   \eqref{eq_maxnorm_zetaT_bound}.

To prove \eqref{eq_maxnorm_zetaTOfpi1_min_norm_zetaTOfpi2_bound}, we first derive an upper bound on the term $\|d_{t}\odot \Qfunc{t}{T}{\pi_1} -  f_{t}\odot \Qfunc{t}{T}{\pi_2}\|_\infty $ for a given $t$ satisfying $0\leq t\leq T-1$. We have 
\begin{eqnarray}
&&\|d_{t}\odot \Qfunc{t}{T}{\pi_1} -  f_{t}\odot \Qfunc{t}{T}{\pi_2}\|_\infty 
% &=& {\color{cyan}\|d_{t}\odot(\Qfunc{t}{T}{\pi_1}-\Qfunc{t}{T}{\pi_2})+ \Qfunc{t}{T}{\pi_2}\odot( d_{t}-  f_{t})\|_\infty }\nonumber\\
\leq  \|d_{t}\odot(\Qfunc{t}{T}{\pi_1}-\Qfunc{t}{T}{\pi_2})\|_\infty + \| \Qfunc{t}{T}{\pi_2}\odot( d_{t}-  f_{t})\|_\infty \nonumber\\
&\leq & \|d_{t}\|_\infty\cdot\|\Qfunc{t}{T}{\pi_1}-\Qfunc{t}{T}{\pi_2}\|_{\rm max} + \| \Qfunc{t}{T}{\pi_2}\|_{\rm max}\cdot\| d_{t}-  f_{t}\|_\infty \label{eq_split_norminfSprod_norminftimesnormmax1}\\
&\leq & \|\Qfunc{t}{T}{\pi_1}-\Qfunc{t}{T}{\pi_2}\|_{\rm max} + \| \Qfunc{t}{T}{\pi_2}\|_{\rm max}\cdot\| d_{t}-  f_{t}\|_\infty, \label{eq_norm_dTmintzetatmin1pi1_minus_fTmintzetatmin1pi2_bdd}
\end{eqnarray}
where the inequality \eqref{eq_split_norminfSprod_norminftimesnormmax1} follows from  \eqref{eq_normAschurprodB_leq_normmaxA_times_norminfB}.

 Next, let $(s,a)\in\mathcal{S}\times \mathcal{A}$ and $1\leq t\leq T-1$. 
%As before, temporarily denote $R(s,a)$ by $R$. 
From \eqref{zetat_interms_zetatmin1}, we  have 
\begin{eqnarray}
&&|\Qfunc{t-1}{T}{\pi_1}(s,a)-\Qfunc{t-1}{T}{\pi_2}(s,a)|\nonumber\\
&=&|e^{\gamma\beta^{t-1}R(s,a)}p(.|s,a)^\prime( d_{t}\odot \Qfunc{t}{T}{\pi_1} - f_{t}\odot \Qfunc{t}{T}{\pi_2})\one_{n\times 1}| \nonumber\\
&\leq&  e^{|\gamma|\beta^{t-1}C}\|p(.|s,a)\|_1 \cdot\|(d_{t}\odot \Qfunc{t}{T}{\pi_1} -  f_{t}\odot \Qfunc{t}{T}{\pi_2})\one_{n\times 1}\|_\infty\label{eq_splittingnorms1}\\
&\leq & e^{|\gamma|\beta^{t-1} C}\|d_{t}\odot \Qfunc{t}{T}{\pi_1} -  f_{t}\odot \Qfunc{t}{T}{\pi_2}\|_\infty\cdot \|\one_{n\times 1}\|_\infty \label{eq_splitting_Sprod_norminf} \\
&= & e^{|\gamma|\beta^{t-1} C} \|d_{t}\odot \Qfunc{t}{T}{\pi_1} -  f_{t}\odot \Qfunc{t}{T}{\pi_2}\|_\infty \label{eq_splitting_Sprod_norminf2}\\
&\leq & e^{|\gamma|\beta^{t-1} C} (\|\Qfunc{t}{T}{\pi_1}-\Qfunc{t}{T}{\pi_2}\|_{\rm max} + \| \Qfunc{t}{T}{\pi_2}\|_{\rm max}\cdot\| d_{t}-  f_{t}\|_\infty). \label{eq_zetapi1_minus_zetap12_sath_bound}
\end{eqnarray}
The inequalities \eqref{eq_splittingnorms1}, \eqref{eq_splitting_Sprod_norminf} and \eqref{eq_splitting_Sprod_norminf2}  follow by  applying \eqref{eq_Abs_xtrnspsy_leq_norm1x_times_norminfy}, \eqref{eq_norminfAx_leq_norminfA_mult_norminfx} and \eqref{eq_normAschurprodB_leq_normmaxA_times_norminfB}, respectively, along with $|R(s,a)|\leq C$ and $\|p(\cdot|s,a)\|_{1}=\|\one_{n\times 1}\|_{\infty}=1$. The last inequality \eqref{eq_zetapi1_minus_zetap12_sath_bound} follows from \eqref{eq_norm_dTmintzetatmin1pi1_minus_fTmintzetatmin1pi2_bdd}. 

Note that the  bound in \eqref{eq_zetapi1_minus_zetap12_sath_bound} is independent of $s$ and $a$. Also, we had chosen $1\leq t\leq T-1$ arbitrarily. Therefore, it follows that 
\begin{equation}
\|\Qfunc{t-1}{T}{\pi_1}-\Qfunc{t-1}{T}{\pi_2}\|_{\rm max} \leq e^{|\gamma|\beta^{t-1} C} \left(\|\Qfunc{t}{T}{\pi_1}-\Qfunc{t}{T}{\pi_2}\|_{\rm max} + \| \Qfunc{t}{T}{\pi_2}\|_{\rm max}\cdot \| d_{t}-  f_{t}\|_\infty\right), \label{eq_norm_zetapi1_minus_zetap12_bound}
\end{equation}
for all $1\leq t\leq T-1$. 

On recalling 
from \eqref{eq_zeta0_eql_ones_matrix} 
that $\Qfunc{T-1}{T}{\pi_{1}}=\Qfunc{T-1}{T}{\pi_{2}}$, we see that \eqref{eq_maxnorm_zetaTOfpi1_min_norm_zetaTOfpi2_bound} holds for $t=T-1$. To complete the proof by induction, suppose \eqref{eq_maxnorm_zetaTOfpi1_min_norm_zetaTOfpi2_bound} holds for $t=k$, where $1\leq k\leq T-1$. Letting $t=k$ in \eqref{eq_norm_zetapi1_minus_zetap12_bound} and using the induction hypothesis along with \eqref{eq_maxnorm_zetaT_bound} gives
\begin{eqnarray}
 \|\Qfunc{k-1}{T}{\pi_1}-\Qfunc{k-1}{T}{\pi_2}\|_{\rm max}  &\leq&  e^{|\gamma|\beta^{k-1} C} \left(\|\Qfunc{k}{T}{\pi_1}-\Qfunc{k}{T}{\pi_2}\|_{\rm max} + \| \Qfunc{k}{T}{\pi_2}\|_{\rm max}\cdot \| d_{k}-  f_{k}\|_\infty\right) \nonumber\\
&\leq& e^{|\gamma|\beta^{k-1} C}  e^{|\gamma| K(\beta^k-\beta^T)}\left(\sum_{i=k+1}^{T-1}\| d_{i}-  f_{i}\|_\infty + \| d_{k}-  f_{k}\|_\infty\right) \nonumber\\%\label{eq_induchypoth_result1}
&=& e^{|\gamma| K(\beta^{k-1}-\beta^T)} \sum_{i=k}^{T-1}\| d_{i}-  f_{i}\|_\infty.   
\end{eqnarray}
 
% \begin{eqnarray}
% \|\Qfunc{k-1}{T}{\pi_1}-\Qfunc{k-1}{T}{\pi_2}\|_{\rm max} & \leq & e^{|\gamma|\beta^{k-1} C} \left(\|\Qfunc{k}{T}{\pi_1}-\Qfunc{k}{T}{\pi_2}\|_{\rm max} + \| \Qfunc{k}{T}{\pi_2}\|_{\rm max}\cdot \| d_{k}-  f_{k}\|_\infty\right) \nonumber\\
% &\leq& e^{|\gamma|\beta^{k-1} C}  e^{|\gamma| K(\beta^k-\beta^T)}\left(\sum\limits_{i=k+1}^{T-1}\| d_{i}-  f_{i}\|_\infty + \| d_{k}-  f_{k}\|_\infty\right)\nonumber %\label{eq_induchypoth_result1}
% \\ 
% &=& e^{|\gamma| K(\beta^{k-1}-\beta^T)} \sum\limits_{i=k}^{T-1}\| d_{i}-  f_{i}\|_\infty. \nonumber
% \end{eqnarray}
In other words, \eqref{eq_maxnorm_zetaTOfpi1_min_norm_zetaTOfpi2_bound} holds for $t=k-1$. It now follows by induction that \eqref{eq_maxnorm_zetaTOfpi1_min_norm_zetaTOfpi2_bound} holds for all $t$ satisfying $0\leq t\leq T-1$. 
\eop

\textit{Proof of Theorem \ref{thm_unif_cont_of_USTOCost_prodtop}. }
 Let $\pi_{1}=\{d_{t}\}_{t=0}^{\infty}$ and $\pi_{2}=\{f_{t}\}_{t=0}^{\infty}$ be policies in $\Pi_{\rm MR}$. Recall from the proof of Lemma \ref{lemma_zeta_and_its_differ_bounds} that 
 \eqref{eq_norm_dTmintzetatmin1pi1_minus_fTmintzetatmin1pi2_bdd} holds for every $t$ satisfying $0\leq t\leq T-1$. Letting $t=0$ in \eqref{eq_norm_dTmintzetatmin1pi1_minus_fTmintzetatmin1pi2_bdd} and using \eqref{eq_maxnorm_zetaT_bound} and \eqref{eq_maxnorm_zetaTOfpi1_min_norm_zetaTOfpi2_bound} for $t=0$ gives 
 \begin{eqnarray}
     \|d_{0}\odot \Qfunc{0}{T}{\pi_1} - f_{0}\odot \Qfunc{0}{T}{\pi_2}\|_\infty &\leq & \|\Qfunc{0}{T}{\pi_1} - \Qfunc{0}{T}{\pi_2}\|_{\rm max} + \|\Qfunc{0}{T}{\pi_2}\|_{\rm max}\cdot\|d_{0} -f_{0}\|_\infty \nonumber \\ &\leq& e^{|\gamma| K(1-\beta^{T})}\sum\limits_{t=0}^{T-1}\|d_t-f_t\|_\infty. \label{eq_diff_d0zetaTpi1_f0zetaTpi2}
 \end{eqnarray}
 
 Next, starting from \eqref{eq_JpiT_In_Matrix_Form_eql_expgammar_d0_zetaTmin1} and using \eqref{eq_norminfAx_leq_norminfA_mult_norminfx}, \eqref{eq_normAschurprodB_leq_normmaxA_times_norminfB} and \eqref{eq_diff_d0zetaTpi1_f0zetaTpi2} along with $\|\one_{n\times 1}\|_{\infty}=1$ and $\|R\|_{\rm max}\leq C$,  we have 
 \begin{eqnarray}
  \|\RSCostTVec{\pi_1}{\gamma}{R}{T}-\RSCostTVec{\pi_2}{\gamma}{R}{T}\|_\infty &=&  \|(d_0\odot \Qfunc{0}{T}{\pi_1})\one_{n\times 1} -  (f_0\odot \Qfunc{0}{T}{\pi_2})\one_{n\times 1}\|_\infty\nonumber\\
&\leq& \| d_0\odot \Qfunc{0}{T}{\pi_1} - f_0\odot \Qfunc{0}{T}{\pi_2}\|_\infty \|\one_{n\times 1}\|_\infty \nonumber\\
&\leq&   e^{|\gamma| K(1-\beta^{T})}\sum_{t=0}^{T-1}\|d_t-f_t\|_\infty 
= e^{|\gamma| K(1-\beta^{T})}\sum_{t=0}^{T-1}\frac{\delta^t\|d_t-f_t\|_\infty}{\delta^t}\nonumber\\
&\leq&  e^{|\gamma| K(1-\beta^{T})}\mu(\pi_{1},\pi_{2})\sum_{t=0}^{T-1}\delta^{-t} \nonumber\\
&\leq&
e^{|\gamma| K(1-\beta^{T})}(1-\delta)^{-1}\cdot\delta^{-(T-1)}\mu(\pi_{1},\pi_{2}).
 \end{eqnarray}
  
%  \begin{eqnarray}
% \|\RSCostTVec{\pi_1}{\gamma}{R}{T}-\RSCostTVec{\pi_2}{\gamma}{R}{T}\|_\infty &=&  \|(d_0\odot \Qfunc{0}{T}{\pi_1})\one_{n\times 1} -  (f_0\odot \Qfunc{0}{T}{\pi_2})\one_{n\times 1}\|_\infty\nonumber\\
% &\leq &\| d_0\odot \Qfunc{0}{T}{\pi_1} - f_0\odot \Qfunc{0}{T}{\pi_2}\|_\infty \|\one_{n\times 1}\|_\infty\nonumber\\
% % &\leq&\|e^{\odot\gamma R}\|_{\rm max}\cdot\|d_0\odot \zeta_{T-1}^{\pi_1} - f_0\odot \zeta_{T-1}^{\pi_2}\|_\infty\nonumber\\
% &\leq &  e^{|\gamma| K(1-\beta^{T})}\sum\limits_{t=0}^{T-1}\|d_t-f_t\|_\infty 
% %\label{eq_diff_USTOCost_leq_d0zetaTpi1_min_f0zetaTpi2_2}\\
% = e^{|\gamma| K(1-\beta^{T})}\sum\limits_{t=0}^{T-1}\frac{\delta^t\|d_t-f_t\|_\infty}{\delta^t}\nonumber\\
% &\leq & e^{|\gamma| K(1-\beta^{T})}\mu(\pi_{1},\pi_{2})\sum\limits_{t=0}^{T-1}\delta^{-t} 
% \leq
% e^{|\gamma| K(1-\beta^{T})}\frac{\delta^{-(T-1)}}{1-\delta}\mu(\pi_{1},\pi_{2}). \nonumber
% %\label{eq_Jpi1T_min_Jpi2T_leq_mupi1pi2_2}
% \end{eqnarray}
This proves the first assertion of the theorem. 

To prove the second assertion, recall from Lemma \ref{lem_RSCost_finhor_cgsto_infinhor} that $\RSCostTVec{\pi}{\gamma}{R}{T}$ converges to $\RSCostVec{\pi}{\gamma}{R}$ uniformly in $\pi$ as $T\rightarrow \infty$. As shown above, for each $T$, the function $\pi\mapsto \RSCostTVec{\pi}{\gamma}{R}{T}$ is Lipschitz continuous, and hence uniformly continuous, on the metric space $(\Pi_{\rm MR},\mu)$. The second assertion now follows from the fact 
%uniform limit theorem
% \cite[Theorem 21.6]{JamesMunkres} and \cite[Prop 23, Ch 9, Pg 202]{RoyFtz} or \cite[Theorem 4.6, Pg 234]{JDug}
%\cite[Thm. 4.6, p. 234]{JDug},
%\cite[Thm. 21.6]{JamesMunkres}, \cite[Prop. 23, p. 202]{RoyFtz},  
%which states 
that the uniform limit of a sequence of uniformly continuous functions on a compact metric space is uniformly continuous. Refer the following literature by \cite[Thm. 4.6, p. 234]{JDug},
\cite[Thm. 21.6]{JamesMunkres}, \cite[Prop. 23, p. 202]{RoyFtz} for more details.  \hspace{1cm}
\eop

\textit{Proof of Lemma \ref{lem_fesblereg_prop}.} 
We begin by showing that the sets  $\FeasibleSetOrgnal$, $\FeasibleSetTrunBlw{T}$ and $\FeasibleSetTrunAbv{T}$ are compact, which also serves to prove (i). Fix $T\in \{1,2,\ldots\}$. 
%If the sets $\FeasibleSetOrgnal$, $\FeasibleSetTrunBlw{T}$ and $\FeasibleSetTrunAbv{T}$ are empty, then they are compact trivially. Suppose they are non-empty. We can observe that
By theorems \ref{thm_contuty_lin_disc_cost_prod_toplogy} and  \ref{thm_unif_cont_of_USTOCost_prodtop}, each of the sets defined in \eqref{eq_Infin_Hor_Lin_Const}-\eqref{eq_defn_RSTrunAbv} is a finite intersection of inverse images of closed intervals
% \footnote{\label{footnt_bddcost} By finite state-action space and discount factor $\beta\in (0,1)$, the cost functions involved in defining constraint sets \eqref{eq_Infin_Hor_Lin_Const}-\eqref{eq_defn_RSTrunAbv} are uniformly upper and lower bounded.} 
under continuous functions, and are hence closed subsets of the compact space $\Pi_{\rm MR}$. This  proves statement \ref{itm_compact_fsblsets} as well as compactness in statements \ref{itm_FT_nondec} and \ref{itm_barFT_noninc}. 

\ref{itm_FT_nondec} For each $i\in\{1,\ldots, M\}$ and $T\in\{1,2,\ldots\}$, define $\FeasibleSetLDTruncBlwAti{i}{T} :=\{\pi\in\Pi_{\rm MR}:\LDCostT{\pi}{C_i}{T}\leq b_i-K\beta^T \}$ and $\FeasibleSetLDBlwAti{i} :=\{\pi\in\Pi_{\rm MR}:\LDCost{\pi}{C_i}\leq b_i \}$. Similarly, define $\FeasibleSetRSTruncBlwAtj{T} :=\left\{\pi\in\Pi_{\rm MR}: \RSCostT{\pi}{\gamma}{\hat{C}_j}{T}\leq \frac{\hat{b}_j}{K_{T}}\right\}$ and $\FeasibleSetRSBlwAtj:=\{\pi\in\Pi_{\rm MR}:\RSCost{\pi}{\gamma}{\hat{C}_j}\leq \hat{b}_j \}$ for each $j\in\{1,\ldots,\hat{M}\}$ and $T\in\{1,2,\ldots\}$. We can observe from
\eqref{eq_Infin_Hor_Lin_Const}, \eqref{eq_Infin_Hor_RSConst}, \eqref{eq_defn_LDTrunBlw} and  \eqref{eq_defn_RSTrunBlw}  that 
$\FeasibleSetLD=\cap_{i=1}^{M}\FeasibleSetLDBlwAti{i}$, $\FeasibleSetRS=\cap_{j=1}^{\hat{M}}\FeasibleSetRSBlwAtj$
and, for each $T=1,2\ldots$,
$ \FeasibleSetLDTruncBlw{T}=\cap_{i=1}^{M}\FeasibleSetLDTruncBlwAti{i}{T}$ and $\FeasibleSetRSTruncBlw{T}=\cap_{j=1}^{\hat{M}} \FeasibleSetRSTruncBlwAtj{T}$. It is now easy to see from the definitions of $\FeasibleSetTrunBlw{T}$ and $\mathcal{F}$ that statement \ref{itm_FT_nondec} follows if the following two claims are established.
\begin{description} 
\item[Claim 1:] For each $i\in\{1,\ldots, M\}$, $\{\FeasibleSetLDTruncBlwAti{i}{T}\}_{T=1}^{\infty}$ is a non-decreasing sequence of nested sets contained in $\FeasibleSetLDBlwAti{i}$. 
\item[Claim 2:] For each $j\in\{1,\ldots,\hat{M}\}$, $\{\FeasibleSetRSTruncBlwAtj{T}\}_{T=1}^{\infty}$ is a non-decreasing sequence of nested sets contained in $\FeasibleSetRSBlwAtj$
\end{description} 

To prove Claim 1,  note that 
\begin{equation}
    \left|\LDCostT{\pi}{C_i}{T+1}-\LDCostT{\pi}{C_i}{T}\right|=|\E^\pi_x[\beta^{T}R(X_{T},A_{T})]|\leq \beta^{T}C,
    \label{eq_onesteplcostdiff}
\end{equation}
holds for every $i\in\{1,\ldots,M\}$, $T\in\{1,2,\ldots\}$ and $\pi\in\Pi_{\rm MR}$.
Hence, if $\LDCostT{\pi}{C_i}{T}\leq b_{i}-\beta^{T}K$ holds, then 
$\LDCostT{\pi}{C_i}{T+1}\leq \LDCostT{\pi}{C_i}{T}+\beta^{T}C \leq b_{i}-\beta^{T}K+\beta^{T}C=b_{i}-\beta^{T+1}K$
% $$\LDCostT{\pi}{C_i}{T+1}\leq \LDCostT{\pi}{C_i}{T}+\beta^{T}C \leq b_{i}-\beta^{T}K+\beta^{T}C=b_{i}-\beta^{T+1}K$$ 
also holds (recall that $C=K(1-\beta)$). Thus, $\FeasibleSetLDTruncBlwAti{i}{T}\subseteq \FeasibleSetLDTruncBlwAti{i}{T+1}$ for each $T$ and $i$.  

Furthermore, applying \eqref{eq_onesteplcostdiff} repeatedly yields 
$\left|\LDCostT{\pi}{C_i}{T+k}-\LDCostT{\pi}{C_i}{T}\right|\leq \beta^{T}(1+\beta+\cdots+\beta^{k-1})C,$
% $$\left|\LDCostT{\pi}{C_i}{T+k}-\LDCostT{\pi}{C_i}{T}\right|\leq \beta^{T}(1+\beta+\cdots+\beta^{k-1})C,$$
for each $k\geq 1$. 

Lemma \ref{lem_linCost_finhor_cgsto_infinhor} implies that $\lim_{k\rightarrow \infty}\LDCostT{\pi}{C_i}{T+k}=\LDCost{\pi}{C_i}$. Letting $k\rightarrow \infty$ in the last inequality thus leads to 
\begin{equation}
    |\LDCost{\pi}{C_i}-\LDCostT{\pi}{C_i}{T}|\leq \beta^{T}K.
    \label{eq_infsteplcostdiff}
\end{equation}
Hence, if $\LDCostT{\pi}{C_i}{T}\leq b_{i}-\beta^{T}K$ holds, then $\LDCost{\pi}{C_i}\leq \LDCostT{\pi}{C_i}{T}+\beta^{T}K\leq b_{i}$ holds as well. It immediately follows that $\FeasibleSetLDTruncBlwAti{i}{T}\subseteq \FeasibleSetLDBlwAti{i}$. This proves Claim 1. 

 To prove Claim 2, fix $j\in\{1,\ldots,\check{M}\}$ and $\pi\in\Pi_{\rm MR}$. Note that  $e^{-|\gamma| \beta^T C}\leq e^{\gamma\beta^{T}R(s,a)}\leq e^{|\gamma| \beta^T C}$ holds for every $T\in\{1,2,\ldots,\}$ and every $(s,a)\in\mathcal{S}\times \mathcal{A}$. The last inequality yields
%\begin{eqnarray}
%e^{-|\gamma|\beta^T C} \RSCostT{\pi}{\gamma}{\hat{C}_j}{T}&\leq& \E^\pi_x\left[ e^{\gamma\sum_{t=0}^{T-1}\beta^t \hat{C}_j(X_t,A_t)} e^{\gamma \beta^T \hat{C}_j(X_T,A_T)}\right]\nonumber\\
%&=&\RSCostT{\pi}{\gamma}{\hat{C}_j}{T+1}\leq e^{|\gamma|\beta^T C} \RSCostT{\pi}{\gamma}{\hat{C}_j}{T}
%\end{eqnarray}
\begin{eqnarray}
e^{-|\gamma| \beta^T C} \RSCostT{\pi}{\gamma}{\hat{C}_j}{T}\leq \RSCostT{\pi}{\gamma}{\hat{C}_j}{T+1}
\leq  e^{|\gamma| \beta^T C}\RSCostT{\pi}{\gamma}{\hat{C}_j}{T}, \label{eq_ratio_RSCOSTs_bdd}
\end{eqnarray}
for every $T=1,2,\ldots$. 
It is easy to see from \eqref{eq_ratio_RSCOSTs_bdd}  that if $\RSCostT{\pi}{\gamma}{\hat{C}_j}{T}\leq \frac{\hat{b}_j}{K_{T}}$ holds for some $T$, then $\RSCostT{\pi}{\gamma}{\hat{C}_j}{T+1}\leq e^{|\gamma| \beta^T C} \RSCostT{\pi}{\gamma}{\hat{C}_j}{T} \leq e^{|\gamma| \beta^T C} \frac{\hat{b}_j}{K_{T}} =\frac{\hat{b}_j}{K_{T+1}} $ also holds. We immediately conclude that $\FeasibleSetRSTruncBlwAtj{T}\subseteq \FeasibleSetRSTruncBlwAtj{T+1}$ for each $T$ and $j$.   %\FeasibleSetRSTruncBlwAtj{T}

To complete the proof of Claim 2, fix $T\in\{1,2,\ldots\}$. Applying  \eqref{eq_ratio_RSCOSTs_bdd} repeatedly yields 
$e^{-|\gamma|\beta^{T}K (1-\beta^k)} \RSCostT{\pi}{\gamma}{\hat{C}_j}{T}\leq \RSCostT{\pi}{\gamma}{\hat{C}_j}{T+k}\leq e^{|\gamma|\beta^{T}K (1-\beta^k)}\RSCostT{\pi}{\gamma}{\hat{C}_j}{T}$
% $$e^{-|\gamma|\beta^{T}K (1-\beta^k)} \leq \frac{\RSCostT{\pi}{\gamma}{\hat{C}_j}{T+k}}{\RSCostT{\pi}{\gamma}{\hat{C}_j}{T}}\leq e^{|\gamma|\beta^{T}K (1-\beta^k)}$$ 
for each $k\geq 1$. Lemma \ref{lem_RSCost_finhor_cgsto_infinhor} implies that $\lim_{k\to\infty}\RSCostT{\pi}{\gamma}{\hat{C}_j}{T+k} =\RSCost{\pi}{\gamma}{\hat{C}_j}$. Letting $k\to\infty$ in the last inequality thus leads to 
% \begin{equation}
%  K_{T}^{-1}\RSCostT{\pi}{\gamma}{\hat{C}_j}{T}=e^{-|\gamma| \beta^{T}K}\RSCostT{\pi}{\gamma}{\hat{C}_j}{T} \leq \RSCost{\pi}{\gamma}{\hat{C}_j}\leq e^{|\gamma| \beta^{T}K}\RSCostT{\pi}{\gamma}{\hat{C}_j}{T}=K_{T}\RSCostT{\pi}{\gamma}{\hat{C}_j}{T}.\label{eq_RScostT_div_RScost_bdd}
% \end{equation}
\begin{equation}
 \frac{1}{K_{T}}=e^{-|\gamma| \beta^{T}K} \leq \frac{\RSCost{\pi}{\gamma}{\hat{C}_j}}{\RSCostT{\pi}{\gamma}{\hat{C}_j}{T}}\leq e^{|\gamma| \beta^{T}K}=K_{T}.\label{eq_RScostT_div_RScost_bdd}
\end{equation}
Hence, if $\RSCostT{\pi}{\gamma}{\hat{C}_j}{T}\leq \frac{\hat{b}_j}{K_{T}}$ then
$\RSCost{\pi}{\gamma}{\hat{C}_j}\leq K_{T} \RSCostT{\pi}{\gamma}{\hat{C}_j}{T}\leq \hat{b}_j$. It immediately follows that $\FeasibleSetRSTruncBlwAtj{T} \subseteq\FeasibleSetRSBlwAtj$ for all $T$ and $j$. Claim 2 now follows. This proves  \ref{itm_FT_nondec} of the lemma. 
% This completes the proof of  statement \ref{itm_FT_nondec} of the lemma. 

\ref{itm_barFT_noninc} For each $i\in\{1,2,\ldots M\}$ and $T\in\{1,2,\ldots\}$, define, $\FeasibleSetLDTruncAbvAti{T}=\{\pi\in\Pi_{\rm MR}:\LDCostT{\pi}{C_i}{T}\leq b_i + K\beta^T \}$. Similarly for each $j\in\{1,2,\ldots \hat{M}\}$ and $T\in\{1,2,\ldots\}$, define, $\FeasibleSetRSTruncAbvAtj{T}=\{\pi\in\Pi_{\rm MR}:\RSCostT{\pi}{\gamma}{\hat{C}_j}{T}\leq \hat{b}_j K_{T} \}$. 
We can observe from \eqref{eq_defn_LDTrunAbv} and \eqref{eq_defn_RSTrunAbv} that for each $T=1,2,\ldots$, $\FeasibleSetLDTruncAbv{T}=\cap_{i=1}^{M}\FeasibleSetLDTruncAbvAti{T}$ and $\FeasibleSetRSTruncAbv{T}=\cap_{j=1}^{\hat{M}}\FeasibleSetRSTruncAbvAtj{T}$.
It is now easy to see from the definitions of $\FeasibleSetTrunAbv{T}$ and $\FeasibleSetOrgnal$ that statement \ref{itm_barFT_noninc} follows if the following two claims are established.
\begin{description} 
\item[Claim 3:] For each $i\in\{1,\ldots, M\}$, $\{\FeasibleSetLDTruncAbvAti{T}\}_{T=1}^{\infty}$ is a non-increasing sequence of nested sets converging to $\FeasibleSetLDBlwAti{i}$. 
\item[Claim 4:] For each $j\in\{1,\ldots,\hat{M}\}$, $\{\FeasibleSetRSTruncAbvAtj{T}\}_{T=1}^{\infty}$ is a non-increasing sequence of nested sets converging to $\FeasibleSetRSBlwAtj$.
\end{description} 
To prove Claim 3, fix $i\in\{1,\ldots,M\}$ and $\pi\in\Pi_{\rm MR}$. It can be seen from \eqref{eq_onesteplcostdiff} that, if $\LDCostT{\pi}{C_i}{T+1}\leq b_i + K\beta^{T+1}$ holds for some $T$, then 
$\LDCostT{\pi}{C_i}{T}\leq \LDCostT{\pi}{C_i}{T+1} + C\beta^{T}  \leq b_i + K\beta^{T+1} + C\beta^{T} = b_i +K\beta^{T}$
% $$\LDCostT{\pi}{C_i}{T}\leq \LDCostT{\pi}{C_i}{T+1} + C\beta^{T}  \leq b_i + K\beta^{T+1} + C\beta^{T} = b_i +K\beta^{T}$$ 
also holds. Thus, $\FeasibleSetLDTruncAbvAti{T+1}\subseteq \FeasibleSetLDTruncAbvAti{T}$ for each $T$ and $i$.

To complete the proof of Claim 3, fix $T\in\{1,2,\ldots\}$. It is easy to see from \eqref{eq_infsteplcostdiff} that if $\LDCost{\pi}{C_i}\leq b_i$ holds, then $\LDCostT{\pi}{C_i}{T}\leq \LDCost{\pi}{C_i}  + K\beta^T \leq b_i + K\beta^T$ also holds. It immediately follows that $\FeasibleSetLDBlwAti{i}\subseteq \FeasibleSetLDTruncAbvAti{T}$ for any fixed $T$, implying that $\FeasibleSetLDBlwAti{i}\subseteq \cap_{T}\FeasibleSetLDTruncAbvAti{T}$. 
On the other hand, if $\LDCostT{\pi}{C_i}{T}\leq b_i + K\beta^T$ holds for all $T=1,2,\ldots$, then letting $T\to\infty$ and applying Lemma \ref{lem_linCost_finhor_cgsto_infinhor} shows that $\LDCost{\pi}{C_i} \leq b_i$ holds. We can immediately conclude that $\cap_{T}\FeasibleSetLDTruncAbvAti{T}\subseteq\FeasibleSetLDBlwAti{i}$. Thus Claim 3 is established.

To prove Claim 4, fix $j\in\{1,\ldots,\hat{M}\}$ and $\pi\in\Pi_{\rm MR}$. It can be seen from \eqref{eq_ratio_RSCOSTs_bdd} that if $\RSCostT{\pi}{\gamma}{\hat{C}_j}{T+1}\leq \hat{b}_j K_{T+1} $ for some $T$, then  $\RSCostT{\pi}{\gamma}{\hat{C}_j}{T}\leq \RSCostT{\pi}{\gamma}{\hat{C}_j}{T+1}e^{|\gamma|C\beta^T} $ $  \leq \hat{b}_j K_{T+1}e^{|\gamma|C\beta^T} = \hat{b}_j K_{T}$
% $$\RSCostT{\pi}{\gamma}{\hat{C}_j}{T}\leq \RSCostT{\pi}{\gamma}{\hat{C}_j}{T+1}e^{|\gamma|C\beta^T} \leq \hat{b}_j K_{T+1}e^{|\gamma|C\beta^T} = \hat{b}_j K_{T}$$ 
also holds. We therefore have $\FeasibleSetRSTruncAbvAtj{T+1}\subseteq \FeasibleSetRSTruncAbvAtj{T}$ for each $T\in\{1,2,\ldots\}$.

To complete the proof of Claim 4, fix $T\in\{1,2,\ldots\}$. It is easy to see from \eqref{eq_RScostT_div_RScost_bdd} that if $\RSCost{\pi}{\gamma}{\hat{C}_j}\leq \hat{b}_j$ holds, then $\RSCostT{\pi}{\gamma}{\hat{C}_j}{T}\leq \RSCost{\pi}{\gamma}{\hat{C}_j}K_{T}\leq \hat{b}_jK_{T}$ also holds. Hence, $\FeasibleSetRSBlwAtj\subseteq \FeasibleSetRSTruncAbvAtj{T}$ for every  $T\in\{1,2,\ldots\}$. As a consequence, it follows that    $\FeasibleSetRSBlwAtj\subseteq\cap_{T} \FeasibleSetRSTruncAbvAtj{T}$. On the other hand, if $\RSCostT{\pi}{\gamma}{\hat{C}_j}{T}\leq  \hat{b}_j K_{T}$ holds for all $T$, then letting $T\rightarrow \infty$ and applying Lemma \ref{lem_RSCost_finhor_cgsto_infinhor} shows that $\RSCost{\pi}{\gamma}{\hat{C}_j}\leq \hat{b}_j$ also holds. Hence $\cap_{T} \FeasibleSetRSTruncAbvAtj{T}\subseteq\FeasibleSetRSBlwAtj$. Thus Claim 4 is established and \ref{itm_barFT_noninc} is proved.
 \eop

% closed subsets of compact set $\Pi_{\rm MR}$ (Theorem \ref{thm_metric_on_Pi_metrizes_ProdTop}), because each of them equals finite intersections of continuous inverse image of closed intervals of the form $(-\infty,a]\subset \R$. Thus, making them compact. Since, the feasible regions $\FeasibleSetOrgnal$, $\FeasibleSetTrunBlw{T}$ and $\FeasibleSetTrunAbv{T}$ are finite intersections of sets defined in  \eqref{eq_Infin_Hor_Lin_Const}-\eqref{eq_defn_RSTrunAbv}, makes each of them compact, completing the proof of \ref{itm_compact_fsblsets} .

%%%%%%%%%%%%%%%%%%%%%%%%%% Proof of Theorem \ref{lem_fesblereg_prop} ends %%%%%%%%%%%%%%%%%%%%%%%%%% 
%\line(1,0){360}

\section{Proofs for Results from Subsection \ref{thm2ssec} }\label{sec_apndxB}%
%%%%%%%%%%%%%%%%%%%%% Proof of Theorem \ref{prop_generalised_convgce_offunctions_over_seqofsets} begins %%%%%%%%%%%%%%%%%%%% 

\textit{Proof of Proposition \ref{prop_generalised_convgce_offunctions_over_seqofsets}.} Choose $\epsilon>0$. By uniform convergence, there exists $L_{1}\in\Z^{+}$ such that $|f(z)-f_{T}(z)|<\frac{\epsilon}{3}$ for all $z\in\mathcal{V}$ and $T>L_{1}$. For every $z\in\mathcal{V}$ and $T>L_{1}$, we have  $\inf_{y\in\mathcal{V}}f(y)\leq f(z)\leq f_{T}(z)+\epsilon/3.$
Thus, for every $T>L_{1}$,  $f_{T}$ is bounded below on $\mathcal{V}$, and hence on $\mathcal{V}_{T}\subseteq \mathcal{V}$ as well. The last inequality now yields 
\begin{equation}
\inf_{z\in\mathcal{V}}f(z)\leq \inf_{z\in\mathcal{V}}f_{T}(z)+\frac{\epsilon}{3}< \inf_{z\in\mathcal{V}_{T}}f_{T}(z)+\epsilon,
\label{ineq2}
\end{equation}
for every $T>L_{1}$.

Next, let $z_{1}\in\mathcal{V}$ be such that $f(z_{1})<\inf_{z\in\mathcal{V}}f(z)+\frac{\epsilon}{3}$. Define  the set $\mathcal{O}:=\left\{z\in\mathcal{U}: |f(z)-f(z_{1})|<\frac{\epsilon}{3}\right\}$, and note that $z_{1}\in\mathcal{O}$. By the continuity of $f$, $\mathcal{O}$ is an open set in $\mathcal{U}$. Since $\mathcal{V}$ is the closure of the union $\cup_{T=1}^{\infty}\mathcal{V}_{T}$ and $z_{1}\in\mathcal{V}$, it follows from the definition of closure that every open neighborhood of $z_{1}$ has a nonempty intersection with the union  $\cup_{T=1}^{\infty}\mathcal{V}_{T}$. In particular, we may conclude that there exists $L_{2}\in\Z^{+}$ and $z_{2}\in\mathcal{U}$ such that $z_{2}\in\mathcal{O}\cap\mathcal{V}_{L_{2}}$. Since the sequence of sets $\{\mathcal{V}_{T}\}_{T=1}^{\infty}$  is non-decreasing, it follows that $z_{2}\in\mathcal{O}\cap\mathcal{V}_{T}$ for all $T>L_{2}$. For every $T>\max\{L_{1},L_{2}\}$, we now have 
\begin{equation}
    \inf_{z\in\mathcal{V}_{T}}f_{T}(z)\leq f_{T}(z_{2})<f(z_{2})+\frac{\epsilon}{3}<f(z_{1})+\frac{2\epsilon}{3}<\inf_{z\in\mathcal{V}}f(x)+\epsilon.
    \label{ineq3}
\end{equation}
Note that the first, second, third and last inequalities in (\ref{ineq3}) follow from $z_{2}\in\mathcal{V}_{T}$, our choice of $L_{1}$, $z_{2}\in\mathcal{O}$, and our choice of $z_{1}$, respectively. 

The inequalities (\ref{ineq2}) and (\ref{ineq3}) together imply that $$|\inf_{z\in\mathcal{V}}f(z)-\inf_{z\in\mathcal{V}_{T}}f_{T}(x)|<\epsilon, \mbox{ for all }T>\max\{L_{1},L_{2}\}.$$ Since we chose $\epsilon>0$ arbitrarily, \eqref{eq_limntoinf_infofgnoverFn_eqls_infofgoverF} follows. 
\eop

\textit{Proof of Proposition \ref{prop_cl_uninovrT_intscovri_FTi_eql_intscovrFi}.}
 For each  $z\in\mathcal{H}$, define $I(z):=\{i:1\leq i\leq N, g_i(z)<B_i\}$ and $J(z):=\{i:1\leq i\leq N,g_i(z)=B_i \}$, and note that $I(z)\cup J(z)=\{1,\ldots,N\}$.  For every $i\in I(z)$, our convergence assumptions on the sequences $\{g_{T,i}\}_{T=1}^{\infty}$ and $\{B_{T,i}\}_{T=1}^{\infty}$ imply that the sequence $\{g_{T,i}(z)-B_{T,i}\}_{T=1}^{\infty}$ converges to $g_{i}(z)-B_{i}<0$. Consequently, for every $i\in I(z)$, there exists $\tau_{i}$ such that $g_{T,i}(z)-B_{T,i}<0$ for all $T>\tau_{i}$. On letting $\tau(z)=\max_{i\in I(z)}\tau_{i}$, it follows that $z\in \cap_{i\in I(z)}\mathcal{G}_{T,i}$ for all $T>\tau(z)$.
 
To prove that $\mathcal{H}\subseteq \cl\left(\cup_{T=1}^{\infty}\mathcal{H}_{T}\right)$, choose $y\in\mathcal{H}$.  Consider the case where $J(y)=\varnothing$. In this case, $I(y)=\{1,\ldots,N\}$, and we conclude from the previous paragraph that $y\in \cap_{i=1}^{N}\mathcal{G}_{T,i}=\mathcal{H}_{T}$ for every $T>\tau(y)$. In particular, $y\in\cup_{T=1}^{\infty}\mathcal{H}_{T}\subseteq \cl\left(\cup_{T=1}^{\infty}\mathcal{H}_{T}\right)$. 
 
Next, suppose $0$ is not a local minimum value of $h$. Consider the case where $J(y)\neq \varnothing$, and let $\mathcal{V}\subseteq \mathcal{U}$ be an open set containing $y$. Recall  that, for every $j\in J(y)$, $g_{j}(y)=B_{j}$. Hence it follows from $J(y)\neq \varnothing$ that  $h(y)=0$. Since $0$ is not a local minimum value of $h$, there exists $z\in \mathcal{V}$ such that $h(z)<0$. The condition $h(z)<0$ implies that $z\in\mathcal{H}$ and $J(z)=\varnothing$. Applying the arguments of the previous paragraph to the point $z$ shows that $z\in\cup_{T=1}^{\infty}\mathcal{H}_{T}$. We have thus shown that the arbitrarily chosen open neighborhood  $\mathcal{V}$ of $y\in\mathcal{H}$ has a nonempty intersection with $\cup_{T=1}^{\infty}\mathcal{H}_{T}$. It follows that $y\in\cl\left(\cup_{T=1}^{\infty}\mathcal{H}_{T}\right)$. 
 
Since $y\in\mathcal{H}$ above was chosen arbitrarily, it follows from the conclusions of the previous two paragraphs that $\mathcal{H}\subseteq \cl\left(\cup_{T=1}^{\infty}\mathcal{H}_{T}\right)$. To show the reverse inclusion, we first deduce from elementary set-theoretic arguments that  $\cup_{T=1}^{\infty}\mathcal{H}_{T}=\cup_{T=1}^{\infty}\cap_{i=1}^{N}\mathcal{G}_{T,i}\subseteq \cap_{i=1}^{N}\cup_{T=1}^{\infty}\mathcal{G}_{T,i}$. Next, we have 
%  \begin{eqnarray}
%  \cl\left(\bigcup\limits_{T=0}^{\infty}\mathcal{H}_{T}\right)&\subseteq&
%  \cl\left(\bigcap\limits_{i=1}^{N}\bigcup\limits_{T=0}^{\infty}\mathcal{G}_{T,i}\right) \nonumber \\
%  &\subseteq& \bigcap\limits_{i=1}^{N}\cl\left(\bigcup\limits_{T=0}^{\infty}\mathcal{G}_{T,i}\right)\subseteq \bigcap\limits_{i=1}^{N}\cl(\mathcal{G}_{i})=\bigcap\limits_{i=1}^{N}\mathcal{G}_{i}=\mathcal{H}. %\label{eq_cls_un_ovrT_FTstr_sbst_F_1}
%  \end{eqnarray}
%  {\color{red} Check the equation labels before deleting the above black equations where $T=0,1,$
\begin{eqnarray}
 \cl\left(\cup_{T=1}^{\infty}\mathcal{H}_{T}\right)&\subseteq&
 \cl\left(\cap_{i=1}^{N}\cup_{T=1}^{\infty}\mathcal{G}_{T,i}\right) 
 \nonumber \\
 &\subseteq& \cap_{i=1}^{N}\cl\left(\cup_{T=1}^{\infty}\mathcal{G}_{T,i}\right)\subseteq \cap_{i=1}^{N}\cl(\mathcal{G}_{i})=\cap_{i=1}^{N}\mathcal{G}_{i}=\mathcal{H}. \label{eq_cls_un_ovrT_FTstr_sbst_F_1}
 \end{eqnarray}
The first inclusion in \eqref{eq_cls_un_ovrT_FTstr_sbst_F_1} follows from the fact that the closure operation preserves inclusion, the second from the fact that the closure of an intersection is contained in the intersection of the closures, and the third from our assumption that $\mathcal{G}_{T,i}\subseteq \mathcal{G}_{i}$ for each $T$ and $i$. The first equality in \eqref{eq_cls_un_ovrT_FTstr_sbst_F_1} follows from the fact that, for each $i$,  $\mathcal{G}_{i}$ is the inverse image of the closed interval  $(-\infty,B_{i}]$ under the continuous function $g_{i}$,  and hence closed. This completes the proof.\eop
 
%%%%%%%%%%%%%%%%%%%%%%%
%%%%%%%%%%%%%%%%%%% Proof of Theorem \ref{prop_cl_uninovrT_intscovri_FTi_eql_intscovrFi} ends %%%%%%%%%%%%%%%%%%%%%

%\hspace{-0.4 cm}\line(1,0){430}
\section{Proofs for Results from Subsection \ref{sec_thm3ssec} }\label{sec_apndxC}

\textit{Proof of Proposition \ref{prop_lim_of_inf_overGT_eql_inf_overG}.}
 Choose $\epsilon>0$. By uniform convergence, there exists $N>0$ such that $|g_{T}(y)-g(y)|\leq \epsilon$ for all $y\in\mathcal{G}_{1}$ and all $T\geq N$. Since a continuous function attains its infimum on a compact set, for every $T\geq 1$, there exists  $y_T\in \mathcal{G}_{T}$ such that $g_{T}(y_{T})= \inf_{y\in \mathcal{G}_T} g_{T}(y)$.
 
 For every $T\geq N$ and every  $z\in\mathcal{G}\subseteq\mathcal{G}_{T}\subseteq\mathcal{G}_{1}$, we have $g_{T}(y_{T})\leq g_{T}(z)\leq g(z)+\epsilon$. It follows that 
 \begin{equation}
 \limsup_{T\rightarrow \infty}g_{T}(y_{T})\leq \inf_{y\in\mathcal{G}}g(y)+\epsilon.
 \label{eq_limit_GT_infm_leq_inf_G}
 \end{equation}
 
There exists a  strictly increasing  sequence $\{T_{l}\}_{l=1}^{\infty}$ of integers such that the sequence  $\{g_{T_{l}}(y_{T_{l}})\}_{l=1}^{\infty}$ converges to  $\liminf_{T\rightarrow \infty}g_{T}(y_{T})$. We may choose  $L$ such that $T_{L}>N$ and $g_{T_{l}}(y_{T_{l}})\leq \liminf_{T\rightarrow \infty}g_{T}(y_{T})+\frac{\epsilon}{2}$ for all $l\geq L$. 
 
Next, for each $T\geq 1$, define the set $\mathcal{D}_T:= \{y_t:t\geq T\}\subseteq \mathcal{G}_T$.  Note that, for every $T\geq 1$, $\cl(\mathcal{D}_T)$ is a closed subset of the compact set $\mathcal{G}_T$, and hence compact. The sequence of sets $\{\cl(\mathcal{D}_T)\}_{T=1}^\infty$ is thus a non-increasing sequence of non-empty compact sets. By the Frechet Intersection, Theorem  given by \cite[p. 236]{RoyFtz}, \cite[p. 253]{JDug}, 
%\cite[Ch. 11, prob. 40, p. 236]{RoyFtz}, \cite[Ch. XI, prob. 2, p. 253]{JDug}
$\cap_{T\geq1} \cl(\mathcal{D}_T)$ is non-empty. Also, note that  $\cap_{T\geq1} \cl(\mathcal{D}_T)\subseteq \cap_{T\geq1} \mathcal{G}_T=\mathcal{G}$.

Next, choose, $z\in\cap_{T\geq1} \cl(\mathcal{D}_T)\subseteq\mathcal{G}$. By continuity, there exists an open neighbourhood $\mathcal{V}\subseteq \mathcal{U}$ of $z$ such that $g(y)>g(z)-\frac{\epsilon}{2}$ for all $y\in \mathcal{V}$. Since $z\in\cap_{T\geq1} \cl(\mathcal{D}_T)\subseteq\cl(\mathcal{D}_{T_{L}})$, the set $\mathcal{V}\cap \mathcal{D}_{T_{L}}$ is nonempty. In other words, there exists $l\geq L$ such that 
%$T_{l}\geq T_{L}>N$ and  
$y_{T_l}\in \mathcal{V}$. Note that $T_{l}\geq T_{L}>N$.  We now have,
\begin{equation}
    \inf_{y\in\mathcal{G}}g(y)-2\epsilon\leq g(z)-2\epsilon\leq g(y_{T_{l}})-\frac{3\epsilon}{2}\leq g_{T_{l}}(y_{T_{l}})-\frac{\epsilon}{2}\leq \liminf_{T\rightarrow \infty}g_{T}(y_{T}). 
    \label{eq_limit_GT_infm_leq_inf_G2}
\end{equation} 
Comparing \eqref{eq_limit_GT_infm_leq_inf_G} and \eqref{eq_limit_GT_infm_leq_inf_G2} shows that 
$\inf_{y\in\mathcal{G}}g(y)-2\epsilon \leq \liminf_{T\rightarrow \infty}g_{T}(y_{T}) \leq \limsup_{T\rightarrow \infty}g_{T}(y_{T}) \leq \inf_{y\in\mathcal{G}}g(y)+\epsilon$.
% \[\inf_{y\in\mathcal{G}}g(y)-2\epsilon \leq \liminf_{T\rightarrow \infty}g_{T}(y_{T}) \leq \limsup_{T\rightarrow \infty}g_{T}(y_{T}) \leq \inf_{y\in\mathcal{G}}g(y)+\epsilon.\]
Since $\epsilon>0$ was chosen arbitrarily, we conclude that $\inf_{y\in\mathcal{G}}g(y)=\liminf_{T\rightarrow \infty}g_{T}(y_{T})$ $ = \limsup_{T\rightarrow \infty}g_{T}(y_{T})$. In other words, \eqref{eq_limit_GT_infm_inf_G} holds. This completes the proof.
\eop

 \section{Counter Example}\label{appendixD}

\begin{ex}
Consider an MDP with $\mathcal{S}=\{s\}$ and $\mathcal{A}=\{a_{1},a_{2}\}$. Define two immediate cost functions $C_{i}:\mathcal{S}\times \mathcal{A}\rightarrow \R$, $i\in\{1,2\}$ by $C_{1}(s,a_{1})=1$,  $C_{1}(s,a_{2})=0$, and $C_{2}:=1-C_{1}$. Let $\beta=2^{-1}$. Note that $C=1$, and thus $K=C/(1-\beta)=2$. 

Next, for each $i\in\{1,2\}$, recall the definitions \eqref{eq_LDCost_FinandinfHor} of the infinite-horizon standard discounted costs $\pi\mapsto\LDCostVec{\pi}{C_{i}}$ and the finite-horizon standard discounted costs $\pi\mapsto \LDCostTVec{\pi}{C_{i}}{T}$ for $T\in\{1,2,\ldots\}$. It is easy to see that for all $\pi\in\Pi_{\rm MR}$ and $T\in\{1,2,\ldots\}$,
\begin{equation}
\LDCostTVec{\pi}{C_{1}}{T}(s)+\LDCostTVec{\pi}{C_{2}}{T}(s)=2-\frac{1}{2^{T-1}},\ \mbox{ and } \ \LDCostVec{\pi}{C_{1}}(s)+\LDCostVec{\pi}{C_{2}}(s)=2.\label{eq_ctrex2}
\end{equation}

Next, define $b_{1}=b_{2}=1$, and fix a finite horizon $T> 0$. Recall the set $\FeasibleSetLDTruncBlw{T}$ defined by \eqref{eq_defn_LDTrunBlw} with $x=s$. If $\pi\in \FeasibleSetLDTruncBlw{T}$, then $\pi$ must satisfy $\LDCostTVec{\pi}{C_{1}}{T}(s)+\LDCostTVec{\pi}{C_{2}}{T}(s)\leq b_{1}+b_{2}-2K\beta^{T}=2-(\frac{1}{2})^{T-2}$, which contradicts \eqref{eq_ctrex2}. The contradiction shows that the feasible region $\FeasibleSetTrunBlw{T}=\FeasibleSetLDTruncBlw{T}$ of the problem $\CRSMDPFinHorOne$ is empty for each $T\in\{1,2,\ldots\}$. 

We claim that the feasible region $\FeasibleSetOrgnal$ of the problem \eqref{eq_CRSMDP_Prob} is nonempty. Indeed, let $\phi\in\Pi_{\rm MR}$ denote the stationary policy that selects from the two actions uniformly at random at each decision epoch. It is easy to check that $\LDCostVec{\phi}{C_{1}}(s)=\LDCostVec{\phi}{C_{2}}(s)=1$. It follows that $\phi\in \FeasibleSetOrgnal=\FeasibleSetLD$, where $\FeasibleSetLD$ is defined by \eqref{eq_Infin_Hor_Lin_Const} with $x=s$. 

We conclude from the previous paragraph that the assertion of Theorem \ref{thm_finhorapprx_frmblw} does not hold for the MDP considered in this example. Next we claim that the local minimum condition assumed in the theorem fails to hold. Indeed, let $h:\Pi_{\rm MR}\rightarrow \R$ denote the maximum constraint violation function defined in the theorem, that is, $h(\pi)=\max\{\LDCostVec{\pi}{C_{1}}(s)-1,\LDCostVec{\pi}{C_{2}}(s)-1\}$ for all $\pi\in \Pi_{\rm MR}$. In light of the  second equality in  \eqref{eq_ctrex2}, we have  $h(\pi)=\max\{\LDCostVec{\pi}{C_{1}}(s)-1,1-\LDCostVec{\pi}{C_{1}}(s)\}\geq 0$ for all $\pi\in\Pi_{\rm MR}$. Also, from the previous paragraph, we have $h(\phi)=0$. Thus $0$ is a global (and hence a local) minimum value for $h$.  This example shows that the assertion of Theorem \ref{thm_finhorapprx_frmblw} may not hold if the local minimum condition assumed in the theorem is violated. 
\end{ex}

\end{document}